%% file: polycat.tex
\documentclass{amsart}
\input{aux/style}
\input{aux/packages} 
\input{aux/usualcmds}

\input{aux/commands} 
\title{Framed polytopes and higher categories}

\author[Laplante-Anfossi]{Guillaume Laplante-Anfossi}
\address{Centre for Quantum Mathematics, Syddansk Universitet, Campusvej 55, Odense, Denmark}
\email{\href{mailto:glaplanteanfossi@imada.sdu.dk}{glaplanteanfossi@imada.sdu.dk}}

\author[Medina-Mardones]{Anibal M. Medina-Mardones}
\address{Department of Mathematics, Western University, ON, Canada.}
\email{\href{mailto:anibal.medina.mardones@uwo.ca}{anibal.medina.mardones@uwo.ca}}

\author[Padrol]{Arnau Padrol}
\address{Departament de Matem\`atiques i Inform\`atica, Universitat de Barcelona, and Centre de Recerca Matem\`atica, Barcelona, Spain.}
\email{\href{mailto:arnau.padrol@ub.edu}{arnau.padrol@ub.edu}}

\date{\today}


\begin{document}
	\begin{center}
		\small{---\textsc{PREPRINT}---}
		\vspace*{15pt}
	\end{center}
	\input{sec/abstract}

	\maketitle
	\tableofcontents
	\input{sec/introduction}
	\input{sec/acknowledgment}

	\input{sec/preliminaries}
	\input{sec/framed}
	\input{sec/strings}
	\input{sec/conjecture}
	\input{sec/orientals}
	\input{sec/matroids}
	\input{sec/modulispaces}
	\input{sec/shelling}
	\sloppy
	\printbibliography
\end{document}

Conventions:
Refer to \section as "part", \subsection as "section" and \subsubsection as "subsection".

\(\set{e_1,\dots,e_d}\) is the canonical basis, not the standard basis.

%% file: aux/style.tex
\makeatletter
\renewcommand\subsubsection{\@startsection{subsubsection}{3}%
	\z@{.5\linespacing\@plus.7\linespacing}{0em}%
	{\normalfont\itshape}}
\makeatother

\usepackage[
	backend = biber,
	style = alphabetic,
	backref = true,
	url = false,
	doi = false,
	isbn = false,
	eprint = false]{biblatex}

\renewbibmacro{in:}{} 
\AtEveryBibitem{
	\clearfield{pages}
	\clearfield{publisher}
	\clearfield{location}
}

\DeclareFieldFormat{title}{\myhref{\mkbibemph{#1}}}
\DeclareFieldFormat
[article,inbook,incollection,inproceedings,patent,thesis,unpublished,book]
{title}{\myhref{\mkbibquote{#1\isdot}}}

\newcommand{\doiorurl}{%
	\iffieldundef{url}
	{\iffieldundef{doi}
		{\iffieldundef{eprint}
			{}
			{https://arxiv.org/abs/\strfield{eprint}}}
		{https://doi.org/\strfield{doi}}}
	{\strfield{url}}
}

\newcommand{\myhref}[1]{%
	\ifboolexpr{%
		test {\ifhyperref}
		and
		not test {\iftoggle{bbx:eprint}}
		and
		not test {\iftoggle{bbx:url}}
	}
	{\href{\doiorurl}{#1}}%
	{#1}%
}

\usepackage[
	bookmarks = true,
	linktocpage = true,
	bookmarksnumbered = true,
	breaklinks = true,
	pdfstartview = FitH,
	hyperfigures = false,
	plainpages = false,
	naturalnames = true,
	colorlinks = true,
	pagebackref = false,
	pdfpagelabels]{hyperref}

\hypersetup{
	colorlinks,
	citecolor = blue,
	filecolor = blue,
	linkcolor = blue,
	urlcolor = blue
}

\usepackage[capitalize, noabbrev]{cleveref}
\crefname{subsection}{\S\!\!}{subsections}

\calclayout

\makeatletter
\@namedef{subjclassname@2020}{%
	\textup{2020} Mathematics Subject Classification}
\makeatother

%% file: aux/packages.tex
\usepackage[T1]{fontenc}
\usepackage{lmodern}
\usepackage[dvipsnames]{xcolor} 
\usepackage{microtype}
\usepackage{amssymb}
\usepackage{mathtools}
\usepackage{mathbbol} 
\usepackage{todo}

\usepackage{csquotes}
\usepackage{xargs} 
\usepackage{comment}
\usepackage{tikz,tikz-cd,tikz-3dplot}
\usetikzlibrary{shapes,positioning,calc,arrows}
\usepackage{enumerate,enumitem}

%% file: aux/usualcmds.tex
\newtheorem{theorem}[subsubsection]{Theorem}
\newtheorem*{theorem*}{Theorem}
\newtheorem{proposition}[subsubsection]{Proposition}
\newtheorem*{proposition*}{Proposition}
\newtheorem{lemma}[subsubsection]{Lemma}
\newtheorem*{lemma*}{Lemma}
\newtheorem{corollary}[subsubsection]{Corollary}
\newtheorem*{corollary*}{Corollary}

\theoremstyle{definition}

\newtheorem*{definition*}{Definition}
\newtheorem{remark}[subsubsection]{Remark}
\newtheorem*{remark*}{Remark}
\newtheorem{example}[subsubsection]{Example}
\newtheorem*{example*}{Example}
\newtheorem{construction}[subsubsection]{Construction}
\newtheorem*{construction*}{Construction}

\newtheorem*{convention*}{Convention}

\newtheorem*{terminology*}{Terminology}

\newtheorem*{notation*}{Notation}

\newtheorem*{question*}{Question}

\DeclareMathOperator{\bd}{\partial}
\DeclareMathOperator{\sign}{sign}
\newcommand{\ot}{\otimes}

\newcommand{\N}{\mathbb{N}}
\newcommand{\Z}{\mathbb{Z}}

\newcommand{\R}{\mathbb{R}}

\newcommand{\sym}{\mathbb{S}}

\newcommand{\gsimplex}{\mathbb{\Delta}}

\newcommand{\Cat}{\mathsf{Cat}}

\newcommand{\Set}{\mathsf{Set}}

\newcommand{\Ch}{\mathsf{Ch}}

\newcommand{\cube}{\square}

\DeclareMathOperator{\chains}{N}

\DeclareMathOperator{\gchains}{C}

\DeclarePairedDelimiter\bars{\lvert}{\rvert}

\DeclarePairedDelimiter\angles{\langle}{\rangle}
\DeclarePairedDelimiter\set{\{}{\}}

\newcommand{\id}{\mathsf{id}}
\renewcommand{\th}{\mathrm{th}}

\newcommand{\Hom}{\mathrm{Hom}}

\newcommand{\xla}[1]{\xleftarrow{#1}}
\newcommand{\xra}[1]{\xrightarrow{#1}}
\newcommand{\defeq}{\coloneqq}

\newcommand{\cC}{\mathcal{C}}

\newcommand{\cG}{\mathcal{G}}
\newcommand{\cH}{\mathcal{H}}

\newcommand{\cN}{\mathcal{N}}

\newcommand{\cP}{\mathcal{P}}

\newcommand{\cS}{\mathcal{S}}
\newcommand{\cT}{\mathcal{T}}
\newcommand{\cU}{\mathcal{U}}

\newcommand{\bG}{\mathbb{G}}


%% file: aux/commands.tex
\newcommand{\resp}{resp.\ }

\definecolor{darkblue}{rgb}{0,0,0.7} 
\newcommand{\darkblue}{\color{darkblue}} 
\newcommand{\defn}[1]{{\darkblue \emph{#1}}}

\renewcommand{\corollary}{\noindent\textit{\darkblue Corollary}.\ }
\renewcommand{\theorem}{\noindent\textit{\darkblue Theorem}.\ }
\renewcommand{\lemma}{\noindent\textit{\darkblue Lemma}.\ }
\renewcommand{\proposition}{\noindent\textit{\darkblue Proposition}.\ }
\renewcommand{\remark}{\noindent\textit{\darkblue Remark}.\ }
\renewcommand{\example}{\noindent\textit{\darkblue Example}.\ }

\newcommand{\conjecture}{\noindent\textit{\darkblue Conjecture}.\ }

\renewcommand{\SS}{\subsubsection{}}

\newcommand{\bach}{\Ch_{\mathsf{ba}}}
\newcommand{\interval}{\mathbb{I}}
\DeclareMathOperator{\normal}{n}
\newcommand{\CP}{C}
\newcommand{\NC}{\cN}

\renewcommand{\bd}{\partial}

\DeclareMathOperator{\cone}{Cone}

\newcommand{\wCat}{\omega\Cat}
\newcommand{\wDiag}{\omega\mathsf{Diag}}
\newcommand{\cells}{\nu}

\DeclareMathOperator{\conv}{conv} 

\DeclareMathOperator{\Mol}{Mol}

\DeclareMathOperator{\so}{s} 
\DeclareMathOperator{\ta}{t} 
\DeclareMathOperator{\bso}{\overline{s}} 
\DeclareMathOperator{\bta}{\overline{t}} 

\renewcommand{\chains}{\gchains}

\newcommand{\sprod}[2]{\left\langle #1 \, , \, #2 \right\rangle} 

\newcommand\restr[2]{{
 \left.\kern-\nulldelimiterspace 
 #1 
 \vphantom{\big|} 
 \right|_{#2} 
 }}

\newcommand{\bigslant}[2]{{\raisebox{.3em}{$#1$} \Big/ \raisebox{-.3em}{$#2$}}} 

\newcommandx{\ssimplex}[1][1 = n]{{\gsimplex_{#1}}} 
\newcommandx{\scube}[1][1 = d]{{\cube_{#1}}} 

\newcommandx{\abs}[1]{{|#1|}} 

\DeclareMathOperator{\GL}{GL}
\newcommandx{\GLn}[1][1 = {\R^d}]{\GL(#1)} 

\DeclareMathOperator{\Aff}{Aff}
\DeclareMathOperator{\Lin}{Lin}
\newcommandx{\Affn}[1][1 = {\R^d}]{\Aff(#1)} 

\newcommandx{\frames}[1][1 = P]{\mathcal{F}_{#1}} 

\newcommandx{\faces}[2][1 = P,2 = {}]{\mathcal{L}_{#2}(#1)} 
\newcommandx{\nefaces}[2][1 = P,2 = {}]{\hat{\mathcal{L}}_{#2}(#1)} 

\newcommandx{\ms}[1][1 = B]{\mathcal{F}({#1})} 
\newcommandx{\msp}[2][1 = B,2 = P]{\mathcal{F}_{{#2}}({#1})} 

\newcommandx{\rs}[1][1 = A]{\mathcal{R}({#1})} 
\newcommandx{\rsinf}[1][1 = {\chi}]{\mathcal{R}_\infty({#1})} 

\newcommandx{\chiro}[1][1 = A]{\chi^{#1}} 
\newcommandx{\cyclift}[1][1 = A]{{#1}^{\uparrow}} 

\newcommand{\ball}[2]{B_{#2}(#1)} 

\DeclareMathOperator{\Ob}{Ob}
\DeclareMathOperator{\Mor}{Mor}

\newcommand{\thickmidarrow}{
\begin{tikzpicture}
	\draw[ultra thick] (-2,1) -- (0,0) ;
	\draw[ultra thick] (-2,-1) -- (0,0) ;
\end{tikzpicture}}

%% file: sec/abstract.tex

\begin{abstract}
	In the early 1990s, Kapranov and Voevodsky proposed a geometric method for constructing higher-categorical pasting diagrams from generically framed convex polytopes.
	This work revisits their construction and identifies a convex-geometric condition that is both necessary and sufficient for the procedure to yield a well-defined pasting diagram.
	Our criterion, the absence of cellular loops, relates their construction to the theory of cellular strings, an active area of convex geometry originating in the Baues problem.
	This paper introduces higher-dimensional cellular strings and uses them to disprove the Kapranov--Voevodsky conjecture in the following strong sense.
	Not only do we exhibit framed polytopes admitting cellular loops, but we also construct examples for which every admissible frame produces one.
	As observed by these authors, Street’s orientals arise from canonically framed cyclic simplices.
	We establish that this family is exceptional as any random \(n\)-simplex, canonically framed, almost surely exhibits cellular loops in the large \(n\)-limit.
\end{abstract}

%% file: sec/introduction.tex

\section{Introduction}

Certain polytopes, such as simplices, cubes, and associahedra, hold a central position in higher category theory, providing a geometric framework to structure coherence relations.
The connection between convex geometry and higher category theory was deepened by Kapranov and Voevodsky in \cite{kapranov1991polycategory}, where they introduced a procedure that, starting from an \(n\)-dimensional polytope with a generic ambient frame, purports to produce a special type of \(n\)-category, a so-called \textit{pasting diagram}.

The concept of pasting diagram is a central yet subtle notion that has received multiple formalizations, including those in \cite{johnson1986pasting, street1991paritycomplexes, Steiner93, steiner2004omega, forest2022pasting, Hadzihasanovic24}.
In this work, we provide a new conceptual description of Steiner's pasting formalism \cite{steiner2004omega} and use it to establish a purely convex-geometric criterion under which the above \textit{Kapranov--Voevodsky Conjecture} holds.

\medskip\theorem The oriented faces of a framed polytope \(P\) define a Steiner diagram if and only if \(P\) has no cellular loops.

\medskip The cellular loops referenced in this theorem are special cases of \emph{cellular strings}, a concept that generalizes to higher dimensions the homonymous notion introduced by Billera, Kapranov, and Sturmfels in \cite{billera1994strings}.
Originating in Baues’s geometric approach to Adams’s cobar construction \cite{adams1956cobar, bauesGeometryLoopSpaces1980}, their (0-dimensional) cellular strings have subsequently become a focus of considerable attention in convex geometry \cite{AthanasiadisEdelmanReiner2000, BilleraSturmfels1992, Bjorner1992, RambauZiegler1996, Reiner1999}.
In this work, we propose the study of cellular strings in all dimensions and use them to provide model-independent counterexamples to the Kapranov--Voevodsky Conjecture.

\medskip\theorem From dimension 4 onward there exist framed polytopes with cellular loops.

\medskip\noindent This result is treated in \cref{s:examples-loops}.
Additionally, in \cref{s:loop-inevitability}, we constructively disprove a weaker version of the conjecture with the following.

\medskip\theorem From dimension 4 onward there exist polytopes for which any frame produces a cellular loop.

\medskip In Theorem~2.5 of the same reference, Kapranov and Voevodsky claimed that their procedure could be used to recover the free \(n\)-category generated by the oriented faces of the \(n\)-simplex—a concept introduced by Street in \cite{street1987orientals}, which is fundamental for all applications of simplicial methods in higher category theory; specifically, when applied to the \textit{cyclic realization} of the \(n\)-simplex, i.e., the convex hull of \(n+1\) distinct points on the \textit{moment curve} \(t \mapsto (t, t^2, \dots, t^n)\), together with the canonical frame \(\set{e_1, \dots, e_n}\).
In \cref{s:street_orientals}, we verify this conjecture after replacing the canonical frame with \(\set{e_1, -e_2, e_3, -e_4, \dots}\).

We also establish that this infinite family of embedded simplices is exceptional in the following probabilistic sense.
A \textit{Gaussian \(n\)-simplex} is the convex hull of \(n+1\) independent random points in \(\R^n\), each selected according to an \(n\)-dimensional standard normal distribution.
In \cref{s:random}, we prove the following.

\medskip\theorem The probability that a canonically framed Gaussian \(n\)-simplex has a cellular loop tends to 1 as \(n\) approaches infinity.

\medskip Regarding cubes, we also prove that the cubical orientals \cite{ aitchison2010cubes} can be obtained by following the Kapranov--Voevodsky procedure using their cyclic embedding and the canonical Euclidean frame, a statement conjectured as Theorem~2.8 in the same reference \cite{kapranov1991polycategory}.
Additionally, we establish an analogous probabilistic result that underscores the exceptional nature of cyclic cubes.

\medskip The passage from a framed polytope to an orientation on each of its faces (\textit{\(f\)-orientation}), which provides the data underlying the connection with higher category theory, raises an interesting moduli question.
Let the \textit{moduli space of frames of an \(f\)-oriented polytope} \(P\) be the set of all frames inducing the given orientation on each face of \(P\).
Using a celebrated result of N. Mn\"ev \cite{Mnev1988}, in \cref{part:universality} we prove the following.

\medskip\theorem For every open primary basic semi-algebraic set \(S\) defined over \(\mathbb{Z}\), there exists an \(f\)-oriented simplex such that its moduli space of frames is stably equivalent to~\(S\).

\medskip We highlight a key step in the proof of this theorem which is interesting in its own right, as it expresses \(f\)-orientations on simplices using matroid-theoretic terms, linking them closely to the emerging concept of oriented flag matroids (cf. \cite{JarraLorscheid2022, BoretskyEurWilliams2022}).
In \cref{part:matroids}, we prove the following.

\medskip\theorem \(f\)-orientations of framed simplices are in bijection with uniform, acyclic, realizable full flag chirotopes.

\medskip Relaxing the acyclicity condition may result in empty sources or targets, whereas relaxing the uniformity condition---or the \( P \)-admissibility condition in the polytopal case---allows faces of dimension greater than \( n \) to serve as \( n \)-cells.
A fully developed theory connecting full flag chirotopes and higher categories is left for future work.

\medskip In \cref{part:molecules}, we study a weakening of the \( P \)-admissibility condition, which leads to the construction of \emph{layerings}---specific decompositions of faces that are closely related to Bruggesser--Mani’s line shellings~\cite{BruggesserMani1971}.
These layerings provide a link between framed polytopes and another model for pasting diagrams, originally introduced by Steiner~\cite{Steiner93} and further developed by Hadzihasanovic~\cite{Hadzihasanovic20, Hadzihasanovic24}.
More precisely, we conclude with the following.

\medskip\theorem Any framed polytope canonically defines a regular directed complex whose atoms are its faces.

%% file: sec/acknowledgment.tex

\subsection*{Acknowledgements}

The authors would like to thank Dimitri Ara, Cl\'emence Chanavat, Tobias Dyckerhoff, Simon Forest, Amar Hadzihasanovic, Mikhail Kapranov, Diana Kessler, Muriel Livernet, Gabriel Longpr\'e, Viktoriya Ozornova, Eva Philippe, Marcy Robertson, Hugh Thomas and Nicholas Williams for insightful discussions related to this work.

We would like to thank the Max Planck Institute of Mathematics in Bonn, the Center for Geometry and Topology in Copenhagen, and the Department of Software Science at TalTech where part of this work was carried out.

\subsection*{Funding}

Work of G. Laplante-Anfossi was supported by the Andrew Sisson Fund and the Australian Research Council Future Fellowship FT210100256, the Novo Nordisk Foundation grant NNF20OC0066298, the Villum Fonden Villum Investigator grant 37814, and the Danish National Research Foundation grant DNRF157.

Work of A. M. Medina-Mardones is partially supported by NSERC grants RES000678 and R7444A03.

Work of A. Padrol is supported by grants PID2022-137283NB-C21 and PCI2024-155081-2 funded by MCIN/AEI/10.13039/501100011033/UE, PAGCAP ANR-21-CE48-0020 of the French National Research Agency ANR, SGR GiT-UB (2021 SGR 00697) funded by the Dept. Recerca i Universitats of Generalitat de Catalunya, and the Severo Ochoa and María de Maeztu Program CEX2020-001084-M of the Spanish State Research Agency.

%% file: sec/preliminaries.tex

\subsection*{Preliminaries}\label{s:preliminaries}

We denote the field of real numbers by~\(\R\), the ring of integers by~\(\Z\), and the monoid of non-negative integers by~\(\N\).

All vector spaces are considered over~\(\R\) and assumed finite-dimensional.
The dual \(\Hom(V, \R)\) of a vector space~\(V\) is denoted by \(V^*\), and the dual (or transpose) of a linear map \(L\) by \(L^*\).
If \(\psi \in V^*\) and \(v \in V\), we write their canonical evaluation \(\psi(v) \in \R\) as \(\sprod{\psi}{v}\), so that for any linear map \(L\) we have
\[
\sprod{L^*(\psi)}{v} = \sprod{\psi}{L(v)}.
\]

%% file: sec/framed.tex

\section{Framed polytopes}

After reviewing foundational concepts in the theory of convex polytopes in \cref{s:polytopes}, we introduce the notion of \textit{\(f\)-orientation} on a polytope \(P\) in \cref{s:orientations}, that is, an orientation on each face of \(P\); we then explore a related combinatorial structure arising from the splitting of the facets of every face into \textit{sources} and \textit{targets}.

In \cref{s:framed_polytopes}, we define a \textit{framed polytope} to be a polytope \(P \subset \R^d\) with a generic frame of \(\R^d\).
Framed polytopes are naturally equipped with an \(f\)-orientation, as well as with \textit{extended} sources and targets.
These can be computed convex-geometrically, as detailed in \cref{s:sk_tk}, and are closely related to \textit{tight coherent subdivisions}, discussed in~\cref{s:subdivision}.

\input{sec/polytopes}
\input{sec/orientations}
\input{sec/frames}
\input{sec/sk_tk}
\input{sec/subdivisions}

%% file: sec/polytopes.tex

\subsection{Polytopes}\label{s:polytopes}

\SS The \defn{convex hull} \(\conv(S)\) of a set \(S \subseteq \R^d\) is the intersection of all convex sets containing \(S\).
A \defn{polytope} \(P\) is the convex hull of a finite set of points in some \(\R^d\).
The notation \(P\) is reserved throughout this article for an unspecified polytope.
Two polytopes are \defn{affinely isomorphic} if they are related by an affine isomorphism.

A \defn{\(d\)-simplex} is the convex hull of \(d+1\) affinely independent points.
All \(d\)-simplices are affinely isomorphic.
The \defn{standard \(d\)-simplex} is \(\ssimplex[d] \defeq \conv(e_1,\ldots,e_{d+1}) \subset \R^{d+1}\), the convex hull of the vectors of the canonical basis.

\SS The \defn{affine span} (resp.\ \defn{linear span}) of \(P\) is the smallest affine (resp.\ linear) subspace \(\Aff P\) (resp.\ \(\Lin P\)) containing \(P\) (resp.\ \(\set{x - y \mid x, y \in P}\)).
We say that vectors in \(\Lin P\) are \defn{parallel} to \(P\).
The \defn{dimension} \(\dim P\) of a polytope \(P \subset \R^d\) is the dimension of its affine span \(\Aff P\).
We say that \(P\) is a \defn{\(d\)-polytope} if \(\dim P = d\).

\SS A \defn{face} of a polytope \(P\) is defined as the \textit{maximal locus} of a (linear) functional on \(P\).
Explicitly, \(F \subseteq P\) is a face if there is \(\phi \in (\R^d)^*\), such that
\[
F = P^\phi \defeq \set[\big]{x \in P \mid \phi(x) = \max_{y \in P} \phi(y)}.
\]
Faces of polytopes are themselves polytopes.
A structure on \(P\) is said to be \defn{inheritable} if it restricts to a structure of the same type on every face of \(P\).
A \defn{\(k\)-face} is a face of dimension \(k\), and its \defn{codimension} is \(\dim P - k\).
The empty set is also considered a face; its dimension is set to \(-1\).

\SS We denote the set of \(k\)-faces of \(P\) by \(\faces[P][k]\) and the set of all faces (resp. non-empty faces) of \(P\) by \(\faces\) (resp. \(\nefaces\)).
Ordered by inclusion, the poset \(\faces\) forms a lattice, the \defn{face lattice} of \(P\).
Two polytopes are \defn{combinatorially equivalent} if their face lattices are isomorphic (as posets).
Straightforwardly, affinely isomorphic polytopes are combinatorially equivalent.

\SS The \defn{normal cone} of a face \(F\) of \(P \subseteq \R^d\) is
\[
\NC^P_F \defeq \set{\psi \in (\R^d)^* \mid F = P^\psi}.
\]
We refer to elements in \(\NC^P_F\) as \defn{normal covectors}, reserving the notation \(\normal_F^P\) for an arbitrary one.

Note that we do not require \(P\) to be \(d\)-dimensional to define normal covectors. If~\(P\) is not \(d\)-dimensional, then its normal cone contains the annihilator of \(\Lin P\).

For a facet~\(F\) of \(P\), if \(P\) is \(d\)-dimensional then \(\normal_F^P\) is uniquely defined up to multiplication by a positive scalar. In general, for a facet~\(F\), \(\normal_F^P\) is uniquely defined up to multiplication by a positive scalar and addition of vectors in the annihilator of \(\Lin P\).

\SS\label{ss:inner-pointing vectors}
Faces of codimension \(1\) are called \defn{facets}.
The set of facets of a polytope \(P\) will be denoted \(\bd(P)\) instead of \(\faces[P][\dim(P)-1]\), and referred to as its \defn{boundary}.
We remark that
\[
\bigcup_{\ \mathclap{F \in \bd(P)}} \,F
\]
is the \textit{topological boundary} of \(P\).

Given \(F \in \bd(P)\), a vector \(v \in \Lin(P)\) is said to be \defn{inner-pointing} (resp.\ \defn{outer-pointing}) if
\[
\langle \normal, v \rangle < 0 \quad \text{(resp. } > 0 \text{)}
\]
for any \(\normal \in \NC^P_F\).
This definition is independent of the choice of normal covector \(\normal\), because any \(v \in \Lin(P)\) can be written as
\[
v = \mu(q - p),
\]
for some \(p \in P\), \(q \in F\), and \(\mu \in \mathbb{R}\), and \(\sprod{\normal}{q} \geq \sprod{\normal}{p}\) for any \(\normal \in \NC^P_F\).

\SS The \defn{Minkowski sum} of two subsets \(A, B\) of \(\R^d\) is defined as
\[
A + B = \set{a + b \mid a \in A,\ b \in B}.
\]
If \(P\) and \(Q\) are polytopes, then their Minkowski sum \(P + Q\) is also a polytope.
A Minkowski sum of segments is known as a \defn{zonotope}.
Every zonotope is the image of the \defn{unit cube} \(\scube \defeq [0,1]^d\) under an affine map.

%% file: sec/orientations.tex

\subsection{Orientations}\label{s:orientations}

\SS\label{ss:facial-orientations}
The \defn{exterior algebra} \(\bigwedge V\) of a vector space \(V\) is the quotient of the tensor algebra
\[
\R \oplus V \oplus V^{\ot 2} \oplus \dotsb
\]
by the two-sided ideal generated by all elements \(v \ot v\) with \(v \in V\).
We use the notation \(v_1 \wedge\dots\wedge v_k\) and \(\bigwedge^k V\) respectively for the images of \(v_1 \ot\dotsb\ot v_k\) and \(V^{\ot k}\) in this quotient.
Note that, for \(k = 0\), \(\bigwedge^0 V\cong \R\).
If \(V\) is \(n\)-dimensional, then \(\bigwedge^n V\) is one-dimensional.
We say that two generators \(\beta, \beta' \in \bigwedge^n V\) are \defn{equivalent}, denoted \(\beta \sim \beta'\), if there exists \(\lambda > 0\) such that \(\beta = \lambda \, \beta'\).
An \defn{orientation} is a choice of one of the two resulting equivalence classes.

\SS A (positive) \defn{\(f\)-orientation} on a polytope \(P\), short for full orientation, is a choice of orientation for the linear span of each face of \(P\), such that the orientation assigned to each vertex is the positive generator \(1 \in \R\).

An \(f\)-orientation on a polytope \(P\) restricts to each of its faces, and we refer to the resulting orientations as the \defn{inherited \(f\)-orientations}.

\SS\label{d:sources_and_targets}
The \(f\)-orientation and the notions of inner- and outer-pointing vectors determine a partition of the facets of \(P\) into \textit{sources} and \textit{targets}.
Intuitively, a facet is in the source (\resp target) of \(P\) if its orientation can be completed to that of \(P\) by wedging, on the right, with an inner-pointing (\resp outer-pointing) vector.

More precisely, let \(\beta\) be an \(f\)-orientation on \(P\).
The \defn{source} \(\so(P)\) (\resp \defn{target} \(\ta(P)\)) of \(P\) contains a facet \(F \in \bd(P)\) if and only if there exist \(v \in \Lin(P)\) and \(\normal \in \NC^P_F\) with \(\langle \normal, v \rangle < 0\) (\resp \(> 0\)) such that \(\beta_P \sim \beta_F \wedge v\).

The source \(\so(F)\) and target \(\ta(F)\) of a face \(F\) of \(P\) are defined using the inherited \(f\)-orientation on \(F\).

For completeness, we verify that these definitions are independent of the choice of vector \(v \in \Lin(P)\) and covector \(\normal \in \NC^P_F\).
Let \(v' \in \Lin(P)\) be another vector such that \(\beta_P \sim \beta_F \wedge v'\).
This condition is equivalent to the existence of \(\lambda > 0\) and \(u \in \Lin(F)\) such that
\[
v' = \lambda v + u.
\]
Any \(u \in \Lin(F)\) can be written as \(u = \mu(p - q)\) for some \(\mu \in \R\) and \(p, q \in F\).
Since \(\normal\) attains the same (maximal) value on \(p\) and \(q\) we have
\[
\sprod{\normal}{v'} = \lambda \cdot \sprod{\normal}{v}\!,
\]
which shows that the sign of the pairing is preserved as claimed.

Independence of \(\normal\) follows by the same argument used in the definition of inner- and outer-pointing vectors in \cref{ss:inner-pointing vectors}.

\SS\remark\label{ss:sotagiveforientation}
A (positive) \(f\)-orientation on a polytope \(P\) is completely determined by the partitions \(\bd(F) = \so(F) \sqcup \ta(F)\) on all faces \(F\) of \(P\).

To see this, consider two \(f\)-orientations \(\beta\) and \(\gamma\) which, by definition, agree on the vertices of \(P\).
Assume inductively that \(\beta_E \sim \gamma_E\) for every face \(E\) of dimension at most \(n - 1\) and let \(F\) be a face of dimension \(n\).
We aim to show that \(\beta_F \sim \gamma_F\).

Let \(E\) be a facet of \(F\), which we may assume without loss of generality lies in \(\ta(F)\).
By the inductive hypothesis, we have \(\beta_E \sim \gamma_E\).
Choosing any outward-pointing vector \(v \in \Lin(F)\), we compute:
\[
\beta_F \sim \beta_E \wedge v \sim \gamma_E \wedge v \sim \gamma_F
\]
as claimed.

%% file: sec/frames.tex

\subsection{Frames}\label{s:framed_polytopes}

\SS \label{sec:frame} A \defn{frame} is an ordered basis \((v_1,\dots,v_d)\) of \(\R^d\) for some \(d \in \N\).
The notation \(v_i\) is reserved throughout this article for elements in a given frame.

The \defn{canonical frame} \((e_1,\dots,e_d)\) of \(\R^d\) consists of the canonical basis with its natural order.

A frame of \(\R^d\) defines a natural collection of subspaces \(V_k \defeq \Lin{(v_1,\ldots,v_k)}\) and a \defn{system of projections}, defined as the collection of linear maps
\[
\set{\pi_k \colon \R^d \to V_k}_{k=0}^d
\]
such that
\[
\pi_k(v_i) =
\begin{cases}
	v_i & \text{ if } i \leq k,\\
	\hfil 0 & \text{ if } i > k.
\end{cases}
\]
The notation \(\pi_k\) is reserved for the \(k^\th\) map in the system of projections of a given frame.

\SS A frame is said to be \defn{\(P\)-admissible} if for any \(k\)-face \(F\) of \(P\) the restriction \(\pi_k \colon \Lin F \to V_k\) is a linear isomorphism.
In this case, we denote its inverse by \(\sigma_k^F \colon V_k \to \Lin F\).
We reserve the notation \(\sigma_k^F\) for such inverses throughout the text.

We remark that \(P\)-admissibility is an open and inheritable property.
That is, any sufficiently small perturbation of a \(P\)-admissible frame remains \(P\)-admissible,
and any \(P\)-admissible frame is also \(F\)-admissible for every face \(F\) of \(P\).

\SS A \defn{framed polytope} is a pair consisting of a polytope \(P\) and a \(P\)-admissible frame~\(B\).
In this case we say that \(P\) is framed by \(B\), although we typically omit references to the frame.

The inheritability of \(P\)-admissibility implies that the faces of a framed polytope are themselves framed.

\SS\label{induced_facial_orientation}
The \defn{induced \(f\)-orientation} of a polytope framed by \((v_1,\dots,v_d)\) is given on a \(k\)-face \(F\) by
\[
\beta_F \defeq
\begin{cases}
	\sigma^F_k(v_1 \wedge\dots\wedge v_k) & k > 0, \\
	\hfil 1 & k = 0,
\end{cases}
\]
where, as usual, \(\sigma^F_k(v_1 \wedge\dots\wedge v_k) = \sigma^F_k(v_1) \wedge\dots\wedge \sigma^F_k(v_k)\).

\SS\lemma\label{lem:sources and targets in a framed polytope}
In a framed polytope, let \(E\) be a facet of a \(k\)-face \(F\).
Then, the statements in the following columns are equivalent:
\begin{center}
	\renewcommand{\arraystretch}{2} 
	\begin{tabular}{r@{\hspace{1em}}l@{\hspace{5em}}r@{\hspace{1em}}l}
		(1) & \(E \in \so(F)\), &
		(1') & \(E \in \ta(F)\), \\
		(2) & \(\sprod{\normal^F_E}{\sigma_{k}^F(v_{k})} < 0\), &
		(2') & \(\sprod{\normal^F_E}{\sigma_{k}^F(v_{k})} > 0\), \\
		(3) & \(\sprod{\normal^{\pi_{k}(F)}_{\pi_{k}(E)}}{v_{k}} < 0\). &
		(3') & \(\sprod{\normal^{\pi_{k}(F)}_{\pi_{k}(E)}}{v_{k}} > 0\). \\
	\end{tabular}
\end{center}

\begin{proof}
	We will prove the equivalence of the statements in the left column only, as those in the right column follow analogously.

	\medskip\noindent
	(1) \(\Leftrightarrow\) (2):
	This equivalence amounts to showing that
	\[
	\beta_E \wedge \sigma_k^F(v_k) \sim \beta_F,
	\]
	where \(\beta\) denotes the induced \(f\)-orientation.
	We will in fact prove equality.

	Let \(\iota_{k-1}^k \colon V_{k-1} \to V_k\) be the canonical inclusion, which satisfies \(\sigma_{k-1}^E = \sigma_k^F \circ \iota_{k-1}^k\).
	Applying \(\sigma_k^F\) to the identity
	\[
	\iota_{k-1}^k(v_1 \wedge \dots \wedge v_{k-1}) \wedge v_k = v_1 \wedge \dots \wedge v_k
	\]
	yields
	\[
	\sigma_{k-1}^E(v_1 \wedge \dots \wedge v_{k-1}) \wedge \sigma_k^F(v_k)
	=
	\sigma_k^F(v_1 \wedge \dots \wedge v_k),
	\]
	the desired identity.

	\medskip\noindent
	(2) \(\Leftrightarrow\) (3):
	For simplicity, we write \(\pi = \pi_k\) and \(\sigma = \sigma_k^F\).
	Given that the sign of \(\langle \normal, v_{k} \rangle\) is constant for all \(\normal \in \NC^{\pi(P)}_{\pi(F)}\), it suffices to show that \(\sigma^*(\normal_E^F) \in \NC^{\pi(P)}_{\pi(F)}\).
	Since \(\normal_E^F \in \NC_E^F\), we have
	\[
	E = F^{\normal_E^F} = \left\{ p \in F \mid \langle \normal_E^F, p \rangle = \max_F \langle \normal_E^F, - \rangle \right\}.
	\]
	Applying \(\pi\) and using \(\sigma \circ \pi = \mathrm{id}\) on \(F\), we obtain:
	\begin{align*}
		\pi(E) &=
		\set[\Big]{\pi(p) \in \pi(F) \mid \sprod{\normal_E^F}{p} = \max_F \sprod{\normal_E^F}{-}} \\ &=
		\set[\Big]{\pi(p) \in \pi(F) \mid \sprod{\normal_E^F}{(\sigma\circ\pi)(p)} = \max_F \sprod{\normal_E^F}{(\sigma\circ\pi)(-)}} \\ &=
		\set[\Big]{\pi(p) \in \pi(F) \mid \sprod{\sigma^*(\normal_E^F)}{\pi(p)} = \max_{\pi(F)} \sprod{\sigma^*(\normal_E^F)}{-}} \\ &=
		\pi(F)^{\sigma^*(\normal_E^F)},
	\end{align*}
	as desired.
\end{proof}

\begin{figure}[h!]
	\input{fig/lower_triangular}
	\caption{The two frames \((v,w)\) and \((v',w)\) are \(P\)-equivalent for the regular hexagon \(P\).}
	\label{fig:lower-triangular}
\end{figure}
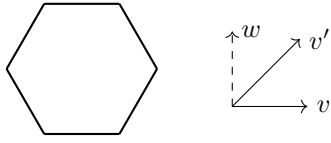

\SS Two \(P\)-admissible frames are said to be \defn{\(P\)-equivalent} if they induce the same \(f\)-orientation on \(P\), please consult \cref{fig:lower-triangular} for an example.

\SS\label{ss:invariancelinearautomorphism}
Transforming a framed polytope and its frame via a linear automorphism preserves its induced \(f\)-orientation.
More precisely, we have the following.

\medskip\lemma
Let \(\phi \colon \R^d \to \R^d\) be a linear automorphism.
If \(P\) is framed by \(\{v_1, \dots, v_d\}\), then \(\phi(P)\) is framed by \(\{\phi(v_1), \dots, \phi(v_d)\}\),
and with respect to the corresponding induced \(f\)-orientations, we have
\[
\beta_{\phi(F)} = \phi(\beta_F)
\]
for any face \(F\) of \(P\).

\begin{proof}
	Let \(\pi_k\) be the \(k^\text{th}\) canonical projection associated with the \(P\)-admissible frame \(B = \{v_1, \dots, v_d\}\).
	Then \(\phi \mathrel{\circ} \pi_k \mathrel{\circ} \phi^{-1}\) is the \(k^\text{th}\) canonical projection for the frame \(\phi(B) = \{\phi(v_1), \dots, \phi(v_d)\}\).
	From this, the \(\phi(P)\)-admissibility of \(\phi(B)\) is straightforward to verify.

	Similarly, for any \(k\)-face \(\phi(F)\), the canonical inverse \(\sigma_k^{\phi(F)}\) agrees with \(\phi \mathrel{\circ} \sigma_k^F \mathrel{\circ} \phi^{-1}\), and we compute:
	\begin{align*}
		\beta_{\phi(F)}
		&= \sigma_k^{\phi(F)}(\phi(v_1) \wedge \dots \wedge \phi(v_k)) \\
		&= (\phi \circ \sigma_k^F \circ \phi^{-1})(\phi(v_1) \wedge \dots \wedge \phi(v_k)) \\
		&= \phi(\sigma_k^F(v_1 \wedge \dots \wedge v_k)) \\
		&= \phi(\beta_F). \qedhere
	\end{align*}
\end{proof}

\SS\label{ss:lowertriangular}

A sufficient condition for two \(P\)-admissible frames to be \(P\)-equivalent is that their change-of-basis matrix is positive and lower triangular.

\medskip\lemma Let \(B = (v_1,\dots,v_d)\) be a \(P\)-admissible frame.
If a frame \(B' = (v_1',\dots,v'_d)\) is obtained from~\(B\) via a positive lower triangular transformation, meaning that there exists \(\lambda_{pq} \in \R\) for \(p > q\) and \(\lambda_i \in \R_+\) such that
\begin{equation}\label{eq:upper_triangular}
	v_q' = \lambda_q v_q + \sum_{p>q} \lambda_{pq} v_p \,,
\end{equation}
then \(B'\) is \(P\)-admissible and \(P\)-equivalent to \(B\).

\begin{proof}
	Consider the systems of projections \(\set{\pi_k}\) and \(\set{\pi_k'}\) corresponding respectively to \(B\) and \(B'\).
	For any face \(F\), say of dimension \(k\), consider the orientations \(\sigma^F_k(v_1 \wedge\dots\wedge v_k)\) and \(\sigma'^F_k(v'_1 \wedge\dots\wedge v'_k)\) as defined in \cref{induced_facial_orientation}.
	To show that these agree, it suffices to verify that \((v_1 \wedge\dots\wedge v_k) \sim (\pi_k \circ \sigma'^F_k)(v'_1 \wedge\dots\wedge v'_k)\).

	Let us start by noticing that \(\ker \pi'_k \subseteq \ker \pi_k\).
	Indeed, \(\ker \pi_k'\) is generated by \(\set{v'_{k+1}, v'_{k+2}, \dots}\) and each of these elements is in \(\ker \pi_k\) by \eqref{eq:upper_triangular}.

	Now, for each \(i \leq k\), we have \(\sigma'^F_k(v'_i) = v'_i + \kappa_i\) where \(\kappa_i \in \ker \pi'_k \subseteq \ker \pi_k\).
	Therefore,
	\[
	(\pi_k \circ \sigma'^F_k)(v'_i) =
	\pi_k(v'_i + \kappa_i) =
	\pi_k(v'_i) =
	\pi_k(\lambda_i v_i + \sum_{j > i} \lambda_{ji} v_j) =
	\lambda_i v_i + \hspace{7pt} \mathclap{ \sum_{j = i+1}^k} \hspace{9pt} \lambda_{ji} v_j.
	\]
	Using that \(v_j \wedge v_j = 0\) for any \(j\) we conclude that
	\[
	(\pi_k \circ \sigma'^F_k)(v'_1 \wedge\dots\wedge v'_k) =
	(\lambda_1 \cdots \lambda_k) \cdot v_1 \wedge\dots\wedge v_k,
	\]
	which is the same orientation as \(v_1 \wedge\dots\wedge v_k\) since each \(\lambda_i > 0\).
\end{proof}

\SS\label{ss:equivalentorthogonal}
\corollary Every \(P\)-admissible frame is \(P\)-equivalent to an orthonormal frame.

\begin{proof}
	This follows directly from \cref{ss:lowertriangular}, as the Gram--Schmidt orthonormalization transformation is positive lower triangular, when applied from last to first element of the basis.
\end{proof}

%% file: fig/lower_triangular.tex
\begin{tikzpicture}
	\draw[->] (2, -.5) -- (3, -.5) node[anchor = west]{$v$};
	\draw[->] (2, -.5) -- (2.9, .4) node[anchor = west]{$v'$};
	\draw[->,dashed] (2, -.5) -- (2, .5) node[anchor = west]{$w$};

	\node (a) at (1, 0) {};

	\node (b) at (0.5, 0.866) {};

	\node (c) at (-0.5, 0.866) {};

	\node (d) at (-1, 0) {};

	\node (e) at (-0.5, -0.866) {};

	\node (f) at (0.5, -0.866) {};

	\draw[thick] (-0.5, 0.866)--(0.5, 0.866);
	\draw[thick] (1, 0)--(0.5, 0.866);
	\draw[thick] (-0.5, 0.866)--(-1, 0);

	\draw[thick] (1,0)--(0.5, -0.866);
	\draw[thick] (-0.5, -0.866)--(0.5, -0.866);
	\draw[thick] (-0.5, -0.866)--(-1, 0) ;
\end{tikzpicture}

%% file: sec/sk_tk.tex

\subsection{Extended sources and targets}\label{s:sk_tk}

We now introduce the extended sources and targets of a framed polytope \(P\).
These are indexed by a non-negative integer \(k < \dim(P)\), and coincide with the usual sources and targets of \(P\) when \(k = \dim(P) - 1\).
In addition, we provide an effective method for computing extended sources and targets.

\SS \label{def:extended-st}
Let \(P\) be a framed polytope and \(k \in \N\) with \(k < \dim(P)\).
The \defn{\(k\)-boundary} of \(P\) is defined as
\[
\bd_k(P) \defeq \set{F \in \faces_k \mid \pi_{k+1}(F) \in \bd(\pi_{k+1}(P))}.
\]
The \defn{\(k\)-source} of \(P\) is defined as
\[
\so_k(P) \defeq \set[\Big]{F \in \bd_k(P) \mid \pi_{k+1}(F) \in \so(\pi_{k+1}(P))}.
\]
By \cref{lem:sources and targets in a framed polytope},
\[
\so_k(P) =
\set[\Big]{F \in \bd_k(P) \mid \sprod{\normal_{\pi_{k+1}(F)}^{\pi_{k+1}(P)}}{v_{k+1}} < 0}.
\]
The \defn{\(k\)-target} of \(P\), denoted \(\ta_k(P)\), is defined analogously.

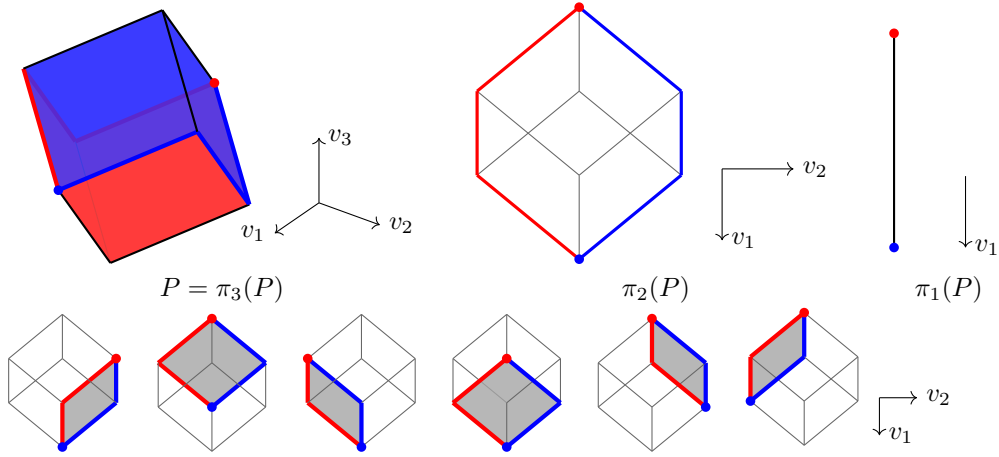
\begin{figure}[htpb]
	\centering
	\input{fig/globular_figure}
	\caption{The first row depicts \(P\) and its projections \(\pi_2(P)\) and \(\pi_1(P)\). The faces in \(s_0(P)\), \(s_1(P)\) and \(s_2(P)\), and their projections, are in red, while the faces in \(t_0(P)\), \(t_1(P)\) and \(t_2(P)\), and their projections, are in blue. The second row shows the \(0\)- and \(1\)-sources and targets of the \(2\)-faces, projected onto the \(\langle v_1, v_2\rangle\) plane. The \(0\)-sources and targets of the \(1\)-faces are computed similarly.}
	\label{fig:globe-poly}
\end{figure}

Extended sources and targets of a face \(F\) of \(P\) are defined by regarding \(F\) as framed polytopes.
In \cref{fig:globe-poly} we presents some examples.

\SS\label{ss:explicitdescription}
\lemma Let \(P \subset \R^d\) be a framed polytope.
Let \(F\) be a \(k\)-face of \(P\) and \(\set{p_1, \dots, p_{k+1}}\) a set of affinely independent vertices of it.
Then, \(F\) is in the \(k\)-boundary of \(P\) if and only if the sign of the determinant
\begin{equation}\label{eq:explicit_determinant}
	\begin{vmatrix}
		1 & \cdots & 1 & 1 & 0 & \cdots & 0 \\
		(p_1)_{1} & \cdots & (p_{k+1})_{1} & q_1 & (v_{k+2})_1 & \cdots & (v_{d})_1 \\
		(p_1)_{2} & \cdots & (p_{k+1})_{2} & q_2 & (v_{k+2})_2 & \cdots & (v_{d})_2 \\
		\vdots & \ddots & \vdots & \vdots & \vdots & \ddots & \vdots \\
		(p_1)_{d} & \cdots & (p_{k+1})_{d} & q_d & (v_{k+2})_d & \cdots & (v_{d})_d \\
	\end{vmatrix}
\end{equation}
is constant for all vertices \(q\) of \(P\) that do not belong to \(F\), where we use $x_i$ to denote the $i$th coordinate with respect to a fixed basis.
In this case, \(F\) will be in the \(k\)-source of \(P\) if the sign of~\eqref{eq:explicit_determinant} coincides with the sign of
\begin{equation}\label{eq:explicit_determinant_2}
	\begin{vmatrix}
		1 & \cdots & 1 & 0 & 0 & \cdots & 0 \\
		(p_1)_{1} & \cdots & (p_{k+1})_{1} & (v_{k+1})_1 & (v_{k+2})_1 & \cdots & (v_{d})_1 \\
		(p_1)_{2} & \cdots & (p_{k+1})_{2} & (v_{k+1})_2 & (v_{k+2})_2 & \cdots & (v_{d})_2 \\
		\vdots & \ddots & \vdots & \vdots & \vdots & \ddots & \vdots \\
		(p_1)_{d} & \cdots & (p_{k+1})_{d} & (v_{k+1})_d & (v_{k+2})_d & \cdots & (v_{d})_d \\
	\end{vmatrix},
\end{equation}
and will be in the \(k\)-target of \(P\) otherwise.

\begin{proof}
	Consider the linear functional \(\psi \in (\R^d)^*\) explicitly defined as
	\[
	\sprod{\psi}{x} =
	\begin{vmatrix}
		1 & \cdots & 1 & 0 & 0 & \cdots & 0 \\
		(p_1)_{1} & \cdots & (p_{k+1})_{1} & x_1 & (v_{k+2})_1 & \cdots & (v_{d})_1 \\
		(p_1)_{2} & \cdots & (p_{k+1})_{2} & x_2 & (v_{k+2})_2 & \cdots & (v_{d})_2 \\
		\vdots & \ddots & \vdots & \vdots & \vdots & \ddots & \vdots \\
		(p_1)_{d} & \cdots & (p_{k+1})_{d} & x_d & (v_{k+2})_d & \cdots & (v_{d})_d \\
	\end{vmatrix}.
	\]

	Since \(\sprod{\psi}{v_j} = 0\) for all \(k+2 \leq j \leq d\), this functional is in the image of \(\pi_{k+1}^*\), say \(\psi = \pi_{k+1}^*(\phi)\) with \(\phi \in (\R^{k+1})^*\).

	Additionally, \(\psi\) is constant on \(F\) with value
	\[
	c = -
	\begin{vmatrix}
		1 & \cdots & 1 & 1 & 0 & \cdots & 0 \\
		(p_1)_{1} & \cdots & (p_{k+1})_{1} & 0 & (v_{k+2})_1 & \cdots & (v_{d})_1 \\
		(p_1)_{2} & \cdots & (p_{k+1})_{2} & 0 & (v_{k+2})_2 & \cdots & (v_{d})_2 \\
		\vdots & \ddots & \vdots & \vdots & \vdots & \ddots & \vdots \\
		(p_1)_{d} & \cdots & (p_{k+1})_{d} & 0 & (v_{k+2})_d & \cdots & (v_{d})_d \\
	\end{vmatrix}.
	\]
	To see this, one needs to simply evaluate it on each of the elements in the linearly independent set \(\set{p_1, \dots, p_{k+1}}\).

	By dimensional reasons, any functional in \((\R^{k+1})^*\) that is constant on \(\pi_{k+1}(F)\) lies in the linear span of \(\phi\), which is not \(0\) by the \(P\)-admissibility of the frame.
	It follows that \(\pi_{k+1}(F) \in \bd(\pi_{k+1}(P))\) if and only if \(\pi_{k+1}(F) = \pi_{k+1}(P)^{\pm \phi}\).
	Explicitly, this occurs when for all vertices \(\pi_{k+1}(q) \in \pi_{k+1}(P) \setminus \pi_{k+1}(F)\), either
	\[
	\sprod{\phi}{\pi_{k+1}(q)} < c \quad \text{or} \quad \sprod{\phi}{\pi_{k+1}(q)} > c.
	\]
	This is equivalent to \(\sprod{\psi}{q} - c\) having the same sign on all vertices \(q \in P \setminus F\), which verifies the first claim, given that \(\sprod{\psi}{q} - c\) equals \eqref{eq:explicit_determinant}.

	The second claim follows directly from \cref{lem:sources and targets in a framed polytope}, since \eqref{eq:explicit_determinant_2} is equal to
	\[
	\sprod{\psi}{v_{k+1}} =
	\sprod{\pi_{k+1}^*(\phi)}{v_{k+1}} =
	\sprod{\phi}{\pi_{k+1}(v_{k+1})} =
	\sprod{\phi}{v_{k+1}}. \qedhere
	\]
\end{proof}

%% file: fig/globular_figure.tex
\begin{tabular}{ccccc}
 \begin{tikzpicture}%
	[x = {(-0.587802cm, -0.404450cm)},
	y = {(0.809005cm, -0.293947cm)},
	z = {(0.000078cm, 0.866034cm)},
	scale = 1.000000,
	back/.style = {thin, color = black!60},
	edge/.style = {color = black, thick},
	sourceedge/.style = {color = red, ultra thick},
	targetedge/.style = {color = blue, ultra thick},
	facet/.style = {fill = blue!95!black,fill opacity = 0.800000},
	targetfacet/.style = {fill = blue!80,fill opacity = 0.800000},
	sourcefacet/.style = {fill = red!80,fill opacity = 0.800000},
	vertex/.style = {inner sep = 1pt,circle,draw = black,fill = black,thick},
	targetvertex/.style = {inner sep = 1pt,circle,draw = blue,fill = blue,thick},
	sourcevertex/.style = {inner sep = 1pt,circle,draw = red,fill = red,thick}]
%
%
\coordinate (-0.58824, 1.42857, -0.83333) at (-0.58824, 1.42857, -0.83333);
\coordinate (0.58824, 1.42857, 0.83333) at (0.58824, 1.42857, 0.83333);
\coordinate (1.76471, 0.00000, 0.00000) at (1.76471, 0.00000, 0.00000);
\coordinate (0.58824, 0.00000, -1.66667) at (0.58824, 0.00000, -1.66667);
\coordinate (-0.58824, -1.42857, -0.83333) at (-0.58824, -1.42857, -0.83333);
\coordinate (-1.76471, 0.00000, 0.00000) at (-1.76471, 0.00000, 0.00000);
\coordinate (-0.58824, 0.00000, 1.66667) at (-0.58824, 0.00000, 1.66667);
\coordinate (0.58824, -1.42857, 0.83333) at (0.58824, -1.42857, 0.83333);

\fill[targetfacet] (-0.58824, 0.00000, 1.66667) -- (0.58824, -1.42857, 0.83333) -- (-0.58824, -1.42857, -0.83333) -- (-1.76471, 0.00000, 0.00000) -- cycle {};
\fill[sourcefacet] (-0.58824, 1.42857, -0.83333) -- (-1.76471, 0.00000, 0.00000) -- (-0.58824, -1.42857, -0.83333) -- (0.58824, 0.00000, -1.66667) -- cycle {};
 \fill[sourcefacet] (0.58824, 0.00000, -1.66667) -- (-0.58824, -1.42857, -0.83333) -- (0.58824, -1.42857, 0.83333) -- (1.76471, 0.00000, 0.00000) -- cycle {};

\draw[edge,back] (0.58824, 0.00000, -1.66667) -- (-0.58824, -1.42857, -0.83333);
\draw[back,sourceedge] (-0.58824, -1.42857, -0.83333) -- (-1.76471, 0.00000, 0.00000);
\draw[back,sourceedge] (-0.58824, -1.42857, -0.83333) -- (0.58824, -1.42857, 0.83333);

\fill[sourcefacet] (0.58824, 0.00000, -1.66667) -- (-0.58824, 1.42857, -0.83333) -- (0.58824, 1.42857, 0.83333) -- (1.76471, 0.00000, 0.00000) -- cycle {};
\fill[targetfacet] (0.58824, -1.42857, 0.83333) -- (1.76471, 0.00000, 0.00000) -- (0.58824, 1.42857, 0.83333) -- (-0.58824, 0.00000, 1.66667) -- cycle {};
 \fill[targetfacet] (-0.58824, 0.00000, 1.66667) -- (0.58824, 1.42857, 0.83333) -- (-0.58824, 1.42857, -0.83333) -- (-1.76471, 0.00000, 0.00000) -- cycle {};
\draw[targetedge] (-0.58824, 1.42857, -0.83333) -- (0.58824, 1.42857, 0.83333);
\draw[edge] (-0.58824, 1.42857, -0.83333) -- (0.58824, 0.00000, -1.66667);
\draw[targetedge] (-0.58824, 1.42857, -0.83333) -- (-1.76471, 0.00000, 0.00000);
\draw[targetedge] (0.58824, 1.42857, 0.83333) -- (1.76471, 0.00000, 0.00000);
\draw[edge] (0.58824, 1.42857, 0.83333) -- (-0.58824, 0.00000, 1.66667);
\draw[edge] (1.76471, 0.00000, 0.00000) -- (0.58824, 0.00000, -1.66667);
\draw[sourceedge] (1.76471, 0.00000, 0.00000) -- (0.58824, -1.42857, 0.83333);
\draw[edge] (-1.76471, 0.00000, 0.00000) -- (-0.58824, 0.00000, 1.66667);
\draw[edge] (-0.58824, 0.00000, 1.66667) -- (0.58824, -1.42857, 0.83333);
\node[targetvertex] at (1.76471, 0.00000, 0.00000) {};
 \node[sourcevertex] at (-1.76471, 0.00000, 0.00000) {};
\begin{scope}[shift = {(0,3,0)}]

		\draw[->] (0, 0,0) -- (1, 0,0) node[anchor = east]{$v_1$};
		\draw[->] (0, 0,0) -- (0, 1,0) node[anchor = west]{$v_2$};
		\draw[->] (0, 0,0) -- (0, 0, 1) node[anchor = west]{$v_3$};
	\end{scope}

\end{tikzpicture}
& &
\begin{tikzpicture}%
	[
	x = {(0.000000cm, -1.000000cm)},
	y = {(1.000000cm, 0.000000cm)},
	z = {(0.000000cm, 0.000000cm)},
	scale = .950000,
	back/.style = {thin, color = black!60},
	edge/.style = {color = black, thick},
	sourceedge/.style = {color = red, very thick},
	targetedge/.style = {color = blue, very thick},
	facet/.style = {fill = blue!95!black,fill opacity = 0.800000},
	targetfacet/.style = {fill = blue,fill opacity = 0.600000},
	sourcefacet/.style = {fill = red,fill opacity = 0.600000},
	vertex/.style = {inner sep = 1pt,circle,draw = black,fill = black,thick},
	targetvertex/.style = {inner sep = 1pt,circle,draw = blue,fill = blue,thick},
	sourcevertex/.style = {inner sep = 1pt,circle,draw = red,fill = red,thick}]
%
%
%
%
\coordinate (-0.58824, 1.42857, -0.83333) at (-0.58824, 1.42857, -0.83333);
\coordinate (0.58824, 1.42857, 0.83333) at (0.58824, 1.42857, 0.83333);
\coordinate (1.76471, 0.00000, 0.00000) at (1.76471, 0.00000, 0.00000);
\coordinate (0.58824, 0.00000, -1.66667) at (0.58824, 0.00000, -1.66667);
\coordinate (-0.58824, -1.42857, -0.83333) at (-0.58824, -1.42857, -0.83333);
\coordinate (-1.76471, 0.00000, 0.00000) at (-1.76471, 0.00000, 0.00000);
\coordinate (-0.58824, 0.00000, 1.66667) at (-0.58824, 0.00000, 1.66667);
\coordinate (0.58824, -1.42857, 0.83333) at (0.58824, -1.42857, 0.83333);
\draw[edge,back] (-0.58824, 1.42857, -0.83333) -- (0.58824, 0.00000, -1.66667);
\draw[edge,back] (1.76471, 0.00000, 0.00000) -- (0.58824, 0.00000, -1.66667);
\draw[edge,back] (0.58824, 0.00000, -1.66667) -- (-0.58824, -1.42857, -0.83333);

\draw[edge,back] (0.58824, 1.42857, 0.83333) -- (-0.58824, 0.00000, 1.66667);
\draw[edge,back] (-1.76471, 0.00000, 0.00000) -- (-0.58824, 0.00000, 1.66667);
\draw[edge,back] (-0.58824, 0.00000, 1.66667) -- (0.58824, -1.42857, 0.83333);


\draw[back,sourceedge] (-0.58824, -1.42857, -0.83333) -- (-1.76471, 0.00000, 0.00000);
\draw[back,sourceedge] (-0.58824, -1.42857, -0.83333) -- (0.58824, -1.42857, 0.83333);
\draw[targetedge] (-0.58824, 1.42857, -0.83333) -- (0.58824, 1.42857, 0.83333);
\draw[targetedge] (-0.58824, 1.42857, -0.83333) -- (-1.76471, 0.00000, 0.00000);
\draw[targetedge] (0.58824, 1.42857, 0.83333) -- (1.76471, 0.00000, 0.00000);
\draw[sourceedge] (1.76471, 0.00000, 0.00000) -- (0.58824, -1.42857, 0.83333);
\node[targetvertex] at (1.76471, 0.00000, 0.00000) {};
\node[sourcevertex] at (-1.76471, 0.00000, 0.00000) {};

\begin{scope}[shift = {(.5,2,0)}]

		\draw[->] (0, 0,0) -- (1, 0,0) node[anchor = west]{$v_1$};
		\draw[->] (0, 0,0) -- (0, 1,0) node[anchor = west]{$v_2$};
	\end{scope}
\end{tikzpicture}&
&
\begin{tikzpicture}%
	[scale = .950000,
	back/.style = {thin, color = black!60},
	edge/.style = {color = black, thick},
	sourceedge/.style = {color = red, very thick},
	targetedge/.style = {color = blue, very thick},
	facet/.style = {fill = blue!95!black,fill opacity = 0.800000},
	targetfacet/.style = {fill = blue,fill opacity = 0.600000},
	sourcefacet/.style = {fill = red,fill opacity = 0.600000},
	vertex/.style = {inner sep = 1pt,circle,draw = black,fill = black,thick},
	targetvertex/.style = {inner sep = 1pt,circle,draw = blue,fill = blue,thick},
	sourcevertex/.style = {inner sep = 1pt,circle,draw = red,fill = red,thick}]
\draw[edge] (0,-2) -- (0,1);
\node[targetvertex] at (0,-2) {};
\node[sourcevertex] at (0,1) {};

\begin{scope}[shift = {(1,-1)}]

		\draw[->] (0, 0) -- ( 0,-1) node[anchor = west]{$v_1$};

	\end{scope}
\end{tikzpicture}
\\
$P = \pi_3(P)$& & $\pi_2(P)$& & $\pi_1(P)$
\end{tabular}

\begin{tikzpicture}%
	[	x = {(0.000000cm, -1.000000cm)},
	y = {(1.000000cm, 0.000000cm)},
	z = {(0.000000cm, 0.000000cm)},
	scale = .5,
	back/.style = {thin, color = black!60},
	edge/.style = {color = black!60, thin},
	sourceedge/.style = {color = red, ultra thick},
	targetedge/.style = {color = blue, ultra thick},
	facet/.style = {fill = black!35,fill opacity = 0.800000},
	targetfacet/.style = {fill = blue!80,fill opacity = 0.800000},
	sourcefacet/.style = {fill = red!80,fill opacity = 0.800000},
	vertex/.style = {inner sep = 1pt,circle,draw = black,fill = black,thick},
	targetvertex/.style = {inner sep = 1pt,circle,draw = blue,fill = blue,thick},
	sourcevertex/.style = {inner sep = 1pt,circle,draw = red,fill = red,thick}]
%
%
\coordinate (-0.58824, 1.42857, -0.83333) at (-0.58824, 1.42857, -0.83333);
\coordinate (0.58824, 1.42857, 0.83333) at (0.58824, 1.42857, 0.83333);
\coordinate (1.76471, 0.00000, 0.00000) at (1.76471, 0.00000, 0.00000);
\coordinate (0.58824, 0.00000, -1.66667) at (0.58824, 0.00000, -1.66667);
\coordinate (-0.58824, -1.42857, -0.83333) at (-0.58824, -1.42857, -0.83333);
\coordinate (-1.76471, 0.00000, 0.00000) at (-1.76471, 0.00000, 0.00000);
\coordinate (-0.58824, 0.00000, 1.66667) at (-0.58824, 0.00000, 1.66667);
\coordinate (0.58824, -1.42857, 0.83333) at (0.58824, -1.42857, 0.83333);


\draw[edge,back] (0.58824, 0.00000, -1.66667) -- (-0.58824, -1.42857, -0.83333);
\draw[edge,back] (-0.58824, -1.42857, -0.83333) -- (-1.76471, 0.00000, 0.00000);
\draw[edge,back] (-0.58824, -1.42857, -0.83333) -- (0.58824, -1.42857, 0.83333);

\fill[facet] (0.58824, 0.00000, -1.66667) -- (-0.58824, 1.42857, -0.83333) -- (0.58824, 1.42857, 0.83333) -- (1.76471, 0.00000, 0.00000) -- cycle {};
\draw[targetedge] (-0.58824, 1.42857, -0.83333) -- (0.58824, 1.42857, 0.83333);
\draw[sourceedge] (-0.58824, 1.42857, -0.83333) -- (0.58824, 0.00000, -1.66667);
\draw[edge] (-0.58824, 1.42857, -0.83333) -- (-1.76471, 0.00000, 0.00000);
\draw[targetedge] (0.58824, 1.42857, 0.83333) -- (1.76471, 0.00000, 0.00000);
\draw[edge] (0.58824, 1.42857, 0.83333) -- (-0.58824, 0.00000, 1.66667);
\draw[sourceedge] (1.76471, 0.00000, 0.00000) -- (0.58824, 0.00000, -1.66667);
\draw[edge] (1.76471, 0.00000, 0.00000) -- (0.58824, -1.42857, 0.83333);
\draw[edge] (-1.76471, 0.00000, 0.00000) -- (-0.58824, 0.00000, 1.66667);
\draw[edge] (-0.58824, 0.00000, 1.66667) -- (0.58824, -1.42857, 0.83333);
 \node[sourcevertex] at (-0.58824, 1.42857, -0.83333) {};
\node[targetvertex] at (1.76471, 0.00000, 0.00000) {};
%

\end{tikzpicture}
\quad
\begin{tikzpicture}%
	[x = {(0.000000cm, -1.000000cm)},
	y = {(1.000000cm, 0.000000cm)},
	z = {(0.000000cm, 0.000000cm)},,
	scale = .5,
	back/.style = {thin, color = black!60},
	edge/.style = {color = black!60, thin},
	sourceedge/.style = {color = red, ultra thick},
	targetedge/.style = {color = blue, ultra thick},
	facet/.style = {fill = black!35,fill opacity = 0.800000},
	targetfacet/.style = {fill = blue!80,fill opacity = 0.800000},
	sourcefacet/.style = {fill = red!80,fill opacity = 0.800000},
	vertex/.style = {inner sep = 1pt,circle,draw = black,fill = black,thick},
	targetvertex/.style = {inner sep = 1pt,circle,draw = blue,fill = blue,thick},
	sourcevertex/.style = {inner sep = 1pt,circle,draw = red,fill = red,thick}]
\coordinate (-0.58824, 1.42857, -0.83333) at (-0.58824, 1.42857, -0.83333);
\coordinate (0.58824, 1.42857, 0.83333) at (0.58824, 1.42857, 0.83333);
\coordinate (1.76471, 0.00000, 0.00000) at (1.76471, 0.00000, 0.00000);
\coordinate (0.58824, 0.00000, -1.66667) at (0.58824, 0.00000, -1.66667);
\coordinate (-0.58824, -1.42857, -0.83333) at (-0.58824, -1.42857, -0.83333);
\coordinate (-1.76471, 0.00000, 0.00000) at (-1.76471, 0.00000, 0.00000);
\coordinate (-0.58824, 0.00000, 1.66667) at (-0.58824, 0.00000, 1.66667);
\coordinate (0.58824, -1.42857, 0.83333) at (0.58824, -1.42857, 0.83333);

\fill[facet] (-0.58824, 1.42857, -0.83333) -- (-1.76471, 0.00000, 0.00000) -- (-0.58824, -1.42857, -0.83333) -- (0.58824, 0.00000, -1.66667) -- cycle {};

\draw[back,sourceedge] (0.58824, 0.00000, -1.66667) -- (-0.58824, -1.42857, -0.83333);
\draw[back,sourceedge] (-0.58824, -1.42857, -0.83333) -- (-1.76471, 0.00000, 0.00000);
\draw[edge,back] (-0.58824, -1.42857, -0.83333) -- (0.58824, -1.42857, 0.83333);

\draw[edge] (-0.58824, 1.42857, -0.83333) -- (0.58824, 1.42857, 0.83333);
\draw[targetedge] (-0.58824, 1.42857, -0.83333) -- (0.58824, 0.00000, -1.66667);
\draw[targetedge] (-0.58824, 1.42857, -0.83333) -- (-1.76471, 0.00000, 0.00000);
\draw[edge] (0.58824, 1.42857, 0.83333) -- (1.76471, 0.00000, 0.00000);
\draw[edge] (0.58824, 1.42857, 0.83333) -- (-0.58824, 0.00000, 1.66667);
\draw[edge] (1.76471, 0.00000, 0.00000) -- (0.58824, 0.00000, -1.66667);
\draw[edge] (1.76471, 0.00000, 0.00000) -- (0.58824, -1.42857, 0.83333);
\draw[edge] (-1.76471, 0.00000, 0.00000) -- (-0.58824, 0.00000, 1.66667);
\draw[edge] (-0.58824, 0.00000, 1.66667) -- (0.58824, -1.42857, 0.83333);
 \node[targetvertex] at (0.58824, 0.00000, -1.66667) {};
 \node[sourcevertex] at (-1.76471, 0.00000, 0.00000) {};
%

\end{tikzpicture}
\quad
\begin{tikzpicture}%
	[x = {(0.000000cm, -1.000000cm)},
	y = {(1.000000cm, 0.000000cm)},
	z = {(0.000000cm, 0.000000cm)},,
	scale = .5,
	back/.style = {thin, color = black!60},
	edge/.style = {color = black!60, thin},
	sourceedge/.style = {color = red, ultra thick},
	targetedge/.style = {color = blue, ultra thick},
	facet/.style = {fill = black!35,fill opacity = 0.800000},
	targetfacet/.style = {fill = blue!80,fill opacity = 0.800000},
	sourcefacet/.style = {fill = red!80,fill opacity = 0.800000},
	vertex/.style = {inner sep = 1pt,circle,draw = black,fill = black,thick},
	targetvertex/.style = {inner sep = 1pt,circle,draw = blue,fill = blue,thick},
	sourcevertex/.style = {inner sep = 1pt,circle,draw = red,fill = red,thick}]
\coordinate (-0.58824, 1.42857, -0.83333) at (-0.58824, 1.42857, -0.83333);
\coordinate (0.58824, 1.42857, 0.83333) at (0.58824, 1.42857, 0.83333);
\coordinate (1.76471, 0.00000, 0.00000) at (1.76471, 0.00000, 0.00000);
\coordinate (0.58824, 0.00000, -1.66667) at (0.58824, 0.00000, -1.66667);
\coordinate (-0.58824, -1.42857, -0.83333) at (-0.58824, -1.42857, -0.83333);
\coordinate (-1.76471, 0.00000, 0.00000) at (-1.76471, 0.00000, 0.00000);
\coordinate (-0.58824, 0.00000, 1.66667) at (-0.58824, 0.00000, 1.66667);
\coordinate (0.58824, -1.42857, 0.83333) at (0.58824, -1.42857, 0.83333);

 \fill[facet] (0.58824, 0.00000, -1.66667) -- (-0.58824, -1.42857, -0.83333) -- (0.58824, -1.42857, 0.83333) -- (1.76471, 0.00000, 0.00000) -- cycle {};

\draw[back,targetedge] (0.58824, 0.00000, -1.66667) -- (-0.58824, -1.42857, -0.83333);
\draw[edge,back] (-0.58824, -1.42857, -0.83333) -- (-1.76471, 0.00000, 0.00000);
\draw[back,sourceedge] (-0.58824, -1.42857, -0.83333) -- (0.58824, -1.42857, 0.83333);

\draw[edge] (-0.58824, 1.42857, -0.83333) -- (0.58824, 1.42857, 0.83333);
\draw[edge] (-0.58824, 1.42857, -0.83333) -- (0.58824, 0.00000, -1.66667);
\draw[edge] (-0.58824, 1.42857, -0.83333) -- (-1.76471, 0.00000, 0.00000);
\draw[edge] (0.58824, 1.42857, 0.83333) -- (1.76471, 0.00000, 0.00000);
\draw[edge] (0.58824, 1.42857, 0.83333) -- (-0.58824, 0.00000, 1.66667);
\draw[targetedge] (1.76471, 0.00000, 0.00000) -- (0.58824, 0.00000, -1.66667);
\draw[sourceedge] (1.76471, 0.00000, 0.00000) -- (0.58824, -1.42857, 0.83333);
\draw[edge] (-1.76471, 0.00000, 0.00000) -- (-0.58824, 0.00000, 1.66667);
\draw[edge] (-0.58824, 0.00000, 1.66667) -- (0.58824, -1.42857, 0.83333);
\node[sourcevertex] at (-0.58824, -1.42857, -0.83333) {};

\node[targetvertex] at (1.76471, 0.00000, 0.00000) {};
%

\end{tikzpicture}
\quad
\begin{tikzpicture}%
	[x = {(0.000000cm, -1.000000cm)},
	y = {(1.000000cm, 0.000000cm)},
	z = {(0.000000cm, 0.000000cm)},,
	scale = .5,
	back/.style = {thin, color = black!60},
	edge/.style = {color = black!60, thin},
	sourceedge/.style = {color = red, ultra thick},
	targetedge/.style = {color = blue, ultra thick},
	facet/.style = {fill = black!35,fill opacity = 0.800000},
	targetfacet/.style = {fill = blue!80,fill opacity = 0.800000},
	sourcefacet/.style = {fill = red!80,fill opacity = 0.800000},
	vertex/.style = {inner sep = 1pt,circle,draw = black,fill = black,thick},
	targetvertex/.style = {inner sep = 1pt,circle,draw = blue,fill = blue,thick},
	sourcevertex/.style = {inner sep = 1pt,circle,draw = red,fill = red,thick}]
\coordinate (-0.58824, 1.42857, -0.83333) at (-0.58824, 1.42857, -0.83333);
\coordinate (0.58824, 1.42857, 0.83333) at (0.58824, 1.42857, 0.83333);
\coordinate (1.76471, 0.00000, 0.00000) at (1.76471, 0.00000, 0.00000);
\coordinate (0.58824, 0.00000, -1.66667) at (0.58824, 0.00000, -1.66667);
\coordinate (-0.58824, -1.42857, -0.83333) at (-0.58824, -1.42857, -0.83333);
\coordinate (-1.76471, 0.00000, 0.00000) at (-1.76471, 0.00000, 0.00000);
\coordinate (-0.58824, 0.00000, 1.66667) at (-0.58824, 0.00000, 1.66667);
\coordinate (0.58824, -1.42857, 0.83333) at (0.58824, -1.42857, 0.83333);


\draw[edge,back] (0.58824, 0.00000, -1.66667) -- (-0.58824, -1.42857, -0.83333);
\draw[edge,back] (-0.58824, -1.42857, -0.83333) -- (-1.76471, 0.00000, 0.00000);
\draw[edge,back] (-0.58824, -1.42857, -0.83333) -- (0.58824, -1.42857, 0.83333);

\fill[facet] (0.58824, -1.42857, 0.83333) -- (1.76471, 0.00000, 0.00000) -- (0.58824, 1.42857, 0.83333) -- (-0.58824, 0.00000, 1.66667) -- cycle {};
\draw[edge] (-0.58824, 1.42857, -0.83333) -- (0.58824, 1.42857, 0.83333);
\draw[edge] (-0.58824, 1.42857, -0.83333) -- (0.58824, 0.00000, -1.66667);
\draw[edge] (-0.58824, 1.42857, -0.83333) -- (-1.76471, 0.00000, 0.00000);
\draw[targetedge] (0.58824, 1.42857, 0.83333) -- (1.76471, 0.00000, 0.00000);
\draw[targetedge] (0.58824, 1.42857, 0.83333) -- (-0.58824, 0.00000, 1.66667);
\draw[edge] (1.76471, 0.00000, 0.00000) -- (0.58824, 0.00000, -1.66667);
\draw[sourceedge] (1.76471, 0.00000, 0.00000) -- (0.58824, -1.42857, 0.83333);
\draw[edge] (-1.76471, 0.00000, 0.00000) -- (-0.58824, 0.00000, 1.66667);
\draw[sourceedge] (-0.58824, 0.00000, 1.66667) -- (0.58824, -1.42857, 0.83333);
\node[targetvertex] at (1.76471, 0.00000, 0.00000) {};
 \node[sourcevertex] at (-0.58824, 0.00000, 1.66667) {};
%

\end{tikzpicture}
\quad
\begin{tikzpicture}%
	[x = {(0.000000cm, -1.000000cm)},
	y = {(1.000000cm, 0.000000cm)},
	z = {(0.000000cm, 0.000000cm)},,
	scale = .5,
	back/.style = {thin, color = black!60},
	edge/.style = {color = black!60, thin},
	sourceedge/.style = {color = red, ultra thick},
	targetedge/.style = {color = blue, ultra thick},
	facet/.style = {fill = black!35,fill opacity = 0.800000},
	targetfacet/.style = {fill = blue!80,fill opacity = 0.800000},
	sourcefacet/.style = {fill = red!80,fill opacity = 0.800000},
	vertex/.style = {inner sep = 1pt,circle,draw = black,fill = black,thick},
	targetvertex/.style = {inner sep = 1pt,circle,draw = blue,fill = blue,thick},
	sourcevertex/.style = {inner sep = 1pt,circle,draw = red,fill = red,thick}]
\coordinate (-0.58824, 1.42857, -0.83333) at (-0.58824, 1.42857, -0.83333);
\coordinate (0.58824, 1.42857, 0.83333) at (0.58824, 1.42857, 0.83333);
\coordinate (1.76471, 0.00000, 0.00000) at (1.76471, 0.00000, 0.00000);
\coordinate (0.58824, 0.00000, -1.66667) at (0.58824, 0.00000, -1.66667);
\coordinate (-0.58824, -1.42857, -0.83333) at (-0.58824, -1.42857, -0.83333);
\coordinate (-1.76471, 0.00000, 0.00000) at (-1.76471, 0.00000, 0.00000);
\coordinate (-0.58824, 0.00000, 1.66667) at (-0.58824, 0.00000, 1.66667);
\coordinate (0.58824, -1.42857, 0.83333) at (0.58824, -1.42857, 0.83333);


\draw[edge,back] (0.58824, 0.00000, -1.66667) -- (-0.58824, -1.42857, -0.83333);
\draw[edge,back] (-0.58824, -1.42857, -0.83333) -- (-1.76471, 0.00000, 0.00000);
\draw[edge,back] (-0.58824, -1.42857, -0.83333) -- (0.58824, -1.42857, 0.83333);

 \fill[facet] (-0.58824, 0.00000, 1.66667) -- (0.58824, 1.42857, 0.83333) -- (-0.58824, 1.42857, -0.83333) -- (-1.76471, 0.00000, 0.00000) -- cycle {};
\draw[targetedge] (-0.58824, 1.42857, -0.83333) -- (0.58824, 1.42857, 0.83333);
\draw[edge] (-0.58824, 1.42857, -0.83333) -- (0.58824, 0.00000, -1.66667);
\draw[targetedge] (-0.58824, 1.42857, -0.83333) -- (-1.76471, 0.00000, 0.00000);
\draw[edge] (0.58824, 1.42857, 0.83333) -- (1.76471, 0.00000, 0.00000);
\draw[sourceedge] (0.58824, 1.42857, 0.83333) -- (-0.58824, 0.00000, 1.66667);
\draw[edge] (1.76471, 0.00000, 0.00000) -- (0.58824, 0.00000, -1.66667);
\draw[edge] (1.76471, 0.00000, 0.00000) -- (0.58824, -1.42857, 0.83333);
\draw[sourceedge] (-1.76471, 0.00000, 0.00000) -- (-0.58824, 0.00000, 1.66667);
\draw[edge] (-0.58824, 0.00000, 1.66667) -- (0.58824, -1.42857, 0.83333);
\node[targetvertex] at (0.58824, 1.42857, 0.83333) {};
 \node[sourcevertex] at (-1.76471, 0.00000, 0.00000) {};
%

\end{tikzpicture}
\quad
\begin{tikzpicture}%
	[x = {(0.000000cm, -1.000000cm)},
	y = {(1.000000cm, 0.000000cm)},
	z = {(0.000000cm, 0.000000cm)},,
	scale = .5,
	back/.style = {thin, color = black!60},
	edge/.style = {color = black!60, thin},
	sourceedge/.style = {color = red, ultra thick},
	targetedge/.style = {color = blue, ultra thick},
	facet/.style = {fill = black!35,fill opacity = 0.800000},
	targetfacet/.style = {fill = blue!80,fill opacity = 0.800000},
	sourcefacet/.style = {fill = red!80,fill opacity = 0.800000},
	vertex/.style = {inner sep = 1pt,circle,draw = black,fill = black,thick},
	targetvertex/.style = {inner sep = 1pt,circle,draw = blue,fill = blue,thick},
	sourcevertex/.style = {inner sep = 1pt,circle,draw = red,fill = red,thick}]
\coordinate (-0.58824, 1.42857, -0.83333) at (-0.58824, 1.42857, -0.83333);
\coordinate (0.58824, 1.42857, 0.83333) at (0.58824, 1.42857, 0.83333);
\coordinate (1.76471, 0.00000, 0.00000) at (1.76471, 0.00000, 0.00000);
\coordinate (0.58824, 0.00000, -1.66667) at (0.58824, 0.00000, -1.66667);
\coordinate (-0.58824, -1.42857, -0.83333) at (-0.58824, -1.42857, -0.83333);
\coordinate (-1.76471, 0.00000, 0.00000) at (-1.76471, 0.00000, 0.00000);
\coordinate (-0.58824, 0.00000, 1.66667) at (-0.58824, 0.00000, 1.66667);
\coordinate (0.58824, -1.42857, 0.83333) at (0.58824, -1.42857, 0.83333);

\fill[facet] (-0.58824, 0.00000, 1.66667) -- (0.58824, -1.42857, 0.83333) -- (-0.58824, -1.42857, -0.83333) -- (-1.76471, 0.00000, 0.00000) -- cycle {};

\draw[edge,back] (0.58824, 0.00000, -1.66667) -- (-0.58824, -1.42857, -0.83333);
\draw[back,sourceedge] (-0.58824, -1.42857, -0.83333) -- (-1.76471, 0.00000, 0.00000);
\draw[back,sourceedge] (-0.58824, -1.42857, -0.83333) -- (0.58824, -1.42857, 0.83333);

\draw[edge] (-0.58824, 1.42857, -0.83333) -- (0.58824, 1.42857, 0.83333);
\draw[edge] (-0.58824, 1.42857, -0.83333) -- (0.58824, 0.00000, -1.66667);
\draw[edge] (-0.58824, 1.42857, -0.83333) -- (-1.76471, 0.00000, 0.00000);
\draw[edge] (0.58824, 1.42857, 0.83333) -- (1.76471, 0.00000, 0.00000);
\draw[edge] (0.58824, 1.42857, 0.83333) -- (-0.58824, 0.00000, 1.66667);
\draw[edge] (1.76471, 0.00000, 0.00000) -- (0.58824, 0.00000, -1.66667);
\draw[edge] (1.76471, 0.00000, 0.00000) -- (0.58824, -1.42857, 0.83333);
\draw[targetedge] (-1.76471, 0.00000, 0.00000) -- (-0.58824, 0.00000, 1.66667);
\draw[targetedge] (-0.58824, 0.00000, 1.66667) -- (0.58824, -1.42857, 0.83333);


 \node[sourcevertex] at (-1.76471, 0.00000, 0.00000) {};
 \node[targetvertex] at (0.58824, -1.42857, 0.83333) {};
%
\begin{scope}[shift = {(.5,2,0)}]

		\draw[->] (0, 0,0) -- (1, 0,0) node[anchor = west]{$v_1$};
		\draw[->] (0, 0,0) -- (0, 1,0) node[anchor = west]{$v_2$};
	\end{scope}
\end{tikzpicture}

%% file: sec/subdivisions.tex

\subsection{Coherent subdivisions}\label{s:subdivision}

\SS A \defn{polytopal complex} $\cP$ is a finite collection of polytopes in $\R^d$ satisfying:
\begin{enumerate}
	\item\label{item:downwards closure of subdivisions} If $P \in \cP$ then all the faces of $P$ are in $\cP$.
	\item If $P, P' \in \cP$ then $P \cap P'$ is a face of both $P$ and $P'$.
\end{enumerate}
Since a polytopal complex is determined by its inclusion-maximal elements, we identify it with that set.
A \defn{subdivision} of a polytope~$Q \in \R^d$ is a polytopal complex $\cS$ whose union is $Q$.

Let $\pi \colon P \subset \R^p \to Q \subset \R^q$ be a \defn{projection of polytopes}, i.e., an affine map \(\pi \colon \R^p \to \R^q\) with \(\pi(P) = Q\), where \(P \subseteq \R^p\) and \(Q \subseteq \R^q\) are polytopes.
A \defn{\(\pi\)-induced subdivision} is a subdivision of \(Q\) of the form \(\set{\pi(F) \mid F \in \cS}\) for some polytopal complex \(\cS \subseteq \faces\) satisfying:
\begin{enumerate}[start=3]
	\item If $F, F' \in \cS$ and $\pi(F) \subseteq \pi(F')$ then $F = F' \cap \pi^{-1}(\pi(F))$.
\end{enumerate}
Since each polytope in this subdivision corresponds to a unique face \(F \in \cS\), we will refer to \(\cS \subseteq \faces\) as the \(\pi\)-induced subdivision itself.
If $\dim F = \dim \pi(F)$ for all $F \in \cS$ the subdivision $\cS$ is said to be \defn{tight}.

To each functional $\psi \in (\R^p)^*$ in the dual of the domain of $\pi$ we can assign a \(\pi\)-induced subdivision $\cS^\psi$ consisting of the faces $F$ of \(P\) such that the maximal locus of $\psi$ in $\pi^{-1}(x) \cap P$ is $\pi^{-1}(x) \cap F$ for all $x \in Q$.
Such subdivisions are called \defn{\(\pi\)-coherent subdivisions}.

Given a frame, let us denote its dual basis by \( v_1^*, v_2^*, \dots \).

\SS\theorem\label{t:coherent_subdivisions}
For any framed polytope $P$, the sets $\so_k(P)$ and $\ta_k(P)$ are the tight $\pi_k$-coherent subdivisions respectively minimizing and maximizing the functional~$v_{k+1}^*$.

\begin{proof}
	We factor the projection
	\[
	\pi_k \colon P \xra{\pi_{k+1}} \pi_{k+1}(P) \xra{\pi} \pi_k(P).
	\]
	For any \( y \in \R^d \), we have \( \psi(y) = \psi(\pi_{k+1}(y)) \); thus, for any \( x \in \pi_k(P) \), it follows that \( \pi^{-1}(x)^\psi = \pi_{k+1}(\pi_k^{-1}(x)^\psi) \).
	Thus, if \( F \) is the face of \( P \) that maximizes \( \psi \) in \( \pi_k^{-1}(x) \), then \( \pi_{k+1}(F) \) maximizes \( \psi \) in \( \pi^{-1}(x) \).

	Note that \( \pi^{-1}(x) \) is parallel to \( v_{k+1} \); this fiber is the intersection of the zero loci of \( v_1^*, \dots, v_k^* \).
	The maximal locus of \( \psi \) in \( \pi^{-1}(x) \) is the face whose normal covector is positively proportional to \( v_{k+1}^* \).
	Now, this face will be of the form \( \pi_{k+1} F \cap \pi^{-1}(x) \) for some face \( F \) of \( P \).
	Its normal covector is the projection of \( \normal_{\pi_{k+1}(F)}^{\pi_{k+1}(P)} \) onto the span of \( v_{k+1}^* \) along \( v_1^*, \dots, v_k^* \).
	It will be positively proportional to \( v_{k+1}^* \) whenever \( \sprod{\normal_{\pi_{k+1}(F)}^{\pi_{k+1}(P)}}{v_{k+1}} > 0 \).
	This is equivalent to \( F \) belonging to \( \ta_k(P) \).

	The same argument, applied with \( v_{k+1} \) replaced by \( -v_{k+1} \), establishes the result for~\( \so_k(P) \).

	Tightness follows from the fact that the dimension of the faces in \( \so_k(P) \) and \( \ta_k(P) \) equals the dimension of~\( \pi_k(P) \).
\end{proof}

%% file: sec/strings.tex

\section{Cellular strings and loops}\label{part:strings}

In a framed polytope, \(k\)-cellular strings are sequences of faces glued along their \(k\)-sources and \(k\)-targets.
These strings can form loops, which, as we will see in \cref{part:connection}, constitute the only obstruction to the faces of a framed \(n\)-polytope freely generating an \(n\)-category.
We start by showing in \cref{ss:low-dimensions} that loops cannot occur in polytopes of dimensions less than 4.
Then, we explicitly exhibit in \cref{s:examples-loops} a framed \(4\)-polytope with a cellular \(1\)-loop and a framed \(6\)-polytope with a cellular \(2\)-loop.
In fact, loops are abundant in framed polytopes.
As we show in \cref{sec:bootstrap}, every simplicial or simple polytope of dimension \(\geq 6\) admits a frame inducing a loop.
\cref{s:loop-inevitability} is devoted to the construction of a polytope for which all admissible frames induce a cellular loop.
Finally, we show in \cref{s:random} that, in contrast to the cyclic cases treated in \cref{part:orientals}, a canonically framed Gaussian \(n\)-simplex and the standard \(n\)-cube endowed with a random frame have cellular loops with probability tending to \(1\) as \(n \to \infty\).

\input{sec/strings_loops}
\input{sec/lowdim}
\input{sec/loops}
\input{sec/bootstrapping}
\input{sec/noframe}
\input{sec/random}

%% file: sec/strings_loops.tex

\subsection{Cellular strings}\label{ss:cellular_strings}

\begin{figure}[htpb]
	\centering
	\input{fig/string_figure}
	\caption{A cellular \(1\)-string and a cellular \(0\)-string.}
	\label{fig:string-poly}
\end{figure}
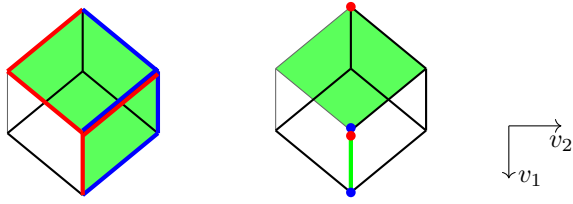

\SS\label{def:cellular-string}
Let \(P\) be a framed polytope.
A \defn{cellular \(k\)-string} is a sequence \(F_1,\dots,F_m\) of faces of \(P\) of dimension greater than~\(k\) such that for each \(i \in \set{1,\dots,m-1}\), we have \(\ta_k(F_i) \cap \so_k(F_{i+1}) \neq \emptyset\).
Examples are depicted in \cref{fig:string-poly}.

\SS\label{lem:intersection-in-string}
\lemma The intersection of two consecutive faces \(F_i\) and \(F_{i+1}\) in a cellular \(k\)-string is a \(k\)-face.

\begin{proof}
	A hyperplane \textit{weakly separates} two polytopes if each of their interiors are in distinct half-spaces defined by \(H\).
	Let \(G\) be a \(k\)-face of \(\ta_k (F_i) \cap \so_k (F_{i+1})\).
	Then \(\pi_{k+1}G\) is a common facet of \(\pi_{k+1} F_i\) and \(\pi_{k+1} F_{i+1}\).
	By \cref{lem:sources and targets in a framed polytope}, the normal covectors \(\normal_{\pi_{k+1} G}^{\pi_{k+1} F_i}\) and \(\normal_{\pi_{k+1} G}^{\pi_{k+1} F_{i+1}}\) have opposite orientations.
	Thus, the hyperplane~\(H\) normal to \(\normal_{\pi_{k+1} G}^{\pi_{k+1} F_i}\) and passing through \(\pi_{k+1}G\) weakly separates \(\pi_{k+1} F_{i}\) and \(\pi_{k+1} F_{i+1}\).
	Therefore, its preimage \(H' = \pi_{k+1}^{-1}(H)\) weakly separates \(F_i\) and \(F_{i+1}\), and we have \(F_i \cap F_{i+1}\subseteq H'\).
	Now, the intersections \(F_i \cap H'\) and \(F_{i+1} \cap H'\) are faces of \(F_i\) and \(F_{i+1}\), respectively, but they cannot be of dimension larger than \(k\) because this would contradict the \(P\)-admissibility assumption.
	Therefore \(F_i \cap H' = F_{i+1} \cap H' = F_i \cap F_{i+1} = G\).
\end{proof}

\SS\remark A cellular \(0\)-string from \(\so_0(P)\) to \(\ta_0(P)\) is a cellular string in the sense of~\cite{billera1994strings}.
More precisely, assuming that the frame \((v_1,v_2,\ldots)\) is orthogonal (\cref{ss:equivalentorthogonal}), the \emph{stippling} defined in \cite[Sec.~1]{billera1994strings} is given by the linear functional \(v_1^*\).

\SS\label{ss:inhomogeneous-cellular-string}
Let \(P\) be a framed polytope.
An \defn{inhomogeneous cellular string} is a sequence \(F_1,\ldots,F_m\) of faces of \(P\) such that for each \(i \in \set{1,\ldots, m-1}\) we have either \(F_i \in \so(F_{i+1})\) or \(\ta(F_i) \ni F_{i+1}\).

The two notions of cellular \(k\)-string and inhomogeneous cellular string are closely related, as we now show.

\SS\label{l:target-in-the-target}
\lemma Let \(P\) be a framed polytope.
Let \(F\) and \(H\) be \(k\)- and \(\ell\)-dimensional faces of \(P\), respectively, such that
\(k > \ell+1\) and \(H \in \ta_\ell(F)\).
Then, there exists a face \(G\) in \(\ta(F)\) such that \(H \in \ta_\ell(G)\).

\begin{proof}
	As we have seen in \cref{s:subdivision}, the projection \(\pi_{k-1}(\ta(F))\) is a (tight coherent) subdivision of \(\pi_{k-1}(F)\) by \((k-1)\)-dimensional polytopes.
	Combined with the fact that \(\pi_{\ell+1} (\pi_{k-1} (F)) = \pi_{\ell+1} (F)\), this implies that there must be a face~\(G\) in~\(\ta(F)\) whose \(\ell\)-boundary \(\bd_{\ell}(G)\) contains \(\pi_{\ell+1}(H)\).
	Moreover, we must have \(H \subset G\), otherwise we would have a higher-dimensional face projecting onto \(\pi_{\ell+1}(H)\), which cannot happen since the basis is \(F\)-admissible and \(H\) is by assumption in the \(\ell\)-boundary \(\bd_{\ell}(F)\) of~\(F\) (see \cref{s:sk_tk}).
	Therefore, we have found a face~\(G\) in \(\ta(F)\) such that \(H \in \ta_\ell(G)\).
\end{proof}

\SS\label{l:inhomogeneous strings}
\proposition
Let \(P\) be a framed polytope.
If \(F_1,\dots,F_m\) is a cellular \(k\)-string, then there exists an inhomogeneous cellular string \(G_1,\dots,G_\ell\) and an increasing injection \(\phi \colon \set{1,\dots,m} \to \set{1,\dots,\ell}\) such that \(F_i = G_{\phi(i)}\).

\begin{proof}
	Consider two consecutive faces \(F_i, F_{i+1}\) in a cellular \(k\)-string.
	Applying recursively \cref{l:target-in-the-target}, we can find a family of faces \(G_1,\ldots,G_\ell\) such that \(G_1 = F_i\), \(G_\ell = F_i \cap F_{i+1}\), and for each \(j \in [\ell-1]\) we have \(\ta(G_j) \ni G_{j+1}\).
	Reasoning similarly by replacing targets with sources, we can find a family of faces \(G_1',\ldots,G_n'\) such that \(G_1' = F_i \cap F_{i+1}\), \(G_n' = F_{i+1}\), and for each \(j \in [n-1]\) we have \(G_j' \in \so(G_{j+1}')\).
	Combining these two sequences of faces, we get an inhomogeneous cellular string between \(F_i\) and \(F_{i+1}\).
	Repeating this process for every pair of consecutive faces in the cellular string, we get the desired inhomogeneous cellular string.
\end{proof}

\subsection{Cellular loops}

\SS\label{def:cellular-loop}
A \defn{cellular \(k\)-loop} is a cellular \(k\)-string \(F_1,\dots,F_m\) with \(F_i = F_j\) for some \(i \neq j\).
The notion of \defn{inhomogeneous cellular loop} is defined analogously.
We say that a framed polytope is \defn{\(k\)-loop-free} if it does not admit any cellular \(k\)-loops, and \defn{loop-free} if it is \(k\)-loop-free for any \(k \in \N\).
Moreover, we say that it is \defn{strongly loop-free} if it does not admit any inhomogeneous cellular loop.

By \cref{l:inhomogeneous strings}, the existence of a cellular \(k\)-loop implies that of an inhomogeneous cellular loop.
However, as the next example shows, there are framed polytopes which admit inhomogeneous cellular loops but no cellular loops.
Thus, any strongly loop-free framed polytope is also loop-free, but the converse is not true.

\begin{figure}[h!]
	\tdplotsetmaincoords{75}{75}
	\input{fig/cross_poly}
	\caption{The \(3\)-dimensional cross-polytope \(P\) admits a basis which is loop-free, but not strongly loop-free.}
	\label{fig:cross-poly}
\end{figure}
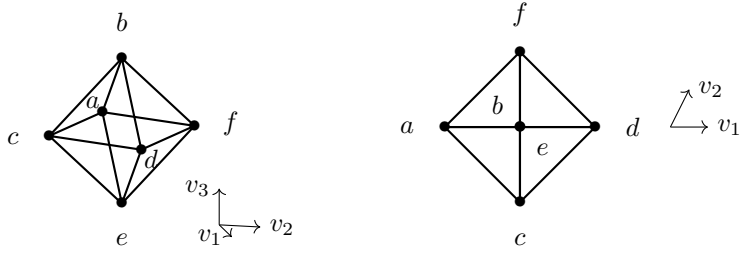

\medskip\example Consider the \(3\)-dimensional cross-polytope \(P\) in \(\R^3\), defined as the convex hull of the points \((a,b,c,d,e,f) \defeq (-e_1,e_3,-e_2,e_1,-e_3,e_2)\), together with the frame \((v_1,v_2,v_3) \defeq (e_1,\tfrac{1}{2}e_1+e_2,e_3)\).
It is clear, by inspection, that this frame is \(P\)-admissible.
By the upcoming \cref{c:up-to-dim-3}, we have that \(P\) is loop-free.
However, it admits the following inhomogeneous loop, represented in \cref{fig:cross-poly}:
\[
[abc] , \ [bc] , \ [c] , \ [cd] , \ [cde] , \ P , \ [abc].
\]

\SS\label{cor:projected-loop-free}
In \cref{ss:cubes} and \cref{ss:simplex} we will need the following.

\medskip\lemma If~\(P\) is framed and (strongly) loop-free, then so is~\(\pi_k(P)\) for all~\(k\).

\begin{proof}
	By construction, the projection \(\pi_k(P)\) of \(P\) is framed by \((v_1,\ldots,v_k)\).
	Moreover, if \(F\) is a face of \(P\) such that \(G \defeq \pi_k(F)\) is a face of \(\pi_k(P)\), then \(\ta_\ell(G) = \pi_k(\ta_\ell(F))\) and \(\so_\ell(G) = \pi_k(\so_\ell(F))\) for all \(\ell \leq k\).
	Therefore, any (inhomogeneous) loop in \(\pi_k(P)\) can be lifted to an (inhomogeneous) loop in \(P\).
\end{proof}

%% file: fig/string_figure.tex
 \begin{tikzpicture}%
	[
	x = {(0.000000cm, -1.000000cm)},
	y = {(1.000000cm, 0.000000cm)},
	z = {(0.000000cm, 0.000000cm)},
	scale = .700000,
 	back/.style = {thin, color = black!60},
	edge/.style = {color = black, thick},
	sourceedge/.style = {color = red, ultra thick},
	targetedge/.style = {color = blue, ultra thick},
	selectededge/.style = {color = green, ultra thick},
	facet/.style = {fill = white,fill opacity = 0.00000},
	targetfacet/.style = {fill = blue!80,fill opacity = 0.800000},
	sourcefacet/.style = {fill = red!80,fill opacity = 0.800000},
	selectedfacet/.style = {fill = green!80,fill opacity = 0.800000},
	vertex/.style = {inner sep = 1pt,circle,draw = black,fill = black,thick},
	targetvertex/.style = {inner sep = 1pt,circle,draw = blue,fill = blue,thick},
	sourcevertex/.style = {inner sep = 1pt,circle,draw = red,fill = red,thick}]
%
%
\coordinate (-0.58824, 1.42857, -0.83333) at (-0.58824, 1.42857, -0.83333);
\coordinate (0.58824, 1.42857, 0.83333) at (0.58824, 1.42857, 0.83333);
\coordinate (1.76471, 0.00000, 0.00000) at (1.76471, 0.00000, 0.00000);
\coordinate (0.58824, 0.00000, -1.66667) at (0.58824, 0.00000, -1.66667);
\coordinate (-0.58824, -1.42857, -0.83333) at (-0.58824, -1.42857, -0.83333);
\coordinate (-1.76471, 0.00000, 0.00000) at (-1.76471, 0.00000, 0.00000);
\coordinate (-0.58824, 0.00000, 1.66667) at (-0.58824, 0.00000, 1.66667);
\coordinate (0.58824, -1.42857, 0.83333) at (0.58824, -1.42857, 0.83333);

\fill[facet] (-0.58824, 0.00000, 1.66667) -- (0.58824, -1.42857, 0.83333) -- (-0.58824, -1.42857, -0.83333) -- (-1.76471, 0.00000, 0.00000) -- cycle {};
\fill[selectedfacet] (-0.58824, 1.42857, -0.83333) -- (-1.76471, 0.00000, 0.00000) -- (-0.58824, -1.42857, -0.83333) -- (0.58824, 0.00000, -1.66667) -- cycle {};
 \fill[facet] (0.58824, 0.00000, -1.66667) -- (-0.58824, -1.42857, -0.83333) -- (0.58824, -1.42857, 0.83333) -- (1.76471, 0.00000, 0.00000) -- cycle {};

\draw[edge,back] (0.58824, 0.00000, -1.66667) -- (-0.58824, -1.42857, -0.83333);
\draw[edge,back] (-0.58824, -1.42857, -0.83333) -- (-1.76471, 0.00000, 0.00000);
\draw[edge,back] (-0.58824, -1.42857, -0.83333) -- (0.58824, -1.42857, 0.83333);

\fill[selectedfacet] (0.58824, 0.00000, -1.66667) -- (-0.58824, 1.42857, -0.83333) -- (0.58824, 1.42857, 0.83333) -- (1.76471, 0.00000, 0.00000) -- cycle {};
\fill[facet] (0.58824, -1.42857, 0.83333) -- (1.76471, 0.00000, 0.00000) -- (0.58824, 1.42857, 0.83333) -- (-0.58824, 0.00000, 1.66667) -- cycle {};
 \fill[facet] (-0.58824, 0.00000, 1.66667) -- (0.58824, 1.42857, 0.83333) -- (-0.58824, 1.42857, -0.83333) -- (-1.76471, 0.00000, 0.00000) -- cycle {};
\draw[edge] (-0.58824, 1.42857, -0.83333) -- (0.58824, 1.42857, 0.83333);
\draw[edge] (-0.58824, 1.42857, -0.83333) -- (0.58824, 0.00000, -1.66667);
\draw[edge] (-0.58824, 1.42857, -0.83333) -- (-1.76471, 0.00000, 0.00000);
\draw[edge] (0.58824, 1.42857, 0.83333) -- (1.76471, 0.00000, 0.00000);
\draw[edge] (0.58824, 1.42857, 0.83333) -- (-0.58824, 0.00000, 1.66667);
\draw[edge] (1.76471, 0.00000, 0.00000) -- (0.58824, 0.00000, -1.66667);
\draw[edge] (1.76471, 0.00000, 0.00000) -- (0.58824, -1.42857, 0.83333);
\draw[edge] (-1.76471, 0.00000, 0.00000) -- (-0.58824, 0.00000, 1.66667);
\draw[edge] (-0.58824, 0.00000, 1.66667) -- (0.58824, -1.42857, 0.83333);
%

\draw[back,sourceedge] (0.58824, 0.00000, -1.66667) -- (-0.58824, -1.42857, -0.83333);
\draw[back,sourceedge] (-0.58824, -1.42857, -0.83333) -- (-1.76471, 0.00000, 0.00000);

\draw[targetedge] (-0.58824, 1.42857, -0.83333) -- (-1.76471, 0.00000, 0.00000);

 \begin{scope}[shift = {(-0.05,-0.0,-0.0)}]
\draw[targetedge] (-0.58824, 1.42857, -0.83333) -- (0.58824, 0.00000, -1.66667);
 	\end{scope}

\draw[targetedge] (-0.58824, 1.42857, -0.83333) -- (0.58824, 1.42857, 0.83333);
\draw[targetedge] (0.58824, 1.42857, 0.83333) -- (1.76471, 0.00000, 0.00000);
\draw[sourceedge] (1.76471, 0.00000, 0.00000) -- (0.58824, 0.00000, -1.66667);
\begin{scope}[shift = {(0.05,+0,+0)}]
\draw[sourceedge] (-0.58824, 1.42857, -0.83333) -- (0.58824, 0.00000, -1.66667);
\end{scope}
\end{tikzpicture}
\qquad\qquad
 \begin{tikzpicture}%
	[
	x = {(0.000000cm, -1.000000cm)},
	y = {(1.000000cm, 0.000000cm)},
	z = {(0.000000cm, 0.000000cm)},
	scale = .700000,
 	back/.style = {thin, color = black!60},
	edge/.style = {color = black, thick},
	sourceedge/.style = {color = red, ultra thick},
	targetedge/.style = {color = blue, ultra thick},
	selectededge/.style = {color = green, ultra thick},
	facet/.style = {fill = white,fill opacity = 0.00000},
	targetfacet/.style = {fill = blue!80,fill opacity = 0.800000},
	sourcefacet/.style = {fill = red!80,fill opacity = 0.800000},
	selectedfacet/.style = {fill = green!80,fill opacity = 0.800000},
	vertex/.style = {inner sep = 1pt,circle,draw = black,fill = black,thick},
	targetvertex/.style = {inner sep = 1pt,circle,draw = blue,fill = blue,thick},
	sourcevertex/.style = {inner sep = 1pt,circle,draw = red,fill = red,thick}]
%
%
\coordinate (-0.58824, 1.42857, -0.83333) at (-0.58824, 1.42857, -0.83333);
\coordinate (0.58824, 1.42857, 0.83333) at (0.58824, 1.42857, 0.83333);
\coordinate (1.76471, 0.00000, 0.00000) at (1.76471, 0.00000, 0.00000);
\coordinate (0.58824, 0.00000, -1.66667) at (0.58824, 0.00000, -1.66667);
\coordinate (-0.58824, -1.42857, -0.83333) at (-0.58824, -1.42857, -0.83333);
\coordinate (-1.76471, 0.00000, 0.00000) at (-1.76471, 0.00000, 0.00000);
\coordinate (-0.58824, 0.00000, 1.66667) at (-0.58824, 0.00000, 1.66667);
\coordinate (0.58824, -1.42857, 0.83333) at (0.58824, -1.42857, 0.83333);

\fill[facet] (-0.58824, 0.00000, 1.66667) -- (0.58824, -1.42857, 0.83333) -- (-0.58824, -1.42857, -0.83333) -- (-1.76471, 0.00000, 0.00000) -- cycle {};
\fill[selectedfacet] (-0.58824, 1.42857, -0.83333) -- (-1.76471, 0.00000, 0.00000) -- (-0.58824, -1.42857, -0.83333) -- (0.58824, 0.00000, -1.66667) -- cycle {};
 \fill[facet] (0.58824, 0.00000, -1.66667) -- (-0.58824, -1.42857, -0.83333) -- (0.58824, -1.42857, 0.83333) -- (1.76471, 0.00000, 0.00000) -- cycle {};

\draw[edge,back] (0.58824, 0.00000, -1.66667) -- (-0.58824, -1.42857, -0.83333);
\draw[edge,back] (-0.58824, -1.42857, -0.83333) -- (-1.76471, 0.00000, 0.00000);
\draw[edge,back] (-0.58824, -1.42857, -0.83333) -- (0.58824, -1.42857, 0.83333);

\fill[facet] (0.58824, 0.00000, -1.66667) -- (-0.58824, 1.42857, -0.83333) -- (0.58824, 1.42857, 0.83333) -- (1.76471, 0.00000, 0.00000) -- cycle {};
\fill[facet] (0.58824, -1.42857, 0.83333) -- (1.76471, 0.00000, 0.00000) -- (0.58824, 1.42857, 0.83333) -- (-0.58824, 0.00000, 1.66667) -- cycle {};
 \fill[facet] (-0.58824, 0.00000, 1.66667) -- (0.58824, 1.42857, 0.83333) -- (-0.58824, 1.42857, -0.83333) -- (-1.76471, 0.00000, 0.00000) -- cycle {};
\draw[edge] (-0.58824, 1.42857, -0.83333) -- (0.58824, 1.42857, 0.83333);
\draw[edge] (-0.58824, 1.42857, -0.83333) -- (0.58824, 0.00000, -1.66667);
\draw[edge] (-0.58824, 1.42857, -0.83333) -- (-1.76471, 0.00000, 0.00000);
\draw[edge] (0.58824, 1.42857, 0.83333) -- (1.76471, 0.00000, 0.00000);
\draw[edge] (0.58824, 1.42857, 0.83333) -- (-0.58824, 0.00000, 1.66667);
\draw[selectededge] (1.76471, 0.00000, 0.00000) -- (0.58824, 0.00000, -1.66667);
\draw[edge] (1.76471, 0.00000, 0.00000) -- (0.58824, -1.42857, 0.83333);
\draw[edge] (-1.76471, 0.00000, 0.00000) -- (-0.58824, 0.00000, 1.66667);
\draw[edge] (-0.58824, 0.00000, 1.66667) -- (0.58824, -1.42857, 0.83333);
\node[targetvertex] at (1.76471, 0.00000, 0.00000) {};
%

 \node[sourcevertex] at (-1.76471, 0.00000, 0.00000) {};

 \begin{scope}[shift = {(-0.05,-0.0,-0.0)}]
 \node[targetvertex] at (0.58824, 0.00000, -1.66667) {};

 	\end{scope}

\begin{scope}[shift = {(0.1,+0,+0)}]
\node[sourcevertex] at (0.58824, 0.00000, -1.66667) {};
\end{scope}
\begin{scope}[shift = {(.5,3,0)}]

		\draw[->] (0, 0,0) -- (1, 0,0) node[anchor = west]{$v_1$};
		\draw[->] (0, 0,0) -- (0, 1,0) node[anchor = north]{$v_2$};
	\end{scope}
\end{tikzpicture}

%% file: fig/cross_poly.tex
\begin{tikzpicture}
	[tdplot_main_coords,
	cube/.style = {very thick,black},
	grid/.style = {very thin,gray},
	axis/.style = {->}]

	\node (a) at (-1,0,0) {$\bullet$};
	\node at (-1.5,0,0) {$a$};

	\node (b) at (0,0,1) {$\bullet$};
	\node at (0,0,1.5) {$b$};

	\node (c) at (0,-1,0) {$\bullet$};
	\node at (0,-1.5,0) {$c$};

	\node (d) at (1,0,0) {$\bullet$};
	\node at (1.5,0,0) {$d$};

	\node (e) at (0,0,-1) {$\bullet$};
	\node at (0,0,-1.5) {$e$};

	\node (f) at (0,1,0) {$\bullet$};
	\node at (0,1.5,0) {$f$};

	\draw[axis] (5,0,0) -- (5.65,0,0) node[anchor = east]{$v_1$};
	\draw[axis] (5,0,0) -- (5.25,0.5,0) node[anchor = west]{$v_2$};
	\draw[axis] (5,0,0) -- (5,0,0.5) node[anchor = east]{$v_3$};


	\draw[-,thick] (1,0,0)--(0,1,0);
	\draw[-,thick] (1,0,0)--(0,-1,0);
	\draw[-,thick] (1,0,0)--(0,0,1);
	\draw[-,thick] (1,0,0)--(0,0,-1);
	\draw[-,thick] (-1,0,0)--(0,1,0);
	\draw[-,thick] (-1,0,0)--(0,-1,0);
	\draw[-,thick] (-1,0,0)--(0,0,1);
	\draw[-,thick] (-1,0,0)--(0,0,-1);
	\draw[-,thick] (0,1,0)--(0,0,1);
	\draw[-,thick] (0,-1,0)--(0,0,1);
	\draw[-,thick] (0,1,0)--(0,0,-1);
	\draw[-,thick] (0,-1,0)--(0,0,-1);


\end{tikzpicture}
\quad \quad
\tdplotsetmaincoords{0}{0}
\begin{tikzpicture}
	[tdplot_main_coords,
	cube/.style = {very thick,black},
	grid/.style = {very thin,gray},
	axis/.style = {->}]

	\node (a) at (-1,0,0) {$\bullet$};
	\node at (-1.5,0,0) {$a$};

	\node (b) at (0,0,1) {$\bullet$};
	\node at (-0.3,0.3,1.5) {$b$};

	\node (c) at (0,-1,0) {$\bullet$};
	\node at (0,-1.5,0) {$c$};

	\node (d) at (1,0,0) {$\bullet$};
	\node at (1.5,0,0) {$d$};

	\node (e) at (0,0,-1) {$\bullet$};
	\node at (0.3,-0.3,-1.5) {$e$};

	\node (f) at (0,1,0) {$\bullet$};
	\node at (0,1.5,0) {$f$};

	\draw[axis] (2,0,0) -- (2.5,0,0) node[anchor = west]{$v_1$};
	\draw[axis] (2,0,0) -- (2.25,0.5,0) node[anchor = west]{$v_2$};

	\draw[-,thick] (1,0,0)--(0,1,0);
	\draw[-,thick] (1,0,0)--(0,-1,0);
	\draw[-,thick] (1,0,0)--(0,0,1);
	\draw[-,thick] (1,0,0)--(0,0,-1);
	\draw[-,thick] (-1,0,0)--(0,1,0);
	\draw[-,thick] (-1,0,0)--(0,-1,0);
	\draw[-,thick] (-1,0,0)--(0,0,1);
	\draw[-,thick] (-1,0,0)--(0,0,-1);
	\draw[-,thick] (0,1,0)--(0,0,1);
	\draw[-,thick] (0,-1,0)--(0,0,1);
	\draw[-,thick] (0,1,0)--(0,0,-1);
	\draw[-,thick] (0,-1,0)--(0,0,-1);
\end{tikzpicture}

%% file: sec/lowdim.tex

\subsection{Low dimensions} \label{ss:low-dimensions}

\SS\label{prop:dim0}
Looping behavior of cellular strings is a higher-dimensional phenomenon, as shown in the following.

\medskip\proposition There are no cellular $0$-loops on any framed polytope.

\begin{proof}
	Note that for any face $F$, we have $\sprod{\pi_1(\so_0(F))}{v_1^*}<\sprod{\pi_1(\ta_0(F))}{v_1^*}$, which shows that there cannot be repeated elements in a $0$-string.
\end{proof}

\SS\label{prop:dimP-2}

For a framed polytope of dimension \(d\), the largest \(k\) for which there are non-trivial cellular \(k\)-strings is \(d - 2\), and we have the following.

\medskip\proposition A framed polytope $P$ has no cellular $(\dim P-2)$-loops.

\begin{proof}
	Denote by $d$ the dimension of $P$.
	Note first that any $(d-2)$-face $F$ that belongs to the $(d-2)$-boundary $\bd_{d-2} P$ (i.e.\ which projects onto a facet of $\pi_{d-1} P$) cannot be both a $(d-2)$-source and target, because all the facets that contain $F$, once projected via $\pi_{d-1}$, lie on the same side of the affine hyperplane $\Aff F \subset V_{d-1}$. Therefore a $(\dim P-2)$-loop can only contain faces either of the $(d-1)$-source or of the $(d-1)$-target, but not both. Assume without loss of generality that it is only the $(d-1)$-target.
	By~\cref{s:subdivision}, the projection $\pi_{d-1}(\ta_{d-1}(P))$ of the hemisphere~$\ta_{d-1}(P)$ gives a coherent subdivision $\mathcal{F}$ of $\pi_{d-1}(P)$. By \cite[Thm.~7.1]{mcmullenFibreTilings2003}, the $\pi_{d-1}$-coherent subdivision $\mathcal{F}$ admits a ``strong dual'' \cite[Sec.~6]{mcmullenFibreTilings2003}.
	That is, there is a dual tiling $\mathcal{F}^{*}$ of $V_{d-1}$ whose $0$-faces and $1$-faces are in bijection with the $(d-1)$- and $(d-2)$-faces of $\mathcal{F}$, respectively, and whose $1$-faces are moreover orthogonal to the $(d-2)$-faces of $\mathcal{F}$.
	Then, any $(d-1)$-string $F_1,\ldots,F_m$ in $\mathcal{F}$ defines a $0$-string in $\mathcal{F}^{*}$.
	Since the frame $(v_1,\ldots,v_d)$ is $P$-admissible, the function $v^*_{d-1}$ is strictly increasing along this $0$-string, which is therefore not a cellular loop.
\end{proof}

\medskip\remark Note that the same proof shows that if for every $k$-cellular string $F_1,\dots,F_m$ there is a $\pi$-coherent subdivision $\mathcal{F}$ of $\pi_{k+1}(P)$ by $(k+1)$-polytopes, such that $\pi_k(F_i) \in \mathcal{F}$ for all $1 \leq i \leq m$; then, $P$ is $k$-loop-free.

\SS\label{c:up-to-dim-3}

Up to dimension~$3$, all possible cases are covered by a combination of \cref{prop:dim0} and \cref{prop:dimP-2}, so we have the following.

\medskip\corollary Any framed polytope of dimension less than or equal to $3$ is loop-free.

%% file: sec/loops.tex

\subsection{Examples of cellular loops}\label{s:examples-loops}

In this section we give explicit examples of framed polytopes that admit cellular loops.
An accompanying \texttt{SageMath} notebook reproducing the constructions is available at \href{https://github.com/apadrol/polycat}{\texttt{github.com/apadrol/polycat}}.

\SS\label{ss:1loop} We will describe first a $5$-simplex $P_5$ for which the canonical frame is admissible and construct a cellular $1$-loop in the resulting framed polytope.
Consider the $6$ points $p_1,\dots,p_6$ in $\R^5$ whose coordinates are the columns of the following matrix:
\[
\begin{pmatrix}
-3 & -2 & -1 & 1 & 2 & 3 \\
-1 & 1 & 0 & 0 & 1 & -1 \\
-1 & 1 & 0 & 1 & -1 & 1 \\
0 & 0 & 1 & 1 & 0 & 0 \\
1 & 1 & 1 & 0 & 0 & 0
\end{pmatrix}.
\]
These points are affinely independent, and therefore their convex hull~$P_5$ is a $5$-simplex.
Furthermore, the canonical frame of $\R^5$ is $P_5$-admissible since the defining $6$ points are in general position and the projection of any $k$-face of $P_5$ by $\pi_k$ is itself a $k$-simplex.

\medskip\proposition The following sequence of $2$-faces of the canonically framed $5$-simplex $P_5$ forms a $1$-loop:
\begin{equation}\label{eq:1loop}
	[p_1p_2p_3],\ [p_2p_3p_6],\ [p_2p_4p_6],\ [p_4p_5p_6],\ [p_1p_4p_5],\ [p_1p_3p_5],\ [p_1p_2p_3].
\end{equation}

\begin{figure}[htpb]
	\centering
	\input{fig/loop_figure}
	\caption{A cellular 1-loop in $P_5$ formed by 2-faces. It depicts the image of the vertices of $P_5$ and some of its edges under the projection $\pi_2 \colon \R^5 \to \R^2$.}\label{fig:counterexample}
\end{figure}
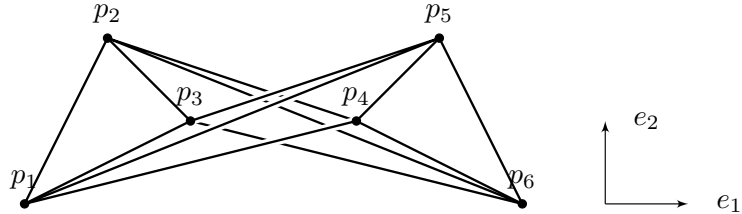

\begin{proof}
	Consulting \cref{fig:counterexample}, it is straightforward to check that for each of the triangles~$t_i$, the edge $t_i \cap t_{i+1}$ lies in the $1$-target $\ta_1(t_i)$ of $t_i$ and in the $1$-source~$\so_1(t_{i+1})$ of $t_{i+1}$.
	The detailed computation can be found in our \texttt{Sagemath} notebook.
\end{proof}

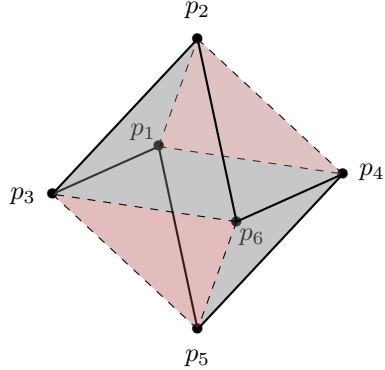
\begin{figure}[h!]
	\tdplotsetmaincoords{75}{75}
	\input{fig/cylinder}
	\caption{The loop (\ref{eq:1loop}) forms the boundary of an octahedron, where two opposite faces have been removed.}
	\label{fig:cylinder}
\end{figure}

Combinatorially, the faces involved in this loop form the boundary of an octahedron minus two opposite triangles, see \cref{fig:cylinder}.

All the $2$-faces involved in the loop \eqref{eq:1loop} are also faces of the $4$-dimensional polytope $P_4 \defeq \pi_4(P_5)$, and hence, with the canonical frame, $P_4$ has a cellular 1-loop.
For the interested reader we mention that $P_4$ is combinatorially equivalent to the \textit{cyclic polytope} $\CP(6,4)$, part of a family of polytopes that we will study in detail in later sections (see \cref{ss:cyclic_loops}).

\SS\label{ss:2loop} We now present a cellular $2$-loop on a framed $6$-simplex.
The loop appears on cells of the non-regular triangulation of Lee's twisted capped prism \cite[Sec.~6]{Lee1991} (see also \cite[Fig.~16.3.3]{LeeSantos2018}), for which we provide an explicit realization. The cyclic nature of this triangulation was already observed by Rambau \cite[Rmk.~5.10]{Rambau1997}. This is a relative of the so-called \emph{mother of all examples} \cite[Sec.~7.1]{deloeraTriangulations2010}.

In contrast with our previous example (\cref{ss:1loop}), the projections of the simplices involved in the loop do not overlap, and their vertices are in convex position.

Consider the $7$ points $q_0,q_1,\dots,q_6$ in $\R^6$ whose coordinates are given by the columns of the following matrix:
\[
\begin{pmatrix}
	0 & 10 & 0 & 0 & 7 & 2 & 3 \\
	0 & 0 & 10 & 0 & 3 & 7 & 2 \\
	0 & 0 & 0 & 10 & 2 & 3 & 7 \\
	1 & 1 & 1 & 0 & 1 & 0 & 0 \\
	0 & 0 & 0 & 1 & 1 & 0 & 1 \\
	0 & 0 & 1 & 0 & 0 & 1 & 0
\end{pmatrix}.
\]
They are affinely independent, and therefore form the vertex set of a simplex~$Q_6$.
As can be easily checked, the frame $B$ of~$\R^6$ given by the columns $v_1,\dots,v_6$ of the matrix
\[
\begin{pmatrix}
	-1 & 2 & 1 & 0 & 0 & 0 \\
	1 & 4 & 1 & 0 & 0 & 0 \\
	-1 & -1 & 1 & 0 & 0 & 0 \\
	0 & 0 & 0 & 1 & 0 & 0 \\
	0 & 0 & 0 & 1 & 1 & 0 \\
	0 & 0 & 0 & 1 & 1 & 1
\end{pmatrix}
\]
is $Q_6$-admissible.
We refer the interested reader to our \texttt{Sagemath} notebook for the detailed computation.

\medskip\proposition The following sequence of $3$-faces of the framed $6$-simplex $(Q_6,B)$ forms a $2$-loop:
\begin{multline} \label{eq:2loop}
	[q_0 q_1 q_4 q_5] ,\ [q_0 q_1 q_3 q_4] ,\ [q_0 q_3 q_4 q_6] ,\ [q_0 q_2 q_3 q_6] ,\ [q_0 q_2 q_5 q_6] ,\ [q_0 q_1 q_2 q_5] ,\ [q_0 q_1 q_4 q_5].
\end{multline}

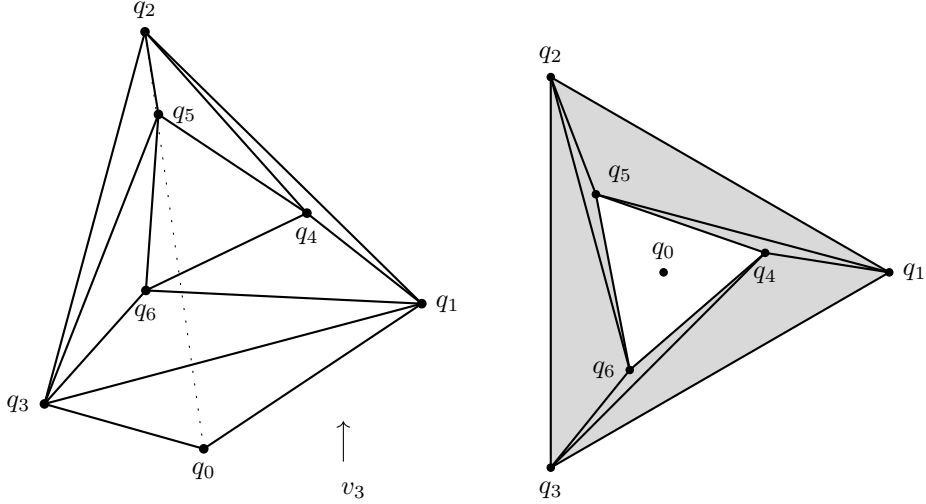
\begin{figure}[htpb]
	\centering
	\input{fig/loop2_figure3d}\qquad
	\input{fig/loop2_figure2d}
	\caption{A cellular $2$-loop on $Q_6$. The convex hull of $\pi_3(Q_6)$ is depicted on the left. On the right we see $\pi_2(Q_6)$ and the edges more relevant in the loop (note that they are not the same as in the convex hull).}
	\label{fig:counterexample2}
\end{figure}

\begin{proof}
	The claim can be verified directly through a straightforward computation.
	However, instead of presenting it completely, we provide the intuition behind it with the help of \cref{fig:counterexample2}.
	Since $\Lin(v_1, v_2, v_3) = \Lin(e_1,e_2,e_3)$ and $\Lin(v_4,v_5,v_6) = \Lin(e_4,e_5,e_6)$, the projection $\pi_3$ is given by forgetting the last three coordinates.
	This is a $2$-loop so it is apparent on $\pi_3(Q_6)$, which is depicted on the left of \cref{fig:counterexample2}.
	The vector~$v_3$, which determines the $2$-sources and $2$-targets, goes in the direction from $q_0$ to the center of the equilateral triangle $[q_4 q_5 q_6]$.
	As the points $q_1,\dots,q_6$ are very close to being coplanar, it is somehow easier to understand the loop in the $2$-dimensional picture on the right of \cref{fig:counterexample2}.
	Here, one should think of the point $q_0$ as being behind the ``plane'' spanned by the other points, and to $v_3$ as being perpendicular to this plane.

	The loop consists of the six tetrahedra arising as the cone of $q_0$ over each of the shaded triangles in the picture on the right.
	For each tetrahedron, the facets (triangles) pointing ``downwards'', towards $q_0$ are in the source, while those pointing ``upwards'', away from $q_0$ are in the target.
	For example, for the tetrahedron $[q_0 q_1 q_4 q_5]$, the source is the triangle $[q_0 q_1 q_5]$, and the target is formed by the triangles $[q_1 q_4 q_5]$, $[q_0 q_1 q_4]$, and $[q_0 q_4 q_5]$.
	Similarly, for the tetrahedron $[q_0 q_1 q_3 q_4]$, the source are the triangles $[q_0, q_1, q_4]$ and $[q_0, q_1, q_3]$, and the target are the triangles $[q_0, q_3, q_4]$ and $[q_1, q_3, q_4]$.
	A similar analysis holds for the other tetrahedra.
	Then, we see that the triangle $[q_0 q_1 q_4]$ is in the target of $[q_0 q_1 q_4 q_5]$ and in the source of $[q_0 q_1 q_3 q_4]$; the triangle $[q_0 q_3 q_4]$ is in the target of $[q_0 q_1 q_3 q_4]$ and in the source of $[q_0 q_3 q_4 q_6]$; and so on, enabling the entire sequence to form the loop (\ref{eq:2loop}).
\end{proof}

%% file: fig/loop_figure.tex
\resizebox{0.8\linewidth}{!}{
\begin{tikzpicture}[auto, node distance = 2cm, > = latex',
	point/.style = {draw, circle, fill = black, inner sep = 1pt}]

	\coordinate (a) at (-3,-1);
	\coordinate (b) at (-2,1);
	\coordinate (c) at (-1,0);
	\coordinate (d) at (1,0);
	\coordinate (e) at (2,1);
	\coordinate (f) at (3,-1);

	\draw[thick] (c) -- (f) -- (b) -- (d);
	\draw[line width = 3pt,white] (d) -- (a) -- (e) -- (c);
	\draw[thick] (d) -- (a) -- (e) -- (c);
	\draw[thick] (a) -- (b) -- (c) -- (a);
	\draw[thick] (d) -- (e) -- (f) -- (d);

	\foreach \i/\l in {a/p_1,b/p_2,c/p_3,d/p_4,e/p_5,f/p_6}
	{
		\node[point,label = {above:$\l$}] at (\i) {};
	}

	\draw[->] (4,-1)-- (4,0);
	\draw[->] (4,-1)-- (5,-1);
	\node[] at (5.5, -1) {$e_1$};
	\node[] at (4.5, 0) {$e_2$};
\end{tikzpicture}}

%% file: fig/cylinder.tex
\begin{tikzpicture}
	[scale = 2,
	tdplot_main_coords,
	cube/.style = {very thick,black},
	grid/.style = {very thin,gray},
	axis/.style = {->}]

	\node (a) at (-1,0,0) {$\bullet$};
	\node at (-1.4,0,0) {$p_1$};

	\node (b) at (0,0,1) {$\bullet$};
	\node at (0,0,1.2) {$p_2$};

	\node (c) at (0,-1,0) {$\bullet$};
	\node at (0,-1.2,0) {$p_3$};

	\node (f) at (0,1,0) {$\bullet$};
	\node at (0,1.2,0) {$p_4$};

	\node (e) at (0,0,-1) {$\bullet$};
	\node at (0,0,-1.2) {$p_5$};

	\node (d) at (1,0,0) {$\bullet$};
	\node at (1.4,0,0) {$p_6$};

	\draw[-,dashed] (-1,0,0)--(0,1,0);
	\draw[-,thick] (-1,0,0)--(0,-1,0);
	\draw[-,dashed] (-1,0,0)--(0,0,1);
	\draw[-,thick] (-1,0,0)--(0,0,-1);


	\fill[red, opacity = 0.15] (-1,0,0) -- (0,0,1) -- (0,1,0) -- (-1,0,0);
	\fill[black!50, opacity = 0.25] (-1,0,0) -- (0,0,-1) -- (0,1,0) -- (-1,0,0);
	\fill[black!50, opacity = 0.25] (-1,0,0) -- (0,0,-1) -- (0,-1,0) -- (-1,0,0);
	\fill[black!50, opacity = 0.25] (-1,0,0) -- (0,0,1) -- (0,-1,0) -- (-1,0,0);
	\fill[black!50, opacity = 0.25] (1,0,0) -- (0,0,1) -- (0,1,0) -- (1,0,0);
	\fill[black!50, opacity = 0.25] (1,0,0) -- (0,0,-1) -- (0,1,0) -- (1,0,0);
	\fill[red, opacity = 0.15] (1,0,0) -- (0,0,-1) -- (0,-1,0) -- (1,0,0);
	\fill[black!50, opacity = 0.25] (1,0,0) -- (0,0,1) -- (0,-1,0) -- (1,0,0);
	\draw[-,thick] (1,0,0)--(0,1,0);
	\draw[-,dashed] (1,0,0)--(0,-1,0);
	\draw[-,thick] (1,0,0)--(0,0,1);
	\draw[-,dashed] (1,0,0)--(0,0,-1);
	\draw[-,thick] (0,1,0)--(0,0,-1);
	\draw[-,dashed] (0,-1,0)--(0,0,-1);
	\draw[-,dashed] (0,1,0)--(0,0,1);
	\draw[-,thick] (0,-1,0)--(0,0,1);
\end{tikzpicture}

%% file: fig/loop2_figure3d.tex
\begin{tikzpicture}%
	[
x = {(0.966340cm, -0.190271cm)},
	y = {(0.257268cm, 0.714691cm)},
	z = {(-0.000001cm, 0.673063cm)},
	xscale = .300000,
	yscale = .400000,
	back/.style = {loosely dotted, thin},
	edge/.style = {color = black, thick},
	facet/.style = {draw},
	vertex/.style = {inner sep = 1pt,circle,draw = black,fill = black,thick}]
%
%

\coordinate (-3.00000, 3.46400, 12.00000) at (-3.00000, 3.46400, 12.00000);
\coordinate (-5.00000, -8.66000, 10.00000) at (-5.00000, -8.66000, 10.00000);
\coordinate (-5.00000, 8.66000, 10.00000) at (-5.00000, 8.66000, 10.00000);
\coordinate (-1.50000, -4.33000, 12.00000) at (-1.50000, -4.33000, 12.00000);
\coordinate (0.00000, 0.00000, 0.00000) at (0.00000, 0.00000, 0.00000);
\coordinate (10.00000, 0.00000, 10.00000) at (10.00000, 0.00000, 10.00000);
\coordinate (4.50000, 0.86600, 12.00000) at (4.50000, 0.86600, 12.00000);
\draw[edge,back] (-5.00000, 8.66000, 10.00000) -- (0.00000, 0.00000, 0.00000);
\draw[edge] (-3.00000, 3.46400, 12.00000) -- (-5.00000, -8.66000, 10.00000);
\draw[edge] (-3.00000, 3.46400, 12.00000) -- (-5.00000, 8.66000, 10.00000);
\draw[edge] (-3.00000, 3.46400, 12.00000) -- (-1.50000, -4.33000, 12.00000);
\draw[edge] (-3.00000, 3.46400, 12.00000) -- (4.50000, 0.86600, 12.00000);
\draw[edge] (-5.00000, -8.66000, 10.00000) -- (-5.00000, 8.66000, 10.00000);
\draw[edge] (-5.00000, -8.66000, 10.00000) -- (-1.50000, -4.33000, 12.00000);
\draw[edge] (-5.00000, -8.66000, 10.00000) -- (0.00000, 0.00000, 0.00000);
\draw[edge] (-5.00000, -8.66000, 10.00000) -- (10.00000, 0.00000, 10.00000);
\draw[edge] (-5.00000, 8.66000, 10.00000) -- (10.00000, 0.00000, 10.00000);
\draw[edge] (-5.00000, 8.66000, 10.00000) -- (4.50000, 0.86600, 12.00000);
\draw[edge] (-1.50000, -4.33000, 12.00000) -- (10.00000, 0.00000, 10.00000);
\draw[edge] (-1.50000, -4.33000, 12.00000) -- (4.50000, 0.86600, 12.00000);
\draw[edge] (0.00000, 0.00000, 0.00000) -- (10.00000, 0.00000, 10.00000);
\draw[edge] (10.00000, 0.00000, 10.00000) -- (4.50000, 0.86600, 12.00000);

\node[vertex, label = {right:$q_5$}] at (-3.00000, 3.46400, 12.00000) {};
\node[vertex, label = {left:$q_3$}] at (-5.00000, -8.66000, 10.00000) {};
\node[vertex, label = {above:$q_2$}] at (-5.00000, 8.66000, 10.00000) {};
\node[vertex, label = {below:$q_6$}] at (-1.50000, -4.33000, 12.00000) {};
\node[vertex, label = {below:$q_0$}] at (0.00000, 0.00000, 0.00000) {};
\node[vertex, label = {right:$q_1$}] at (10.00000, 0.00000, 10.00000) {};
\node[vertex, label = {below:$q_4$}] at (4.50000, 0.86600, 12.00000) {};

	\draw[->] (8,-6,8)-- (8,-6,10);
	\node[] at (9,-8,9) {$v_3$};
\end{tikzpicture}

%% file: fig/loop2_figure2d.tex
\begin{tikzpicture}[auto, node distance = 2cm, > = latex', scale = .3,
	point/.style = {draw, circle, fill = black, inner sep = 1pt}]

	\coordinate (p0) at (0.0000, 0.0000);
	\coordinate (p1) at (10.00, 0.0000);
	\coordinate (p2) at (-5.000, 8.660);
	\coordinate (p3) at (-5.000, -8.660);
	\coordinate (p4) at (4.500, 0.8660);
	\coordinate (p5) at (-3.000, 3.464);
	\coordinate (p6) at (-1.500, -4.330);

	\fill[black!15] (p1) -- (p4) -- (p5) -- (p1);
	\fill[black!15] (p1) -- (p4) -- (p3) -- (p1);
	\fill[black!15] (p6) -- (p4) -- (p3) -- (p6);
	\fill[black!15] (p6) -- (p2) -- (p3) -- (p6);
	\fill[black!15] (p6) -- (p2) -- (p5) -- (p6);
	\fill[black!15] (p1) -- (p2) -- (p5) -- (p1);

	\draw[thick] (p1) -- (p2) -- (p3) -- (p1);
	\draw[thick] (p4) -- (p5) -- (p6) -- (p4);
	\draw[thick] (p1) -- (p5)-- (p2) -- (p6)-- (p3) -- (p4) -- (p1);

%

	\node[point,label = {above:$q_0$}] at (p0) {};
	\node[point,label = {right:$q_1$}] at (p1) {};
	\node[point,label = {above:$q_2$}] at (p2) {};
	\node[point,label = {below:$q_3$}] at (p3) {};
	\node[point,label = {below:$q_4$}] at (p4) {};
	\node[point,label = {above right:$q_5$}] at (p5) {};
	\node[point,label = {left:$q_6$}] at (p6) {};

\end{tikzpicture}

%% file: sec/bootstrapping.tex

\subsection{Bootstrapping loops}\label{sec:bootstrap}

In this section we use our previous examples of framed polytopes with loops to construct new ones.
In particular, we do so for any simplicial or simple polytope of dimension at least \(6\), with a better dimensional bound for cyclic polytopes.

\SS\label{cor:simplicial}
\lemma Let \(P \subset \R^d\) be a polytope and \(F\) one of its \(\ell\)-dimensional faces.
Given an \(F\)-admissible frame \(B = (v_1,\dots,v_\ell)\), there is a \(P\)-admissible frame \(B' = (v'_1,\dots,v'_d)\) such that \((v_1,\dots,v_\ell)\) and \((v'_1,\dots,v'_\ell)\) are \(F\)-equivalent.

\begin{proof}
	This is due to the fact that \(P\)-equivalence is an open condition and that \(F\)-equivalence is preserved by lower triangular transformations (\cref{ss:lowertriangular}).
\end{proof}

Using the notation of the previous lemma, we have that if \((F,B)\) has a loop, then the same loop will exist in \((P,B')\). Therefore, by the example in \cref{ss:1loop}, we obtain the following.

\medskip\theorem
Any polytope containing a 5-dimensional simplex as a face admits a frame inducing a loop.

\medskip Polytopes with the property that their facets are all simplices are referred to as \defn{simplicial polytopes}.
The following is immediate from the previous result.

\medskip\corollary
Every simplicial polytope \(P\) of dimension at least \(6\) admits a frame inducing a loop.

\SS\label{ss:cyclic_loops}

An important family of simplicial polytopes is given by the cyclic polytopes, which we now define.
The \defn{moment curve} \(v \colon \R \to \R^d\) is defined by \(v(t) \defeq (t,t^2,\dots,t^d)\).
For any set of \(n\) distinct real numbers \(t_1 < \cdots < t_n\) with \(n > d\), we define the \defn{cyclic polytope} \(\CP(n,d)\) as the convex hull of the points \(\set{v(t_1),\dots,v(t_n)}\).
The combinatorial type of this polytope, and its behavior with respect to the canonical frame (see \cref{l:facets-cyclic-polytope}), does not depend on the exact choice of \(t_1 < \dots < t_n\), which motivates the use of the notation \(\CP(n,d)\) without specifying the parameters \(t_i\).
Further details can be found, for example, in \cite[Ex.~0.6]{Ziegler95}.
For any \(d \in \N\), the polytope \(\CP(d+1,d)\) is a \(d\)-simplex, and we will refer to it as the \defn{cyclic \(d\)-simplex}.
We will see in \cref{ss:simplex} that, when combined with the canonical frame, these simplices give rise to Street's orientals, one of the most prominent examples.
In particular, canonically framed cyclic polytopes are always loop-free.
The following result shows how the geometry plays a central role, as there are polytopes with the same combinatorial type that do have loops.

\medskip\theorem
For any \(d \geq 4\) and \(n \geq \max(6, d+1)\), there are polytopes combinatorially equivalent to~\(\CP(n,d)\) for which the canonical frame induces a loop.

\begin{proof}
	For \(d \geq 6\), this follows directly from \cref{cor:simplicial}, showing that all simplicial polytopes — of which cyclic polytopes are examples — admit a loop-inducing frame if their dimension is at least \(6\).
	The linear map sending that frame to the canonical frame gives a combinatorial cyclic polytope for which the canonical frame induces a loop (see \cref{ss:invariancelinearautomorphism}).
	Therefore, we need a particular argument for the remaining cases: \(\CP(n,4)\) and \(\CP(n,5)\).

	Choose \(t_1<t_2<\cdots<t_n\) in \(\R\) and label the vertices of the cyclic polytope \(\CP(n,d)\) by \(a_i \defeq (t_i,t_i^2,\dots,t_i^d)\), for \(1 \leq i \leq n\).
	The 5-polytope \(P_5\) from~\cref{ss:1loop} is a simplex, and any ordering of its vertices induces a combinatorial isomorphism with \(\CP(6,5)\), which is also a simplex.
	The 4-polytope \(P_4\) from \cref{ss:1loop} is combinatorially equivalent to a direct sum of two triangles: one with vertices labeled by \(p_1, p_2, p_4\), and one with vertices labeled by \(p_3, p_5, p_6\).
	The latter triangle is combinatorially equivalent to the cyclic polytope \(\CP(6,4)\).
	For simplicity, we relabel the vertices of \(P_4\) as follows:
	\[
	p_1 \mapsto \tilde p_1,\
	p_2 \mapsto \tilde p_3,\
	p_3 \mapsto \tilde p_2,\
	p_4 \mapsto \tilde p_5,\
	p_5 \mapsto \tilde p_4,\
	p_6 \mapsto \tilde p_6.
	\]
	Then, \(P_4\) is combinatorially equivalent to the direct sum of the two triangles \(\tilde p_1, \tilde p_3, \tilde p_5\) and \(\tilde p_2, \tilde p_4, \tilde p_6\), and the combinatorial isomorphism with \(\CP(6,4)\) can be realized by the map \(\tilde p_i \mapsto a_i\).

	Now, given a polytope combinatorially equivalent to \(\CP(n,d)\), it is always possible to construct a \emph{cyclic extension} adding an additional vertex to obtain a polytope combinatorially equivalent to \(\CP(n+1,d)\); cf.\ \cite[Thm.~5.1]{LeeMenzel2010} and \cite[Prop.~3.16]{Padrol2013} for details.
	More precisely, if \(\conv(\tilde p_1,\dots, \tilde p_n)\) is combinatorially equivalent to \(\CP(n,d)\) with the vertices cyclically ordered, then there is some \(\varepsilon>0\) such that for
	\[
	\tilde p_{n+1} \defeq \tilde p_n - \varepsilon (\tilde p_{n-1}- \tilde p_n) + \cdots + (-\varepsilon)^{d} (\tilde p_{n-d}- \tilde p_{n-d+1})
	\]
	we have that \(\conv(\tilde p_1,\dots, \tilde p_n, \tilde p_{n+1})\) is combinatorially equivalent to \(\CP(n+1,d)\).
	Repeating this operation \(n-6\) times on \(P_4\), we obtain a combinatorial cyclic 4-polytope with \(n\) vertices, of which the first six are \(\tilde p_i\) with \(1 \leq i \leq 6\).
	We can check with \emph{Gale's evenness condition} (see \cite[Thm~0.7]{Ziegler95}) that the faces \([\tilde p_1\tilde p_2\tilde p_3\tilde p_4]\), \([\tilde p_1\tilde p_2\tilde p_4\tilde p_5]\), \([\tilde p_2\tilde p_3\tilde p_5\tilde p_6]\), and \([\tilde p_3\tilde p_4\tilde p_5\tilde p_6]\) are still facets of \(\CP(n,4)\), and therefore the triangles involved in the loop~(\ref{eq:1loop}) are also faces of this polytope.
	An analogous argument starting with \(P_5\) instead of \(P_4\) gives combinatorial cyclic 5-polytopes with 1-loops for all \(n \geq 6\).
\end{proof}

\SS\label{ss:hyperplane}
We now study the relationship between loops and hyperplane intersections.
We say that an affine $(d-k)$-dimensional subspace~$L$ is generic (with respect to~$P$), whenever its intersection with every $e$-dimensional face of~$P$ is either empty or $(e-k)$-dimensional.
In particular, a hyperplane $H$ is generic if it does not contain any vertex of~$P$.

\medskip\lemma
Let $P$ be a (strongly) loop-free framed polytope with frame $(v_1,\ldots,v_d)$.
Then, the intersection of~$P$ with any generic hyperplane parallel to $\ker\pi_1 = \Lin(v_2,\ldots,v_d)$ is framed by $(v_2,\ldots,v_d)$ and (strongly) loop-free.

\begin{proof}
	Let $H$ be a generic hyperplane of $\R^d$ parallel to $\ker\pi_1$, and let $Q \defeq P \cap H$ denote its intersection with $P$.
	By the genericity condition, each $k$-face $F_i$ of $Q$ is of the form $F_i = G_i \cap H$ for some $(k+1)$-face $G_i$ of $P$.
	Observe first that by the definition of the frame on~$Q$, the system of projections $\pi_k^Q$ of $Q$ is inherited from the system of projections $\pi_k^P$ of $P$, and that the former is just the restriction of the latter to the hyperplane~$H$.
	That is, we have $\pi_k^Q = \pi_{k+1}^P \big |_H$.
	Since for each $(k+1)$-face $G$ of $P$, the projection $\pi_{k+1}^P$ induces a linear isomorphism~$\Lin G \cong V_{k+1}^P = \Lin(v_1,\dots,v_{k+1})$, the projection $\pi_k^Q$ induces a linear isomorphism~$\Lin (G \cap H) \cong V_{k+1}^P \cap \ker\pi_1 = V_k^Q = \Lin(v_2,\dots,v_k)$ for each $k$-face $F = G \cap H$ of $Q$.
	Thus, the frame $(v_2,\ldots,v_d)$ is~$Q$-admissible.

	Now let $Q'$ be a $k$-face of $Q$ and $P'$ be the $(k+1)$-face of $P$ such that $Q'=P' \cap H$; let $F$ be a facet of $Q'$ and $G$ be the $k$-face of $P$ such that $F = G \cap H$.
	According to \cref{lem:sources and targets in a framed polytope}, to decide whether $F$ is in the source or the target of $Q'$, we need to compute the sign of \(\sprod{\normal^{\pi_{k}(F)}_{\pi_{k}(Q')}}{v_{k}}\).
	Note that $\normal^{\pi_{k}(G)}_{\pi_{k}(P')}$ is also a normal covector to $\pi_{k}(F)$ in $\pi_{k}(Q')$, because $\pi_{k}(F)=\pi_{k}(G)\cap \pi_{k}(H)$ and $\pi_{k}(Q')=\pi_{k}(P')\cap \pi_{k}(H)$.
	Hence, $F$ is in the source (resp.\ target) of $Q'$ if and only if $G$ is in the source (resp.\ target) of~$P'$.

	Therefore, any inhomogeneous cellular loop $F_1,\ldots, F_l$ in $Q$ induces an inhomogeneous cellular loop $G_1,\ldots,G_l$ in $P$.
	The statement for cellular loops then follows from \cref{l:inhomogeneous strings}.
\end{proof}

A converse of the previous result is the following.

\medskip\theorem
Let $P$ be a $d$-polytope and let $L$ be a generic $\ell$-dimensional subspace.
If there is a $(P \cap L)$-admissible frame inducing a $k$-loop on $P \cap L$, then there is a $P$-admissible frame inducing a $k+(d-\ell)$-loop on $P$.

\begin{proof}
	Let us write $L$ as an intersection of generic hyperplanes $L = H_1 \cap \cdots \cap H_{d-\ell}$, and extend the frame $B = (v_1,\dots,v_\ell)$ to a basis $\hat B = (v_1,\dots,v_d)$ of $V_d$ by taking for each $1 \leq i \leq d-\ell$ a vector $v_{\ell+i}$ in $\left(\bigcap_{j = i+1}^{d-\ell} H_j\right)\setminus \left(\bigcap_{j = i}^{d-\ell} H_j\right)$.
	The proof of \cref{ss:hyperplane} shows that a $k$-loop in $P \cap \left(\bigcap_{j = i}^{d-\ell} H_j\right)$ induces a $(k+1)$-loop in $P \cap \left(\bigcap_{j = i+1}^{d-\ell} H_j\right)$.
	By induction we get that it induces a $k+(d-\ell)$-loop in~$P$.
\end{proof}

\SS\label{ss:vertex_figure}
The intersection of a polytope \(P\) with a hyperplane \(H\) separating a vertex~\(v\) from all the other vertices is a polytope called a \defn{vertex figure} of \(P\) at \(v\).
While the construction depends on the choice of \(H\), the combinatorial type does not.
Its face lattice is the interval of \(\faces\) between \(v\) and \(P\); see \cite[Lec.~2]{Ziegler95}.
A direct application of the previous result gives the following.

\medskip\corollary\label{previous-result}
If a vertex figure of a polytope~$P$ admits a frame inducing a $k$-loop, then~$P$ admits a frame inducing a $(k+1)$-loop.

\medskip\defn{Simple polytopes} are defined by the property that their vertex figures are all simplices.
Simplices, cubes, associahedra, and permutohedra of any dimension are all examples of simple polytopes.

\medskip\corollary Every simple polytope of dimension at least $6$ admits a frame inducing a loop.

\begin{proof}
	Any simple polytope $P$ of dimension at least \(6\) has a $5$-dimensional iterated vertex figure $F$ which is a simplex.
	Using the unique affine isomorphism between the canonically framed $5$-simplex~$P_5$ from~\cref{ss:1loop} and $F$, we can use \cref{ss:invariancelinearautomorphism} to obtain a frame inducing a $1$-loop on~$F$.
	Then, we conclude with the previous result \cref{previous-result}.
\end{proof}

%% file: sec/noframe.tex

\subsection{Loop inevitability}\label{s:loop-inevitability}

The goal of this section is to construct polytopes for which \emph{every} admissible frame induces a cellular loop (\cref{cor:noframe}).

\begin{figure}[h!]
	\centering
	\input{fig/flattening}
	\caption{A regular hexagon $P$ for which the frames $(v_1,v_2)$ and $(v_1,v_2')$ are not $P$-equivalent, and a flattened version $\Phi_\varepsilon(P)$ for which they are.
		The set of vectors $w$ for which $(v_1,w)$ is $P$-equivalent to $(v_1,v_2)$ is depicted with a blue cone, and the set of those that are $\Phi_\varepsilon(P)$-equivalent is depicted as a larger red cone.}
	\label{fig:flattening}
\end{figure}
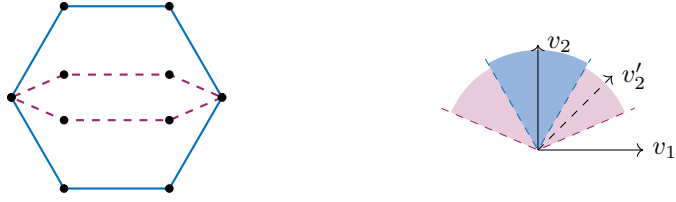

\begin{figure}[h!]
	\input{fig/squashing}
	\caption{A rhombus $P$ is squashed on top of the facet~$F$ of a rectangle~$Q$.}
	\label{fig:squashing}
\end{figure}
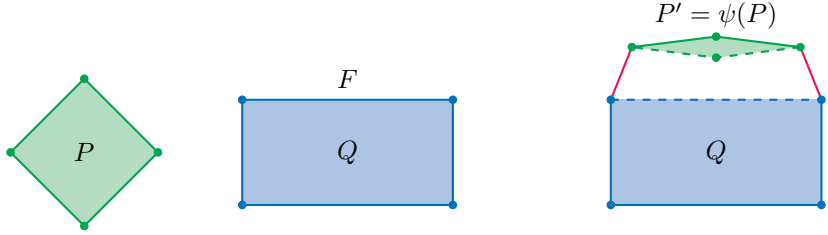

The main idea is to transform framed polytopes with a loop via an operation called \defn{flattening} (\cref{ss:flattening}, see \cref{fig:flattening}) that expands the space of loop-inducing frames.
Since $P$-admissibility is an open condition, if $B$ is a frame inducing a loop on~$P$, then there is an open neighborhood of~$B$ which contains frames that are $P$-equivalent to $B$.
Flattening~$P$ makes this neighborhood large enough to almost cover the space of $P$-admissible frames.
One then combines several reflected copies of a flattened polytope in order to completely cover the full space of admissible frames by loop-inducing ones.
This is done via an operation called \defn{squashing} (\cref{def:squashing}, see \cref{fig:squashing}), which is used to guarantee the convex position of the loop-inducing faces of the copies.
The whole construction, which is quite technical, is decomposed into several independent lemmas.
The first one defines a family of vectors~$\tilde v_i'$ which will be used often in the sequel.

In this section, we consider $\R^d$ as a Euclidean space endowed with the standard inner product.
We continue to use the notation $\sprod{\cdot}{\cdot}$ for the inner product, identifying vectors and covectors via the canonical isomorphism with the dual space induced by the inner product.

\SS\label{lem:directsum}

Let $B = (v_1,\dots,v_d)$ and $B' = (v_1',\dots,v'_d)$ be two orthonormal frames of~$\R^d$. We denote the subspaces and the system of projections associated to $B$ (resp. $B'$) by $V_i$ and $\pi_i$ (resp. $V_i'$ and $\pi_i'$), as in \cref{sec:frame}. Recall in particular that
$V_i = \Lin(v_1,\dots,v_i)$ is the subspace spanned by the first~$i$ vectors of $B$ and that $\ker(\pi_i')= \Lin(v_{i+1}',\dots,v_d')$ is the subspace spanned by the last $d-i$ vectors of $B'$.

\medskip\lemma
Suppose that $C = (v_1,\dots,v_i,v_{i+1}',\dots,v_d')$ forms a basis of $\R^d$ such that the projection~$\tilde v_i'$ of~$v_i'$ onto~$V_i$ along $\ker(\pi_i')$ is well defined.
Then, if ${\sprod{v_i}{\tilde v_i'}} > 0$ we have that $C' = (v_1,\dots,v_{i-1},v_i',v_{i+1}',\dots,v_d')$ forms a basis of $\R^d$ too.

\begin{proof}
	Observe that $\tilde v_i'\in V_i$ does not belong to $\Lin(v_1,\dots,v_{i-1})$ because $B$ is orthonormal and ${\sprod{v_i}{\tilde v_i'}} > 0$.
	Thus we have that $V_i = \Lin(v_1,\dots,v_{i-1},\tilde v_i')$, and the family $(v_1,\dots,v_{i-1},\tilde v_i',v_{i+1}',\dots,v_d')$ forms a basis.
	Since $\tilde v_i'$ and $v_i'$ differ by an element of $\ker(\pi_i')$, it follows that $C'$ forms a basis as well.
\end{proof}

\SS\label{ss:flattening}
We now discuss the flattening operation.
Given a frame $B = (v_1,\dots, v_d)$ and a small positive real tuple $\varepsilon = (\varepsilon_1,\dots,\varepsilon_d)$ in $\R^d$, let us define the linear automorphism $\Phi_\varepsilon \colon v_i \mapsto \varepsilon_i v_i$ which scales the $v_i$ coordinate by~$\varepsilon_i$.
A \defn{flattening} of a framed polytope~$(P,B)$ is the polytope~$\Phi_\varepsilon(P)$ obtained as the image of $P$ under the map $\Phi_\varepsilon$ defined by~$B$ and some $\varepsilon \in \R^d$.

The interest of flattening a framed polytope lies in the following result, illustrated in \cref{fig:flattening}.
Here, the vector~$\tilde v_i'$ is defined as in \cref{lem:directsum} and moreover rescaled to have length $1$, and the map $\Phi_\varepsilon$ is taken with respect to the frame $B$.

\medskip\theorem
Let $P \subset \R^d$ be framed by $B = (v_1,\dots, v_d)$ assumed orthonormal.
For every $0 < \delta < 1$ there is a positive real tuple $\varepsilon \in \R_{>0}^d$ such that if $B' = (v_1',\dots,v'_d)$ is an orthonormal frame of $\R^d$ with ${\sprod{v_i}{\tilde v_i'}} > \delta$ for all $1 \leq i \leq d$, then the frames $B$ and $B'$ are $\Phi_\varepsilon(P)$-equivalent.

\medskip\example Before proceeding to the proof, let us illustrate it with an example.
\cref{fig:flattening} represents a regular hexagon $P$, for which the frame $(v_1,v_2)\defeq ((1,0),(0,1))$ and the frame $(v_1',v_2')\defeq ((1,0),(1,1))$ induce distinct $f$-orientations.
However, for~$\varepsilon \defeq (1,\frac 1 4) \in \R^2$ both frames are $\Phi_\varepsilon(P)$-equivalent.

\begin{proof}
	Recall from \cref{ss:explicitdescription} that the $k$-sources and $k$-targets of the framed polytope~$(P,B)$ are determined by the signs of determinants of the form:
	\begin{equation}\label{eq:det1noframe}
		\begin{vmatrix}
			1 & \cdots &1&1&0&\cdots&0\\
			(p_1)_{1} & \cdots &(p_{k+1})_{1}&q_1&(v_{k+2})_1&\cdots&(v_{d})_1\\
			(p_1)_{2} & \cdots &(p_{k+1})_{2}&q_2&(v_{k+2})_2&\cdots&(v_{d})_2\\
			\vdots & \ddots &\vdots&\vdots&\vdots&\ddots&\vdots\\
			(p_1)_{d} &\cdots &(p_{k+1})_{d}&q_d&(v_{k+2})_d&\cdots&(v_{d})_d\\
		\end{vmatrix}
	\end{equation}
	and
	\begin{equation}\label{eq:det2noframe}
		\begin{vmatrix}
			1 & \cdots &1&0&0&\cdots&0\\
			(p_1)_{1} & \cdots &(p_{k+1})_{1}&(v_{k+1})_1&(v_{k+2})_1&\cdots&(v_{d})_1\\
			(p_1)_{2} & \cdots &(p_{k+1})_{2}&(v_{k+1})_2&(v_{k+2})_2&\cdots&(v_{d})_2\\
			\vdots & \ddots &\vdots&\vdots&\vdots&\ddots&\vdots\\
			(p_1)_{d} & \cdots &(p_{k+1})_{d}&(v_{k+1})_d&(v_{k+2})_d&\cdots&(v_{d})_d\\
		\end{vmatrix},
	\end{equation}
	where $p_1,\dots,p_{k+1}$ are affinely independent vertices in a $k$-face~$F$ of~$P$, and $q$ is another vertex of $P$.
	In what follows, we will consider coordinates expressed in the frame $B$.
	Note that, in this case, the lower right part of the matrices becomes an identity and the determinants \eqref{eq:det1noframe} and \eqref{eq:det2noframe} are reduced to
	\[
	\begin{vmatrix}
		1 & \cdots &1&1\\
		(p_1)_{1} & \cdots &(p_{k+1})_{1}&q_1\\
		(p_1)_{2} & \cdots &(p_{k+1})_{2}&q_2\\
		\vdots & \ddots &\vdots&\vdots\\
		(p_1)_{k+1} & \cdots &(p_{k+1})_{k+1}&q_{k+1}\\
	\end{vmatrix}
	\quad \text{and} \quad
	\begin{vmatrix}
		1 & \cdots &1&0\\
		(p_1)_{1} & \cdots &(p_{k+1})_{1}&(v_{k+1})_1\\
		(p_1)_{2} & \cdots &(p_{k+1})_{2}&(v_{k+1})_2\\
		\vdots & \ddots &\vdots&\vdots\\
		(p_1)_{k+1} & \cdots &(p_{k+1})_{k+1}&(v_{k+1})_{k+1}\\
	\end{vmatrix}.
	\]
	If we evaluate the same determinants for $(\Phi_\varepsilon(P),B)$ instead of $(P,B)$, the signs do not change, as they only differ by a product of $\varepsilon_i$ factors, which are positive.

	We now proceed to the evaluation of the determinants \eqref{eq:det1noframe} and \eqref{eq:det2noframe} for the framed polytope $(\Phi_\varepsilon(P),B')$.
	Gathering the vectors $\tilde v_i'$ from \cref{lem:directsum} into a frame $\tilde B' \defeq (\tilde v_1',\ldots, \tilde v_d')$, and using \cref{ss:lowertriangular}, we have that the determinants for $(\Phi_\varepsilon(P),\tilde B')$ coincide with the those for $(\Phi_\varepsilon(P),B')$.
	By the definition of $\tilde B'$, the last $i$ coordinates (in the basis~$B$) of $\tilde v_{d-i}'$ are~$0$, and thus the lower right submatrices in the determinants above are upper-triangular.

	Therefore, if we develop \eqref{eq:det1noframe} for $(\Phi_\varepsilon(P),\tilde B')$ with respect to the last $d-k-1$ columns, there is only one term that does not contain any $\varepsilon_{i}$ with $i\geq k+2$, which is
	\begin{equation}\label{eq:termdet1noframe}
		(\tilde v_{d}')_d\cdots(\tilde v_{k+2}')_{k+2}
		\varepsilon_{k+1}\cdots \varepsilon_1
		\begin{vmatrix}
			1 & \cdots &1&1\\
			(p_1)_{1} & \cdots &(p_{k+1})_{1}& q_1\\
			(p_1)_{2} & \cdots &(p_{k+1})_{2}&q_2\\
			\vdots & \ddots &\vdots&\vdots\\
			(p_1)_{k+1} & \cdots &(p_{k+1})_{k+1}&q_{k+1}
		\end{vmatrix}.
	\end{equation}
	Similarly, if we develop \eqref{eq:det2noframe} for $(\Phi_\varepsilon(P),\tilde B')$ with respect to the last $d-k-2$ columns, the only term without $\varepsilon_{i}$ with $i\geq k+1$ is
	\begin{equation}\label{eq:termdet2noframe}
		(\tilde v_{d}')_d(\tilde v_{d-1}')_{d-1}\cdots(\tilde v_{k+2}')_{k+2}(\tilde v_{k+1}')_{k+1}\varepsilon_{k}\cdots \varepsilon_1
		\begin{vmatrix}
			1 & \cdots &1\\
			(p_1)_{1} & \cdots &(p_{k+1})_{1}\\
			(p_1)_{2} & \cdots &(p_{k+1})_{2}\\
			\vdots & \ddots &\vdots\\
			(p_1)_{k} & \cdots &(p_{k+1})_{k}\\
		\end{vmatrix}.
	\end{equation}
	Now, our conditions imply that the coordinates of $\tilde v_i'$, expressed in the basis $B$, verify
	\[ \bars{(\tilde v_i')_j} =
	\begin{cases}
		0&\text{ for }j>i,\\
		> \delta &\text{ for }j = i,\\
		<1&\text{ for }j<i.
	\end{cases}
	\]
	Notice that the coordinates of $P$ are fixed, the absolute values of the off-diagonal coordinates of $\tilde B'$ are uniformly bounded above by~$1$, and the diagonal entries are bounded below by~$\delta$. Therefore, we can choose some $0<\varepsilon_d\ll \varepsilon_{d-1}\ll\cdots \ll \varepsilon_1$ small enough so that the sign of the determinants is dominated by the sign of the terms
	\eqref{eq:termdet1noframe} and 	\eqref{eq:termdet2noframe}, respectively.

	Up to the respective factors $(\tilde v_{d}')_d\cdots(\tilde v_{k+2}')_{k+2}$ and $(\tilde v_{d}')_d\cdots(\tilde v_{k+2}')_{k+2}(\tilde v_{k+1}')_{k+1}$, which are positive by the condition ${{\sprod{v_i}{\tilde v_i'}} > \delta}$, these are precisely the determinants obtained when
	evaluating \eqref{eq:det1noframe} and \eqref{eq:det2noframe} for the polytope $\Phi_\varepsilon(P)$ and the original frame~$B$.
	This shows that $B$ and $\tilde B'$ are $\Phi_\varepsilon(P)$-equivalent, and hence $B$ and $B'$ are also~$\Phi_\varepsilon(P)$-equivalent.
\end{proof}

\SS\label{def:squashing}

We will now introduce our second polytope operation, called \emph{squashing}, which is closely related to \emph{connected sums}, as defined for instance in \cite{RichterGebert1999}.

Let $B = (v_1,\ldots,v_d)$ be a frame of $\R^d$; let $Q \subset \R^d$ be a $d$-polytope framed by~$B$ that has a facet~$F$ for which the dual basis vector $v_d^\ast$ is a normal covector, that is $v_d^\ast\in\NC^Q_F$; and let $P \subset \R^d$ be a $d$-polytope also framed by~$B$.
We say that $P'$ is obtained by \defn{squashing} $P$ on top of $F$ if
\begin{enumerate}[leftmargin = *]
	\item $P'$ is the image $\psi(P)$ of $P$ under an affine map $\psi = \tau \circ \Phi_\lambda \circ \Phi_\varepsilon$ defined by the composition of the following three maps:
	\begin{enumerate}
		\item a compression $\Phi_\varepsilon$ of the last coordinate in the basis~$B$ with $\varepsilon = (1,\dots,1,\varepsilon_d)$ and~$\varepsilon_d>0$,
		\item a homothety $\Phi_\lambda$ with $\lambda = (\delta,\dots,\delta)$ and $\delta>0$,
		\item a translation $\tau:x\mapsto x+v$ for some $v \in \R^d$;
	\end{enumerate}
	\item every facet $G$ of $Q$ different from $F$ is still a face of $\conv(Q\cup P')$; and
	\item for every facet $G\in \ta_{d-1}(P)$ in the $(d-1)$-target of $P$, its image $\psi(G)$ is a facet of $\conv(Q\cup P')$.
\end{enumerate}

Notice that above we imposed that all facets of $Q$ but $F$ are still facets of the convex hull of $\conv(Q\cup P')$, and also that the target facets of $P'$ are also facets of $\conv(Q\cup P')$; but we do not ask these to be all the facets.
Indeed, $P'$ will be hovering on top of~$F$, and $\conv(Q\cup P')$ will also have some new facets joining the boundary of~$F$ with the \defn{equator} of $P'$, i.e.\ the faces that are both faces of a facet in~$\so_{d-1}(P)$ and of a facet in~$\ta_{d-1}(P)$.
An example is depicted in \cref{fig:squashing}.

\SS

Let~$F$ be a facet of a polytope~$P \subset \R^d$, let $H$ be its supporting hyperplane, and let $p$ be a point in $\R^d$.
We say that~$p$ is \defn{beneath} (resp.\ \defn{beyond}) $F$ if it lies on the same side of $H$ as $P\setminus F$ (resp.\ on the opposite side).
The \defn{dome} of $F$ is the set of points that are beyond $F$ and beneath all other facets of $P$.

\medskip\lemma
Let $P,Q\subset \R^d$ be $d$-polytopes framed by $B = (v_1,\ldots,v_d)$. If there is a facet $F$ of $Q$ with normal covector $v_d^\ast$, then it is possible to squash~$P$ on top of~$F$.

\begin{proof}
	We need to show that an affine map~$\psi = \tau \circ \Phi_\lambda\circ \Phi_\varepsilon$ fulfilling the requirements~(2) and~(3) in~\cref{def:squashing} exists.
	We choose first a homothety~$\Phi_\lambda$ and translation~$\tau$ such that~$(\tau \circ \Phi_\lambda) (P)$ fits in the dome of~$F$.
	This is always possible since the dome is an open subset of $\R^d$ and $P$ is bounded.
	Once this is done, we have that for every $\varepsilon=(1,\dots, 1, \varepsilon_d)$ with $\varepsilon_d\leq 1$, the image $\psi(P) = (\tau \circ \Phi_\lambda\circ \Phi_\varepsilon) (P)\subseteq (\tau \circ \Phi_\lambda) (P)$ is in the dome of~$F$.
	It follows that all facets of \(Q\) except for $F$ will be preserved, as their supporting inequality will be valid on all vertices of $P'$ too.
	Thus, condition~(2) is satisfied.

	It remains to prove that we can further choose $0<\varepsilon_d\leq 1$ so that condition~(3) holds.
	For a positive tuple $\delta\defeq (\delta_1,\dots,\delta_d)$ define $\Phi^\ast_{\delta} \colon V^\ast \to V^\ast$ as the map that scales the dual basis by $v_i^\ast \mapsto \delta_i v_i^\ast$.
	It is straightforward to verify that if $G$ is a facet of $P$ with normal covector $\normal_G^P$, and $\delta \defeq (1,\dots,1,\frac{1}{\varepsilon_d})$, then $\Phi^\ast_{\delta}(\normal_G^P)$ is a normal covector for $\Phi_\varepsilon(G)$ in $\Phi_\varepsilon(P)$. Since translations and dilations do not alter the normal cones, $\Phi^\ast_{\delta}(\normal_G^P)$ is also a normal covector for $\psi(G)$ in $\psi(P)$.
	Now, a facet~$G$ is in~$\ta_{d-1}(P)$ whenever $\sprod{\normal_G^P}{v_k}>0$, that is, whenever the $v_k^\ast$ coordinate of $\normal_G^P$ is positive.
	Therefore, by choosing $\varepsilon_d$ small enough, we can get the normal vector $\normal_{\psi(G)}^{\psi(P)}$ of the image of~$G$ by $\psi$ to be arbitrarily close to $v_d^\ast = \normal_F^Q$ (up to positive scaling) for every facet~$G$ in~$\ta_{d-1}(P)$.
	Now, if the supporting hyperplane $H$ for $F$ in $Q$ is of the form $H = \{ x \in \R^d \mid \sprod{\normal_F^Q}{x} = c\}$, then we know that for every vertex $v$ of $Q$ which is not in $F$, we have $\sprod{\normal_F^Q}{v}<c$.
	Since there is a finite number of vertices, we can choose the normal vector $\normal_{\psi(G)}^{\psi(P)}$ to be close enough to $\normal_F^Q$ so that the inequality $\sprod{\normal_{\psi(G)}^{\psi(P)}}{v}<c$ is still valid for every vertex of $Q$ not in $F$.
	Since $P'$ is beyond $F$, the supporting halfspace for $\psi(G)$ is of the form $\sprod{\normal_{\psi(G)}^{\psi(P)}}{x}\leq c'$ with a $c'>c$.
	Therefore, all vertices in $Q$ are also beneath this halfspace and thus $\psi(G)$ is still a facet of~$\conv(Q \cup P')$, finishing the proof.
\end{proof}

\SS\label{lem:reorientationequivalence}
We say that a frame $B' = (v_1',\dots, v_d')$ is a \defn{reorientation} of another frame $B = (v_1,\dots,v_d)$ if $v_i' \in \{ v_i,-v_i\}$ for each $1\leq i\leq d$.

\medskip\lemma
Let $P \subset \R^d$ be a polytope, $B$ be a $P$-admissible frame, and $B'$ a reorientation of~$B$.
Then, the framed polytope $(P,B')$ is loop-free if and only if $(P,B)$ is.

\begin{proof}
	Replacing $v_k$ by $-v_k$ in $B$ only swaps the role of $(k-1)$-sources and targets, thus reversing the order of cellular $k$-strings and preserving loops.
\end{proof}

\SS\label{def:stellar-sub}
We will need one more operation.
Let $Q$ be a polytope and $F$ one of its faces, with $b_Q$ and $b_F$ denoting their respective barycenters.
The \defn{stellar subdivision} of~$F$ is the convex hull of $Q$ and a point $p_F \defeq b_F+\zeta_F(b_F- b_Q)$, for some $\zeta_F>0$ small enough so that $p_F$ is in the dome of~$F$.
Stellar subdivision is a classical operation, see for example \cite[Ex.~3.0]{Ziegler95}.
Combinatorially, this operation replaces each face~$G$ of~$Q$ containing~$F$ by the cone of~$p_F$ over the faces of~$G$ not containing~$F$.
Notice that by taking $\zeta_F$ small enough, we can also impose that the supporting hyperplanes and normal covectors of the new facets are arbitrarily close to those of the deleted facets (more precisely, these hyperplanes share a $(d-2)$-dimensional affine subspace, and we can impose that the angle between their normal covectors is arbitrarily small).

\SS\label{thm:noframe}
\theorem
Let $P$ be a framed $d$-polytope with a cellular $k$-loop for some~$k \leq d-2$.
If all the faces in this loop are faces of faces in~$\ta_{d-1}(P)$, then there is a $d$-polytope $\tilde P$ for which every admissible frame induces a $k$-loop.

\begin{proof}
		Without loss of generality (see \cref{ss:lowertriangular,ss:invariancelinearautomorphism}) we will restrict to orthonormal frames, and assume that the given~frame~ $B$ is the canonical frame $(e_1,\dots,e_d)$. For a permutation $\sigma\in \sym_d$, define $B_\sigma = (e_{\sigma(1)},\dots,e_{\sigma(d)})$. 
		
		We set $\delta\defeq \frac{1}{2\sqrt{d}}$ and use \cref{ss:flattening} to find a positive real tuple $\varepsilon$ and a polytope $P' \defeq \Phi_\varepsilon(P)$
		such that for any orthonormal frame $B' = (v_1',\dots,v'_d)$,  	
		if \begin{equation}\label{eq:bound12sqrtd}
			\sprod{e_{i}}{\tilde v_{i}'}\geq \delta =\frac{1}{2\sqrt{d}}\quad 1 \leq i \leq d,
		\end{equation} 
where the vector~$\tilde v_i'$ is defined as in \cref{lem:directsum} and moreover rescaled to have length $1$,
		then $\tilde B'$ and $B$ are $\Phi_\varepsilon(P)$-equivalent. 
		
		We claim now that for any orthonormal frame $B' = (v_1',\dots,v'_d)$ there is some permutation $\sigma\in \sym_d$ such that
		 \begin{equation}\label{eq:permutationboundsqrtd}\left\lvert{\sprod{e_{\sigma(i)}}{\tilde v_{i}'}}\right\rvert\geq\frac{1}{\sqrt{d}}\quad 1 \leq i \leq d.\end{equation}
		 Here, the vector~$\tilde v_i'$ is again defined as in \cref{lem:directsum} but for the frame~$B_\sigma$, and moreover rescaled to have length $1$.
	Indeed, observe first that it is impossible for a unit vector $v_d'$ to be such that \(\left\lvert{\sprod{e_i}{v_d'}}\right\rvert<\frac{1}{\sqrt{d}}\) for all~$1 \leq i \leq d$. Therefore we can choose $\sigma(d)$ so that \(\left\lvert{\sprod{e_{\sigma(d)}}{v_d'}}\right\rvert\geq \frac{1}{\sqrt{d}}\).
	Similarly, it is impossible to have \(\left\lvert{\sprod{e_j}{\tilde v_{d-1}'}}\right\rvert<\frac{1}{\sqrt{d}}<\frac{1}{\sqrt{d-1}}\) for all~$j\neq \sigma(d)$.
	Therefore we can choose $\sigma(d-1)$ so that \(\left\lvert{\sprod{e_{\sigma(d-1)}}{\tilde v_{d-1}'}}\right\rvert\geq \frac{1}{\sqrt{d}}\).
	Continuing in this fashion, we obtain a permutation~$\sigma$ with the property that \(\left\lvert{\sprod{e_{\sigma(i)}}{\tilde v_{i}'}}\right\rvert\geq \frac{1}{\sqrt{d}}\) for $1\leq i\leq d$.
	Choosing a suitable reorientation~$\tilde B'$ of~$B'$ with $\tilde v_{i}'\in\{v_{i}',-v_{i}'\}$, we can further impose that \({\sprod{e_{\sigma(i)}}{\tilde v_{i}'}}\geq \frac{1}{\sqrt{d}}\) for $1\leq i\leq d$.
	
	For every $\sigma\in \sym_d$, let $P_\sigma'$ be the symmetric copy of $P'$ obtained by permuting the coordinates according to $\sigma$. By \cref{ss:invariancelinearautomorphism} and the defining property of $P'$, we have that $\tilde B'$ and $B_\sigma$ are $P_\sigma '$-equivalent, and thus $\tilde B'$ induces a loop on~$P_\sigma'$. By \cref{lem:reorientationequivalence}, so does~$B'$.
	
	The idea now is to combine all the possible $P_\sigma'$ to construct our never loop-free polytope~$\tilde P$.
	This has to be done carefully, because we need the faces involved in the loops to still be faces of the convex hull of all these reflected copies.
	
	We start with an auxiliary simplex $\Delta_d\defeq \conv(0,-e_1,\ldots,-e_d)$, on which we perform stellar subdivisions on all the proper faces of~$\Delta_d$ containing the origin.
	We do so in an order in which the dimension of the faces is non-increasing (i.e. first all faces of dimension~$d-1$, then those of dimension $d-2$, etc.).
	We get a polytope in which for every permutation $\sigma \colon [d] \to [d]$ there is a facet $$F_\sigma \defeq \conv(p_0,p_{0\sigma{(1)}},\dots,p_{0\sigma(1)\dots\sigma(d-1)})$$ given by the convex hull of the extra points $p_{0\sigma(1)\dots\sigma(k)}$ introduced when performing the stellar subdivisions of the faces $\conv(0,-e_{\sigma(1)},\dots,-e_{\sigma(k)})$.
	As we observed in \cref{def:stellar-sub}, we have freedom in the choice of these extra points, and we can choose them in such a way that the normal covector $\normal_\sigma$ of the facet $F_\sigma$ is as close as needed to the normal covector of the facet $\conv(0,-e_{\sigma(1)},\dots,-e_{\sigma(d-1)})$, which is $e_{\sigma(d)}$ under the standard identification of the primal with the dual using the scalar product (here, by close we mean that their angle is as small as needed).

	To make this more precise, denote by $\rho_\sigma$ the rotation that sends $e_{\sigma(d)}$ to the new normal covector $\normal_\sigma$ while fixing the orthogonal complement of the plane spanned by these two vectors.
Recall from \cref{def:stellar-sub} that when performing the stellar subdivision, we can require the angle between $e_{\sigma(d)}$ and $\normal_\sigma=\rho_\sigma(e_{\sigma(d)})$, which are the normal covectors to old and new facets, to be arbitrarily small.	

We construct the stellar subdivisions so that this angle is sufficiently small to guarantee the following property: for any orthonormal frame $B'$ fulfilling \eqref{eq:permutationboundsqrtd}
we also have 
\begin{equation}\label{eq:permutationbound12sqrtd}
	\left\vert{\sprod{\rho_\sigma(e_{\sigma(i)})}{\tilde v_{i}'}}\right\vert\geq \frac{1}{2\sqrt{d}}\quad 1 \leq i \leq d.\end{equation} 
This is always possible because, for any unit vectors $u, v, w $, if $\alpha$ denotes the angle between $u$ and $v$, then
\[
|\sprod{u}{w} - \sprod{v}{w}| \leq \alpha.
\]
Hence, choosing the rotation angle $\alpha$ small enough ensures the stated inequalities.
	
	Now, for every permutation $\sigma\in \sym_d$, consider $\rho_\sigma(P_\sigma')$. Note that by \cref{ss:invariancelinearautomorphism}, the frame $\rho_\sigma(B_\sigma)$ is $\rho_\sigma(P_\sigma')$-admissible and induces a loop, and so does every frame close enough to it in the sense of \cref{ss:flattening}.
	
	We now squash $\rho_\sigma(P_\sigma')$ on top of the facet $F_\sigma$ for every $\sigma\in \sym_d$.
	This might compress even further the polytope in direction $\normal_\sigma$, but the properties we need from \cref{ss:flattening} still hold.
	Denoting the squashed polytopes by $\tilde P_\sigma'$, the never-loop-free polytope we are looking for is then
	$$\tilde P \defeq \conv\left(\bigcup_{\sigma\in \sym_d} \tilde P_\sigma'\right).$$
	
	Indeed, we have seen first that for every $\tilde P$-admissible frame $B'$ there is some permutation~$\sigma$ and a reorientation $\tilde B'$ of $B'$, such that \eqref{eq:permutationboundsqrtd} holds, and hence also~\eqref{eq:permutationbound12sqrtd}. This implies that $\tilde B'$ and $\rho_\sigma(B_\sigma)$ are $\rho_\sigma(P_\sigma')$-equivalent.
	Since $\rho_\sigma(B_\sigma)$ induces a loop on $\rho_\sigma(P_\sigma')$, so does $\tilde B'$.
	This loop is also a loop for~$\tilde P$, as all the faces involved are faces of~$\tilde P$.
	Therefore, by \cref{lem:reorientationequivalence} also $B'$ induces a loop on~$\tilde P$. This concludes the proof.
\end{proof}

\SS\label{cor:noframe}\corollary
There is a $4$-polytope $\tilde P$ for which every admissible frame induces a cellular loop.

\begin{proof}
	It suffices to show that the cyclic polytope $P_4$ defined in \cref{ss:1loop} together with the loop \eqref{eq:1loop} satisfy the condition of \cref{thm:noframe}.
	Let us recall that~$P_4$ was defined as the $4$-projection of the $5$-simplex~$P_5$.
	Explicitly, $P_4$ is the convex hull of the columns $(p_1,\dots,p_6)$ of the matrix
	\[
	\begin{pmatrix}
		-3 & -2 & -1 & 1 & 2 & 3 \\
		-1 & 1 & 0 & 0 & 1 & -1 \\
		-1 & 1 & 0 & 0 & -1 & 1 \\
		0 & 0 & 1 & 1 & 0 & 0
	\end{pmatrix}.
	\]
	As stated in \cref{ss:1loop} it has a $1$-loop given by
	\begin{equation*}
		[p_1p_2p_3],\ [p_2p_3p_6],\ [p_2p_6p_4],\ [p_4p_5p_6],\ [p_1p_4p_5],\ [p_1p_3p_5],\ [p_1p_2p_3].
	\end{equation*}

	All the facets of $P_4$, except for the $3$-dimensional simplex $[p_1p_2p_5p_6]$, belong to the $3$-target $\ta_{3}(P)$.
	In particular, all the $2$-faces of $P_4$ are the face of a face in $\ta_{3}(P)$, as required by the theorem.
\end{proof}

%% file: fig/flattening.tex
\begin{tikzpicture}[auto, scale = 1.4, point/.style = {draw, circle, fill = black, inner sep = 1pt}]
	\begin{scope}[shift = {(4,-.5)}]
		\fill[RedViolet!20] (0,0) -- (0.9*0.918,0.9*0.397) arc(23.36:156.63:0.9) -- cycle;
		\fill[RoyalBlue!40] (0,0) -- (0.95*0.5,0.95*0.866) arc(60:120:0.95*1) -- cycle;

		\draw[dashed,RedViolet] (-0.918,0.397) -- (0,0) -- (0.918,0.397);
		\draw[dashed,RoyalBlue] (-0.5,0.866) -- (0,0) -- (0.5,0.866);

		\draw[->] (0, 0) -- (1, 0) node[anchor = west]{$v_1$};
		\draw[->] (0, 0) -- (0, 1) node[anchor = west]{$v_2$};
		\draw[->, dashed] (0, 0) -- ( 0.707 , 0.707) node[anchor = west]{$v_2'$};
		%


	\end{scope}

	\coordinate (a) at (1, 0);
	\coordinate (b) at (0.5, 0.866);
	\coordinate (c) at (-0.5, 0.866);
	\coordinate (d) at (-1, 0);
	\coordinate (e) at (-0.5, -0.866);
	\coordinate (f) at (0.5, -0.866);

	\coordinate (a') at (1, 0.25*0);
	\coordinate (b') at (0.5, 0.25*0.866);
	\coordinate (c') at (-0.5, 0.25*0.866);
	\coordinate (d') at (-1, 0.25*0);
	\coordinate (e') at (-0.5, -0.25*0.866);
	\coordinate (f') at (0.5, -0.25*0.866);

	\draw[thick,RoyalBlue] (a)--(b)--(c)--(d)--(e)--(f)--(a);

	\draw[thick,dashed,RedViolet] (a')--(f')--(e')-- (d')--(c')--(b')--(a') ;

	\foreach \i in {a,b,c,d,e,f}
	{
		\node[point] at (\i) {};
		\node[point] at (\i') {};

	}
\end{tikzpicture}

%% file: fig/squashing.tex
\begin{tikzpicture}[auto, scale = 1.4, point/.style = {draw, circle, fill = black, inner sep = 1pt}]

 \begin{scope}[shift = {(0,.5)}, scale = .7]

	\coordinate (e) at (-1, 0);
	\coordinate (f) at (0, 1);
	\coordinate (g) at (1, 0);
	\coordinate (h) at ( 0,-1);

	\draw[thick,Green, fill = Green!30] (e)--(f)--(g)--(h)--(e);

	\foreach \i in {e,f,g,h}
	{
		\node[point,Green] at (\i) {};

	}

		\node at (0,0) {$P$};

 \end{scope}

 \begin{scope}[shift = {(2.5,0)}]

	\coordinate (a) at (-1, 0);
	\coordinate (b) at (-1, 1);
	\coordinate (c) at (1, 1);
	\coordinate (d) at (1, 0);

	\draw[thick,RoyalBlue, fill = RoyalBlue!30] (a)--(b)--(c)--(d)--(a);

	\foreach \i in {a,b,c,d}
	{
		\node[point,RoyalBlue] at (\i) {};

	}

	\node at (0,.5) {$Q$};
	\node at (0,1.2) {$F$};

 \end{scope}

 \begin{scope}[shift = {(6,0)}]

	\coordinate (a') at (-1, 0);
	\coordinate (b') at (-1, 1);
	\coordinate (c') at (1, 1);
	\coordinate (d') at (1, 0);

	\coordinate (e') at (-.8, 1.5);
	\coordinate (f') at (0, 1.6);
	\coordinate (g') at (.8, 1.5);
	\coordinate (h') at ( 0,1.4);

%
%

	\fill[fill = Green!30] (e')--(f')--(g')--(h')--(e');

	\fill[fill = RoyalBlue!30] (a')--(b')--(c')--(d')--(a');

 \draw[thick,Green] (e')--(f')--(g');

	\draw[thick,RoyalBlue] (c')--(d')--(a')-- (b');

	\draw[thick,Green, dashed] (g')--(h')--(e');

	\draw[thick,RoyalBlue, dashed ] (b')--(c');

	\draw[thick,OrangeRed] (b')--(e');

	\draw[thick,OrangeRed] (g')--(c');

	\foreach \i in {e,f,g,h}
	{
		\node[point,Green] at (\i') {};

	}
		\foreach \i in {a,b,c,d}
	{
		\node[point,RoyalBlue] at (\i') {};

	}

	\node at (0,.5) {$Q$};
	\node at (0,1.8) {$P' = \psi(P)$};

 \end{scope}

\end{tikzpicture}

%
%
%
%
%
%
%
%
%

%% file: sec/random.tex

\subsection{Random frames and random polytopes}\label{s:random}

We now approach the existence of cellular loops from a probabilistic viewpoint.
Specifically, we prove that the canonical frame on a Gaussian simplex, and a random frame on the canonical simplex, induce a loop almost surely as the dimension increases.

\SS\label{ss:random_simplices}
A \defn{Gaussian $d$-simplex} is the convex hull of $d+1$ independent random points in~$\R^d$, each chosen according to a $d$-dimensional standard normal distribution.

\medskip\theorem
For every $k \geq 1$, the probability that the canonically framed Gaussian $d$-simplex has a $k$-loop tends to~$1$ as~$d$ tends to~$\infty$.

\begin{proof}
	Let $X$ be a Gaussian $d$-simplex, and let $x_1,\dots,x_{d+1}$ be its vertices.
	We start with the case $k = 1$.
	Consider the point configuration $p_1,\dots,p_6$ from \cref{ss:1loop}, whose $2$-dimensional projection is depicted in \cref{fig:counterexample}.
	We show that with high probability the canonical projection $\pi_2(X)$ of $X$ to $\R^2$ contains a vertex in each of the neighborhoods of the canonical projections~$\pi_2(p_i)$ of the points~$p_i$, for $1\leq i \leq 6$.
	Then, the $2$-faces of $X$ induced by these vertices will define a $1$-loop.

	Indeed, having a loop is an open condition, and therefore there is some $\varepsilon>0$ such that any simplex with vertices $p_1',\dots,p_6'$ with $p_i' \in \ball{p_i}{\varepsilon}$, where $\ball{p_i}{\varepsilon}$ denotes the Euclidean ball of radius~$\varepsilon$ centered at~$p_i$, will also present a $1$-loop.
	The measure of each of these balls is positive, and therefore there is some $\delta>0$ such that the probability that $\pi_2(x_j)$ lies in $\ball{p_i}{\varepsilon}$ is at least~$\delta$.
	Since the $x_j$ are chosen according to independent random distributions, for each $1\leq i\leq 6$, the probability that at least one of the $\pi_2(x_j)$ with $\frac{(i-1)d}{6}\leq j <\frac{id}{6}$ lies in $\ball{p_i}{\varepsilon}$ is at least
	\[
	1-(1-\delta)^{\left\lfloor\frac{d}{6}\right\rfloor},
	\]
	and the probability of having at least one $\pi_2(x_j)$ in each $\ball{p_i}{\varepsilon}$ is at least
	\[
	\Big(1-(1-\delta)^{\left\lfloor\frac{d}{6}\right\rfloor}\Big)^6,
	\]
	which tends to~$1$ as $d\to \infty$.
	The same proof holds for higher~$k$, starting with a framed polytope admitting a $k$-loop instead of a $1$-loop.
	For example, one can use~~\cref{ss:hyperplane}, by virtue of which a $(k-1)$-fold pyramid over the polytope from~\cref{ss:1loop} presents a $k$-loop.
\end{proof}

While the proof is certainly not optimized to find the best convergence rate, very few hypotheses on the distribution are used.
In particular, it suffices that the distribution has support \(\R^d\) and that the vertices are independently sampled.
These hypotheses could also be significantly relaxed.
For example, it is enough that a homothety of the configuration $p_1,\dots,p_6$ lies in the support.
Therefore, for most usual distributions of random simplices the same kind of result should hold.


\SS\label{ss:random_frames}
We obtain similar results if instead of fixing the frame and choosing the simplex, we fix the simplex and choose the frame.
For example, in view of \cref{ss:equivalentorthogonal}, a reasonable approach is to consider a random orthogonal frame chosen with respect to the Haar measure.
We consider the canonical embedding~$\ssimplex[d] = \conv(e_1,\ldots,e_{d+1})$ of the standard $d$-simplex in $\R^{d+1}$.

\medskip\theorem
For every $k \geq 1$, the probability that a uniform random orthonormal frame induces a $k$-loop on the standard $d$-simplex tends to~$1$ as~$d$ tends to~$\infty$.

\begin{proof}
	Let $O$ be an orthogonal matrix chosen uniformly at random, whose columns define a frame~$B$.
	Since the $f$-orientation of a framed polytope is invariant under positive scaling of the frame vectors, we can assume that $O$ is orthonormal.
	In this case, we have $O^{-1} = O^t$.
	By \cref{ss:invariancelinearautomorphism}, applying $O^t$ to both $\ssimplex[d]$ and $B$ does not alter the $f$-orientation.
	Thus $(\ssimplex[d],B)$ has the same $f$-orientation as~$(O^t\ssimplex[d],O^tB)$.
	Now, $O^tB$ is just the canonical frame, and we are in the setup of~\cref{ss:random_simplices}.
	The coordinates of $\pi_k(O^t\ssimplex[d])$ are thus the first $k$ rows of $O^t$, which are the first $k$ vectors of $B$.
	It is known that if these are chosen from a uniform random orthogonal matrix then they can be approximated by independent standard normals when $d\to \infty$ (if $k = o(d/\log d)$, which holds in particular if $k$ is a constant), see for example \cite[Thm.~3]{Jiang2006}.
	Therefore, the same proof as in \cref{ss:random_simplices} works in this setting.
\end{proof}

\medskip\corollary
For every $k \geq 1$, the probability that a uniform random orthonormal frame induces a $k$-loop on the unit $d$-dimensional cube tends to~$1$ as~$d$ tends to~$\infty$.

\begin{proof}
	Let $v$ be the vertex of the unit cube $C$ that maximizes the linear functional associated to the last vector $v_d$ of the frame.
	Let $H$ be an affine hyperplane perpendicular to $v_d$ that separates $v$ from the remaining vertices of $C$. Then $C\cap H$ is a simplex, and $\{v_1,\dots,v_{d-1}\}$ is a random orthogonal frame of its affine hull.
 Almost surely all coordinates of $v_d=(a_1,\dots,a_n)$ are non-zero, and by symmetry, we can assume that $a_i<0$ for $1\leq i\leq n$.
 That is, $v$ is the origin and $H$ is of the form $\sum x_ia_i=\epsilon$ for some $\epsilon<0$. The vertices of $C\cap H$ are $\frac{\epsilon}{a_i}e_i$, and $\{v_1,\dots,v_{d-1}\}$ is an orthogonal frame chosen at random on $H$.

 We can apply the same reasoning of \cref{ss:random_frames} to show that $C\cap H$ and $\{v_1,\dots,v_{d-1}\}$ asymptotically almost surely has a loop. The only difference is that now, the coordinates of $\pi_k(O^t\ssimplex[d])$ are the first $k$ rows of $O^t$ with their coordinates scaled by $\frac{\epsilon}{a_i}$ (instead of just the first $k$ rows). But we can still apply the proof from \cref{ss:random_simplices} to this distribution.

 Finally, we apply \cref{ss:hyperplane} to get that if $(C\cap H,\{v_1,\dots,v_{d-1}\})$ has a loop, then $(C,\{v_1,\dots,v_{d}\})$ also has a loop.
\end{proof}

%% file: sec/conjecture.tex

\section{The Kapranov--Voevodsky Conjecture}\label{part:connection}

Kapranov and Voevodsky conjectured in \cite[Thm.~2.3]{kapranov1991polycategory} that a certain procedure on any polytope equipped with an admissible flag of projections produces a pasting diagram.
We reformulate this conjecture in the language of framed polytopes and use Steiner's diagrams \cite{steiner2004omega} to show that the absence of cellular loops is both necessary and sufficient for its validity.
Since the absence of directed loops is required for the definition of pasting diagrams in any reasonable formalization, the examples from the previous section provide model-independent counterexamples to the conjecture.

\input{sec/categories}
\input{sec/pasting_diagrams}
\input{sec/steiner}
\input{sec/poly_diag}
\input{sec/reformulation}

%% file: sec/categories.tex

\subsection{Higher categories}

We review the definitions of globular sets and \(n\)-categories, emphasizing how they generalize directed graphs and small categories, respectively.

\SS Recall that a directed graph consists of a set of vertices~\(V\), a set of edges~\(E\), and two functions~\(\so, \ta \colon E \to V\) called the source and target maps.
A globular set generalizes this notion by including higher-dimensional directed edges.
In precise terms, a \defn{globular set} is an \(\N\)-graded set \(X = \set{X_n}_{n \geq 0}\) together with \defn{source} and \defn{target} maps
\[
\begin{tikzcd}[column sep = small]
	X_0 &
	\arrow[l,shift right = 3pt,"\;\so"'] \arrow[l,shift left = 3pt,"\;\ta"] X_1 &
	\arrow[l,shift right = 3pt,"\;\so"'] \arrow[l,shift left = 3pt,"\;\ta"] X_2 &
	\arrow[l,shift right = 3pt,"\;\so"'] \arrow[l,shift left = 3pt,"\;\ta"] \dotsb
\end{tikzcd}
\]
satisfying the \defn{globular identities}: \(\so \circ \so = \so \circ \ta\) and \(\ta \circ \so = \ta \circ \ta\).
We refer to elements in \(X_n\) as \defn{\(n\)-cells}, and to \(n\) as their \defn{dimension}.
Using the globular identities, we can define, for any \(0 \leq k < n\), the \defn{\(k\)-source} and \defn{\(k\)-target} maps as the unique morphisms \(\so_k, \ta_k \colon X_n \to X_k\) that factor respectively through \(\so, \ta \colon X_{k+1} \to X_k\).

A \defn{morphism of globular sets} is a morphism of \(\N\)-graded sets \(\set{f_n \colon X_n \to X'_n}\) commuting with source and target maps.
We denote the category of globular sets by \(\bG\Set\).

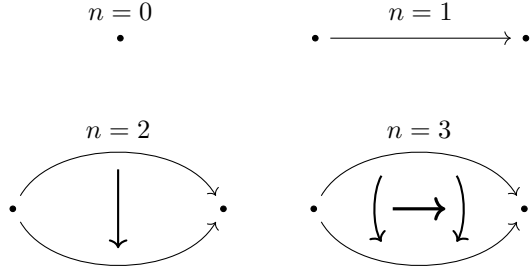
\begin{figure}
	\input{fig/globe}
	\caption{Pictorial representation of the \(n\)-globe \(\bG^n\) for small values of \(n\).}
	\label{f:globes}
\end{figure}

\SS\label{ss:n-globe} Alternatively, the category of globular sets is the category of contravariant functors from a category \(\bG\), termed the \defn{globe category}, to the category of sets.
Without presenting all the details, we mention that the objects of \(\bG\) are the non-negative integers and that the representable globular set \(\bG^n \defeq \bG(-,n)\), termed the \defn{\(n\)-globe}, is given by
\[
\bG^n_n = \set{e_n}
\quad\text{and}\quad
\bG^n_k = \set[\big]{e_k^-,e_k^+}
\quad\text{for}\quad 0 \leq k < n,
\]
with
\begin{align*}
	&\so(e_n) = e_{n-1}^- \ ;\
	\ta(e_n) = e_{n-1}^+
	&&0 < n, \\
	&\so(e_k^+) = \so(e_k^-) = e_{k-1}^- \ ;\
	\ta(e_k^+) = \ta(e_k^-) = e_{k-1}^+
	&&0 < k < n.
\end{align*}
Consult \cref{f:globes} for a pictorial representation of this globular set for small values of \(n\).

\SS A small \defn{category} is a directed graph
\[
\begin{tikzcd}[column sep = 15pt]
	\Ob & \arrow[l,"\ \so"',shift right = 3pt] \arrow[l,"\ \ta",shift left = 3pt] \Mor
\end{tikzcd}
\]
endowed with 1) a \defn{composition} map \(\circ \colon \Mor \times_0 \Mor \to \Mor\), where
\[
\Mor \times_0 \Mor \defeq \set[\big]{(g,f) \in \Mor \times \Mor \mid \so(g) = \ta(f)},
\]
and 2) an \defn{identity} map \(\id \colon \Ob \to \Mor\), satisfying relations of associativity and unitality.

Similarly, an \defn{\(\omega\)-category} is a globular set together with 1) \defn{\(k\)-composition} maps \(\circ_k \colon X_n \times_k X_n \to X_n\) for all \(0 \leq k < n\), where
\[
X_n \times_k X_n \defeq \set[\big]{(y,x) \in X_n \times X_n \mid \so_k(y) = \ta_k(x)},
\]
and 2) \defn{identity} maps \(\id_{n+1} \colon X_n \to X_{n+1}\) for every \(n \in \N\),
satisfying relations of associativity, unitality, and interchange.
Elements in \(X_n\) are called \defn{\(n\)-morphisms} and \(n\) is their \defn{dimension}.
Additionally, we refer to elements in the image of identity maps as \defn{identities}.
Note that one can compose morphisms of different dimensions using these.
Explicitly, if an \(m\)-morphism \(y\) and an \(n\)-morphism \(x\) with \(m < n\) satisfy \(\so_k(y) = \ta_k(x)\), then \(\id_{n}(y) \circ_k x\) is well-defined.
We refer to \mbox{\(\omega\)-categories} whose morphisms of dimension greater than \(n\) are all identities as \defn{\(n\)-categories}.
An equivalent recursive definition is as follows: a \(1\)-category is a category, and an \(n\)-category is a category enriched in \((n-1)\)-categories.

A \defn{morphism of \(\omega\)-categories} is a morphism of underlying globular sets which commutes with identities and compositions.
We denote by \(\wCat\) the category of \(\omega\)-categories.

%% file: fig/globe.tex
\begin{tikzpicture} [scale = .35]
	\draw node at (-4,0){ };
	\draw node at (4,0){ };
	\draw[fill] (-0.4,0) circle [radius = 3pt];

	\draw node at (-0.4,1){$n = 0$} ;
	\draw node at (0,-2.4){} ;
\end{tikzpicture}
\begin{tikzpicture} [scale = .35]
	\draw[fill] (-4,0) circle [radius = 3pt];
	\draw[fill] (4,0) circle [radius = 3pt];

	\draw [->] (-4,0) [shorten > = 0.2cm, shorten < = 0.2cm,->] to (4,0) ;

	\draw node at (0,1){$n = 1$} ;
	\draw node at (-6,-2.4){} ;
\end{tikzpicture}

\begin{tikzpicture} [scale = .35]
	\draw[fill] (-4,0) circle [radius = 3pt];
	\draw[fill] (4,0) circle [radius = 3pt];

	\draw [->] (-4,0) [shorten > = 0.2cm, shorten < = 0.2cm,->, out = 60,in = 120] to (4,0) ;
	\draw [->] (-4,0) [shorten > = 0.2cm, shorten < = 0.2cm,->, out = -60,in = -120] to (4,0) ;

	\draw [->] (0,1.5) [thick] to (0,-1.5) ;

	\draw node at (0,3){$n = 2$} ;
\end{tikzpicture}
\hspace*{25pt}
\begin{tikzpicture} [scale = .35]
	\draw[fill] (-4,0) circle [radius = 3pt];
	\draw[fill] (4,0) circle [radius = 3pt];

	\draw [->] (-4,0) [shorten > = 0.2cm, shorten < = 0.2cm,->, out = 60,in = 120] to (4,0) ;
	\draw [->] (-4,0) [shorten > = 0.2cm, shorten < = 0.2cm,->, out = -60,in = -120] to (4,0) ;

	\draw [->] (-1,2) [thick,shorten > = 0.3cm, shorten < = 0.3cm,->, out = -120,in = 120] to (-1,-2) ;
	\draw [->] (1,2) [thick,shorten > = 0.3cm, shorten < = 0.3cm,->, out = -60,in = 60] to (1,-2) ;

	\draw [->] (-1,0) [very thick] to (1,0) ;

	\draw node at (0,3){$n = 3$} ;
\end{tikzpicture}

%% file: sec/pasting_diagrams.tex

\subsection{Pasting diagrams}\label{ss:diagrams}

In general, due to the relations satisfied by $\omega$-categories, morphisms can be written as compositions in multiple equivalent ways.
For instance, as presented in Forest \cite{forest2022pasting}, both expressions
\begin{gather*}
	a \circ_0 (\alpha \circ_{1} \beta) \circ_{0} ( (\gamma \circ_{0} h) \circ_{1} (\delta \circ_{0} h))
	\shortintertext{and}
	(a \circ_{0} \alpha \circ_{0} e \circ_{0} h) \circ_{1} (a \circ_{0} c \circ_{0} \gamma \circ_{0}
	h) \circ_{1} (a \circ_{0} \beta \circ_{0} \delta \circ_{0} h ).
\end{gather*}
represent the same morphism in any $\omega$-category, as they are related by the associativity and interchange laws.
By passing to higher-dimensional compositional diagrams, we can obtain uniqueness in the representation of certain compositions.
For example, the above morphism is completely determined by the following diagram:
\[
\begin{tikzcd}[ampersand replacement = \&]
	u \arrow[rr,"a"] \& \& v \arrow[rr,"c"{description},""{auto = false,name = c}]
	\arrow[rr,out = 70,in = 110,"b",""{auto = false,name = b}]
	\arrow[rr,out = -70,in = -110,"d"',""{auto = false,name = p}]
	\arrow[phantom,"\Downarrow \alpha",from = b,to = c]
	\arrow[phantom,"\Downarrow \beta",from = c,to = p]
	\& \& w \arrow[rr,"f"{description},""{auto = false,name = f}] \arrow[rr,out = 70,in = 110,"e",""{auto = false,name = e}]
	\arrow[rr,out = -70,in = -110,"g"',""{auto = false,name = g}] \& \& x \arrow[rr,"h",""{auto = false,name = h}] \& \& y.
	\arrow[phantom,"\Downarrow \gamma",from = e,to = f]
	\arrow[phantom,"\Downarrow \delta",from = f,to = g]
\end{tikzcd}
\]
More precisely, once the images of $a,\dots,h$ and $\alpha,\dots,\delta$ are appropriately assigned, there is a unique functor from the free $\omega$-category generated by this diagram into any $\omega$-category.
This property is a key requirement in any theory of \textit{pasting diagrams}, for which multiple formalisms have been proposed, including pasting schemes \cite{johnson1986pasting}, parity complexes \cite{street1991paritycomplexes}, regular directed complexes \cite{Steiner93, Hadzihasanovic24}, augmented directed complexes \cite{steiner2004omega}, and torsion-free complexes \cite{forest2022pasting}.
All models of pasting diagrams disallow directed loops.
For example, it is unclear whether the diagram
\[
\begin{tikzcd}
	a \rar[bend left,"f"] & \lar[bend left,"g"] b
\end{tikzcd}
\]
represents \(f \circ g\) or \(g \circ f\) or any other number of possibilities.

%% file: sec/steiner.tex

\subsection{Steiner diagrams}\label{s:steiner}

We revisit the concept of \textit{augmented directed complex with a loop-free basis} introduced in \cite{steiner2004omega}, which we abbreviate as \textit{Steiner diagram}.
Denoting their category by \(\wDiag\), our goal is to introduce the definitions required to state the following result justifying their role as models for pasting diagrams.

\medskip\theorem(\cite[Thm~5.11]{steiner2004omega})
The functor
\[
\cells \colon \wDiag \to \wCat
\]
is full and faithful, with each resulting \(\omega\)-category freely generated by its atoms.

\medskip For the reader familiar with Steiner's work, we mention that we will provide a novel characterization of the functor \(\cells\) as a globular nerve.

\SS By an \defn{augmented chain complex}, we mean an \(\N\)-graded abelian group \(C = \set{C_n}_{n \in \N}\) with \defn{boundary} and \defn{augmentation} maps
\[
\Z \xla{\,\varepsilon} C_0 \xla{\,\bd} C_1 \xla{\,\bd} C_2 \xla{\,\bd} \dotsb,
\]
satisfying \(\bd \circ \bd = \varepsilon \circ \bd = 0\).
We refer to elements in \(C_n\) as \defn{\(n\)-chains} and to \(n\) as their \defn{degree}.

A \defn{chain map} is a morphism of \(\N\)-graded abelian groups \(f = \set{f_n \colon C_n \to C_n'}\) commuting with boundary maps.
It is said to be \defn{augmentation-preserving} if it also commutes with augmentation maps.

We say that an augmented chain complex \(C\) is \defn{based} when each group \(C_n\) is equipped with a basis, i.e. an \(\N\)-graded set \(B = \set{B_n}_{n\in\N}\) such that \(C_n = \Z\set{B_n}\) for each \(n \in \N\), and \(\varepsilon(b) = 1\) for each \(b \in B_0\).
A \defn{positive chain} in \(C\) is a linear combination of basis elements with positive coefficients.
We say that a chain map between based chain complexes is \defn{non-negative} if the image of each basis element is a positive element or \(0\).

We denote the category of based augmented chain complexes with augmentation-preserving non-negative chain maps as \(\bach\).

\SS Let us consider the functor of \defn{globular chains} \(\chains \colon \bG\Set \to \bach\) sending a globular set \(X\) to the based augmented chain complex with basis \(X\) and boundary \(\ta - \so\).
By a standard argument, this functor is completely determined by its restriction to the full subcategory of \(n\)-globes.
The \defn{globular nerve} \(\nu \colon \bach \to \bG\Set\) is the right adjoint of \(\chains\).
Explicitly, an \(n\)-cell of \(\nu(C)\) is an augmentation-preserving non-negative chain map \(\chains(\bG^n) \to C\).

\SS Let us consider a based augmented chain complex \(C\).
For any \(c \in C\), let \(\bd^-c\) and \(\bd^+c\) be the positive elements determined by the identity
\[
\bd c = \bd^+ c - \bd^-c
\]
and define recursively, for \(0 \leq k < n\) and \(\epsilon \in \{-,+\}\),
\[
\angles{c}_k^\epsilon \defeq
\begin{cases}
	\bd^\epsilon c & k = n-1, \\
	\bd^\epsilon \angles{c}^\epsilon_{k+1} & k < n-1.
\end{cases}
\]

It is straightforward to check that if \(c\) is a positive \(n\)-chain in \(C\), then there is an augmentation-preserving non-negative chain map \(\angles{c} \colon \chains(\bG^n) \to C\) with \(\angles{c}(e_n) = c\) if and only if for all \(\epsilon \in \set{-,+}\) and \(0 \leq k < n\), one has \(\angles{c}(e_k^\epsilon) = \angles{c}^\epsilon_k\) and \(\varepsilon \angles{c}_0^\epsilon = 1\).

If it exists, we refer to \(\angles{c}\) as the \defn{characteristic map} of the positive chain \(c\).
Characteristic maps of basis elements are called \defn{atoms}.
If all basis elements have characteristic maps, i.e., if \(\varepsilon\angles{b}^\varepsilon_0 = 1\) for all \(b \in B\) and \(\varepsilon \in \set{-,+}\), we say that the basis is \defn{unital}.

\SS The globular nerve naturally lifts along the forgetful functor:
\[
\begin{tikzcd}
	& \omega\Cat \dar \\
	\bach \rar \arrow[ru,dashed,bend left]& \bG\Set.
\end{tikzcd}
\]
The unit of a characteristic map \(\angles{c} \colon \chains(\bG^n) \to C\) is the map sending \(e_{n+1}\) to \(0\), both \(e_n^\pm\) to \(c\), and agrees with \(\angles{c}\) in degrees less than \(n\).
The composition of two characteristic maps \(\angles{c}\) and \(\angles{c'}\) with \(\ta_m\angles{c} = \so_m\angles{c'}\) is the map \(\angles{c} + \angles{c'} - \angles{c''}\) where \(c'' = \angles{c}(e_m^+) = \angles{c'}(e_m^-)\).

\SS\label{def:Steiner-diagram}
We introduce a relation \(<_k\) on basis elements of degree greater than \(k\) for \(k \in \N\) as follows:
\(b <_k b'\) if and only if there is a basis element appearing with nonzero coefficient in both \(\angles{b}^+_k\) and \(\angles{b'}^-_k\).
The basis is said to be \defn{loop-free} if for each \(k \geq 0\) the transitive closure of \(<_k\) is antisymmetric, i.e., if it defines a partial order.
An augmented based chain complex is said to be a \defn{Steiner diagram} if its basis is both unital and loop-free.
The full subcategory of Steiner diagrams is denoted by \(\wDiag\).
A property satisfied by the globular nerve of Steiner diagrams is that, if \(\angles{c}\) is a characteristic map, then \(c\) and all \(\angles{c}_k^\epsilon\) are sums of distinct basis elements (\cite[Thm.~4.1]{steiner2012opetopes}).

In many cases, instead of loop-freeness it is easier to check the following stronger condition.
We introduce a relation \(<_\N\) on all basis elements declaring \(b <_\N b'\) if and only if \(b \in \bd^-(b')\) or \(\bd^+(b) \ni b'\).
The basis is said to be \defn{strongly loop-free} if its transitive closure is antisymmetric.
As shown in \cite[Prop.~3.7]{steiner2004omega}, strongly loop-free bases are loop-free.
We refer to a Steiner diagram satisfying this condition as a \defn{strong Steiner diagram}.

\medskip The functor \(\cells \colon \wDiag \to \wCat\) is full and faithful, and the image of a Steiner diagram is an \(\omega\)-category freely generated by its atoms (\cite[Thm~5.11]{steiner2004omega}).

%% file: sec/poly_diag.tex

\subsection{Polytopes and diagrams}\label{s:polytopes_diagrams}

In this section, we construct a based augmented chain complex \(\chains(P,\beta)\) for any polytope \(P\) with an \(f\)-orientation \(\beta\), and show that if \(\beta\) is induced from a \(P\)-admissible frame, the complex \(\chains(P,\beta)\) is unital.
For a framed polytope \(P\) we identify the atoms \(\angles{F}_k^{-}\) and \(\angles{F}_k^{+}\) with \(\so_k(F)\) and \(\ta_k(F)\), and show that \(\chains(P, \beta)\) is a (strong) Steiner diagram, i.e., it is (strongly) loop-free, if and only if \(P\) has no (inhomogeneous) cellular loops.

\SS\label{d:associated_chain_complex}
For any polytope \(P\) with an \(f\)-orientation \(\beta\), its \defn{chain complex} \(\chains(P, \beta)\) is the augmented chain complex
\[
\Z \xla{\,\varepsilon} \Z\set{\faces[P][0]} \xla{\,\bd} \Z\set{\faces[P][1]} \xla{\,\bd} \Z\set{\faces[P][2]} \xla{\,\bd} \dotsb,
\]
generated by its faces, graded by dimension, with
\[
\bd(F) = \sum\ta(F) - \sum \so(F)\,.
\]
The conditions \(\bd \circ \bd = 0\) and \(\varepsilon \circ \bd = 0\) are straightforward to verify; see \cref{p:oriented-thin} for more details.

\SS\label{ss:atoms}
When \(P\) is framed we simply write \(\chains(P)\) instead of \(\chains(P,\beta)\) if \(\beta\) is the induced \(f\)-orientation (\cref{induced_facial_orientation}).
The atoms of \(\chains(P)\) have the following geometric description in terms of extended sources and targets (\cref{s:sk_tk}).

\medskip\lemma For any \(d\)-face \(F\) of a framed polytope
\[
\angles{F}_k^{-} = \sum \so_k(F) \quad \Big(\defeq \sum_{\mathclap{\ E \in \so_k(F)}}\ E\Big)
\qquad\text{and}\qquad
\angles{F}_k^{+} = \sum \ta_k(F)
\]
for any \(0 \leq k < d\).

\begin{proof}
	We proceed by induction on \(\ell = d - k\) focusing on extended sources since targets are handled similarly.
	For \(k = d-1\), we have \(\so_{d-1}(F) = \so(F)\) so there is nothing to prove.
	Now, let us prove the formula for \(k = d-\ell-1\) if
	\[
	\langle F \rangle_{d-\ell}^{-} = \textstyle\sum \so_{d-\ell}(F).
	\]
	Explicitly, we want to show that
	\[
	\bd^{-} \big(\textstyle\sum \so_{d-\ell}(F)\big) = \textstyle\sum \so_{d-\ell-1}(F)
	\]

	From~\cref{s:subdivision}, we have that \(\pi_{d-\ell}(\so_{d-\ell}(F))\) is a polytopal subdivision of \(\pi_{d-\ell}(F)\).
	In particular, any \((d-\ell-1)\)-face of the subdivision either belongs to the boundary or to exactly two \((d-\ell)\)-faces \(\pi_{d-\ell}(F_1)\) and \(\pi_{d-\ell}(F_2)\) with \(F_i\in \so_{d-\ell}(F)\). 	In the latter case, if \(G\) is a \((d-\ell-1)\)-face such that \(\pi_{d-\ell}(G) = \pi_{d-\ell}(F_i) \cap \pi_{d-\ell}(F_j)\), we note that the span of \(\pi_{d-\ell}(G)\) separates \(\pi_{d-\ell}(F_i)\) and \(\pi_{d-\ell}(F_j)\), so their normals point in opposite directions.
	Hence, \(G\) is in the source of one and the target of the other and \(G\) will not appear in \(\bd^{-} \sum\so_{d-\ell}(F)\).
	The remaining faces geometrically form the \((d-\ell-1)\)-boundary of \(F\).
	One such face appears with a negative sign in \(\bd \angles{F}^{-}_{d-\ell}\) if and only if, it is in the source of \(\pi_{d-\ell}(F)\).
	This is exactly a face in \(\so_{d-\ell-1}(F)\) according to \cref{s:sk_tk}.
\end{proof}

\SS\label{t:unital}\theorem
If \(P\) is a framed polytope then the basis of \(\chains(P)\) is unital.

\begin{proof}
	By \cref{ss:atoms}, the \(0\)-dimensional components of the atoms are given by the~\(0\)-source and target.
	That is, for any face~\(F\) of~\(P\) we have \(\langle F \rangle_0^{+} = \ta_0(F)\) and \(\langle F \rangle_0^{-} = \so_0(F)\).
	Since \(P\) is framed, there is a linear isomorphism \(\Lin \pi_1(F) \cong \Lin(v_1)\) between the linear span of the \(1\)-dimensional projection of \(F\) and the line spanned by the first vector of the frame \(v_1\).
	Thus, the \(0\)-source and target \(\so_0(F)\) and \(\ta_0(F)\) are the two vertices of the segment \(\pi_1(F)\), respectively selected by \(v_1\) and \(-v_1\) along the above isomorphism.
	Since both consist of a single vertex, we have \(\varepsilon(\ta_0(F)) = \varepsilon(\so_0(F)) = 1\) and \cref{ss:atoms} concludes the proof.
\end{proof}

The following example shows that \(\chains(P, \beta)\) is not necessarily unital if \(\beta\) is not induced from a \(P\)-admissible frame.

\begin{figure}[h!]
	\input{fig/non_unital}
	\caption{A polytope with an \(f\)-orientation whose chain complex is not unital.}
	\label{fig:not-unital}
\end{figure}
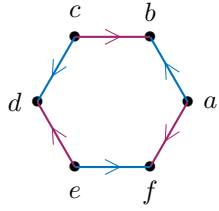

Consider the regular hexagon \(P\) illustrated in \cref{fig:not-unital}.
We give it the following \(f\)-orientation: for the top face, we set \(\beta_P \defeq e_1 \wedge e_2\), while for the edges we set the orientations as given by the arrows drawn.
Let us compute the atom~\(\langle P \rangle\) of~\(P\).
We have \(\langle P \rangle_1^{+} = \bd^{+}(P) = [ab]+[cd]+[ef]\), the blue edges in \cref{fig:not-unital}.
Applying the boundary map again, we get \(\langle P \rangle_0^{+} = \bd^{+}\langle P \rangle_1^{+} = b+d+f\), therefore we have \(\bd_0\langle P \rangle_0^{+} = 3 \neq 1\) and the unital condition does not hold.

\SS\label{t:loop-free}
\theorem A framed polytope \(P\) has no (inhomogeneous) cellular loops if and only if \(\chains(P)\) is a (strong) Steiner diagram.

\begin{proof}
	Since the basis of \(\chains(P)\) is unital by \cref{t:unital}, it suffices to show that \(P\) has no (inhomogeneous) cellular loops if and only if it is (strongly) loop-free.
	Let \(F_1,\ldots,F_m\) be a sequence of faces of \(P\).
	According to \cref{ss:atoms}, we have that \(F_i \in \bd^{-}(F_{i+1})\) (resp \(\bd^{+}(F_i) \ni F_{i+1}\)) if and only if \(F_i \in \so(F_{i+1})\) (resp \(\ta(F_i) \ni F_{i+1}\)).
	Therefore, \(F_1,\ldots,F_m\) is an inhomogeneous cellular loop (\cref{def:cellular-loop}) if and only if the relation \(<_\N\) fails to be antisymmetric.
	Since every cellular loop is part of an inhomogeneous one (\cref{l:inhomogeneous strings}), a similar claim holds for the relations~\(<_k\).
\end{proof}

%% file: fig/non_unital.tex
\begin{tikzpicture}

	\node (a) at (1, 0) {$\bullet$};
	\node at (1.3, 0) {$a$};

	\node (b) at (0.5, 0.866) {$\bullet$};
	\node at (0.5, 1.2) {$b$};

	\node (c) at (-0.5, 0.866) {$\bullet$};
	\node (c) at (-0.5, 1.2) {$c$};

	\node (d) at (-1, 0) {$\bullet$};
	\node at (-1.3, 0) {$d$};

	\node (e) at (-0.5, -0.866) {$\bullet$};
	\node at (-0.5, -1.2) {$e$};

	\node (f) at (0.5, -0.866) {$\bullet$};
	\node at (0.5, -1.2) {$f$};

	\draw[thick,RedViolet] (-0.5, 0.866)--(0.5, 0.866) node[midway,sloped,allow upside down,scale = 0.1]{\thickmidarrow};
	\draw[thick,RoyalBlue] (1, 0)--(0.5, 0.866) node[midway,sloped,allow upside down,scale = 0.1]{\thickmidarrow};
	\draw[thick,RoyalBlue] (-0.5, 0.866)--(-1, 0) node[midway,sloped,allow upside down,scale = 0.1]{\thickmidarrow};

	\draw[thick,RedViolet] (1,0)--(0.5, -0.866) node[midway,sloped,allow upside down,scale = 0.1]{\thickmidarrow};
	\draw[thick,RoyalBlue] (-0.5, -0.866)--(0.5, -0.866) node[midway,sloped,allow upside down,scale = 0.1]{\thickmidarrow};
	\draw[thick,RedViolet] (-0.5, -0.866)--(-1, 0) node[midway,sloped,allow upside down,scale = 0.1]{\thickmidarrow};
\end{tikzpicture}

%% file: sec/reformulation.tex

\subsection{The conjecture}\label{s:conjectures}

In this section, we revisit the Kapranov--Voevodsky conjecture in its original formulation and disprove it by establishing its connection to framed polytopes and cellular loops.

\SS Consider a polytope \(P \subset \R^n\) and a system of affine projections \(p\):
\begin{equation}\label{eq:system_of_projections}
	\R \xla{p_1^2} \R^2 \xla{p_2^3} \dotsb \xla{p_{n-2}^{n-1}} \R^{n-1} \xla{p_{n-1}^n} \R^n.
\end{equation}
We write \(p_k\) for the composition \(\R^n \to \R^k\), and say that \(p\) is \defn{\(P\)-admissible} if \(p_k\) is injective when restricted to any \(k\)-face of \(P\).
In this case, Kapranov and Voevodsky constructed the sources and targets associated to an \(f\)-orientation on \(P\) as follows.
Recall from \cref{ss:sotagiveforientation} that this data uniquely determines the \(f\)-orientation.
First, consider each \(\R^k\) in \eqref{eq:system_of_projections} to be canonically oriented, so that the fiber of \(p_k^{k+1}\) is also oriented.
For a facet \(E\) of a face \(F\), let \(u\) be an affine-linear functional which is zero on \(p_k(E)\) and non-negative on \(p_k(F)\).
Then, \(E \in \so(F)\) or \(E \in \ta(F)\) depending on whether \(u(t) \to -\infty\) or \(u(t) \to \infty\) as \(t \to \infty\) in the direction of the oriented fiber of \(p_k^{k+1}\).

\SS As the following result shows, framed polytopes capture all the complexity of this construction.

\medskip\lemma The system of projections of a framed polytope \(P\) is \(P\)-admissible and its induced \(f\)-orientation agrees with the one constructed by Kapranov--Voevodsky.
Furthermore, for any admissible system of projections on a polytope \(P\) there is a \(P\)-admissible frame inducing the same \(f\)-orientation.

\begin{proof}
	The admissibility of the system of projections of a \(P\)-admissible frame follows directly from the definition of \(P\)-admissibility.
	We need to show that the source of an arbitrary face \(F\) with respect to the \(f\)-orientation \(\beta\) induced by the given frame agrees with the source of \(F\) defined by the Kapranov--Voevodsky \(f\)-orientation \(\gamma\).
	According to \cref{lem:sources and targets in a framed polytope}, the \(\beta\) source of the \(k\)-face \(F\) consists of all \(E\), facets of \(F\), such that
	\[
	\sprod{\normal_{\pi_k E}^{\pi_k F}}{v_k} < 0.
	\]
	The \(\gamma\) source of \(F\) consists of all \(E\), facets of \(F\), for which there exists an affine-linear function \(u(w)\) which is zero on \(p_k(E)\) and non-negative on \(p_k(F)\) such that \(u(w) \to -\infty\) as \(w \to \infty\) in the oriented fiber of \(\pi_{k-1}^k\).
	The fiber of the projection \(\pi_{k-1}^k\) is generated and oriented by \(v_k\).
	The choice
	\[
	u(w) \defeq \max_{x \in P}\sprod{\normal_{\pi_k E}^{\pi_k F}}{x} - \sprod{\normal_{\pi_k E}^{\pi_k F}}{w}
	\]
	shows the agreement of both \(f\)-orientations.

	For the last claim, consider a \(P\)-admissible system of projections \(p\):
	\[
	\R \xla{p_1^2} \R^2 \xla{p_2^3} \dotsb \xla{p_{n-2}^{n-1}} \R^{n-1} \xla{p_{n-1}^n} \R^n.
	\]
	Let \(v_1\) be the canonical generator of \(\R\).
	Let us assume for an induction argument that \((v_1,\dots,v_k)\) have been chosen so that \(\Lin(v_1,\dots,v_j) = \R^j\) as oriented vector spaces and \(p_j^{j+1} = \pi_j^{j+1}\) for any \(1 \leq j < k\).
	Choose an oriented generator of \(\ker p_k^{k+1}\) and define it to be \(v_{k+1}\).
	We have \(\Lin(v_1,\dots,v_{k+1}) = \R^{k+1}\) as oriented vector spaces and \(p_k^{k+1} = \pi_k^{k+1}\), as claimed.
\end{proof}

\SS Using the previous lemma, we can restate the Kapranov--Voevodsky Conjecture in terms of framed polytopes.

\medskip\conjecture(\cite[Thm~2.3]{kapranov1991polycategory})
If \(P\) is a framed polytope then \(\chains(P)\) defines a pasting diagram.

\medskip The model of pasting diagrams cited by the authors is that of Johnson \cite{johnson1986pasting}, which, unfortunately, has certain shortcomings explained in \cite[\S1.5]{forest2022pasting}.
Specifically, Forest provides a refutation of its freeness property as stated in \cite[Theorem~13]{johnson1986pasting}.
However, the chain complex viewpoint used by Kapranov and Voevodsky raised the prospect of a formalization of pasting diagrams in these algebraic terms.
This was accomplished by Steiner in \cite{steiner2004omega}, who constructed a model for pasting diagrams based on chain complexes, a formalism which we reviewed in \cref{s:steiner}.

Recall from \cref{t:loop-free} that the chain complex of a framed polytope defines a Steiner diagram if and only if it has no cellular loops.
Thus, using this formalization of pasting diagrams, we obtain a complete characterization of the framed polytopes that satisfy the Kapranov--Voevodsky conjecture.
We anticipate that this result extends to any other model of pasting diagrams; see \cref{part:molecules} for an additional instance.
This notwithstanding, since all formalizations rule out looping behavior, the absence of cellular loops is a necessary condition for a framed polytope to define a pasting diagram, and the constructions in \cref{part:strings} serve as model-independent counterexamples to the conjecture.

Additionally, the example constructed in \cref{s:loop-inevitability} falsifies a weaker version of the conjecture that naturally arises from the existence of counterexamples to the original.
Explicitly, it shows that there are polytopes which fail to define a pasting diagram for any choice of frame.

\SS \remark It seems difficult to find a positive criterion ensuring that a framed polytope \(P\) is (strongly) loop-free.
For instance, the first counterexample in \cref{s:examples-loops} might suggest the following condition: for every \(k\), require that the projections \(\pi_k(v)\) of all vertices \(v\) of \(P\) lie in convex position.
Our second counterexample in \cref{s:examples-loops} shows that even if this condition could prevent the existence of \(1\)-loops, it does not prevent the existence of \(2\)-loops.

%% file: sec/orientals.tex

\section{Orientals}\label{part:orientals}

As demonstrated in \cref{part:strings}, cellular loops are common in framed polytopes.
For example, recall that the probability that a cellular loop exists in a canonically framed randomly embedded \(n\)-simplex tends to \(1\) as \(n \to \infty\) (\cref{s:random}).
In this part, we study two infinite families of canonically framed polytopes in which no cellular loop exists.
These families are obtained from certain special presentations of cubes and simplices, known as \textit{cyclic embeddings}.
We will prove that, up to certain signs, the associated pasting diagrams recover the simplicial and cubical orientals, as proposed by Kapranov and Voevodsky \cite{street1987orientals,  aitchison2010cubes}.
Their associated \(\omega\)-categories play a pivotal role in multiple areas of higher category theory.
In particular, they define the simplicial and cubical nerve of \(\omega\)-categories, offering an approach to weak higher categories as simplicial or cubical sets equipped with additional structure \cite{Verity2008a, Verity2008b, CampionKapulkin2025, DohertyKapulkin2023}.
Additionally, in the theory of higher Segal spaces \cite{dyckerhoff2012higher}, cyclic simplices provide a geometric model for formulating the higher Segal conditions, as detailed in \cite{Dyckerhoff2025}.

Although Street’s orientals were introduced and studied earlier, we will reverse this historical order in our presentation.
This is because most results in the simplicial theory can be deduced from corresponding statements in the cubical theory by considering vertex figures.

\input{sec/cubes}
\input{sec/simplices}

%% file: sec/cubes.tex

\subsection{Cyclic cubes}\label{ss:cubes}

\SS Consider the \defn{Veronese curve} \(\xi \colon \R \rightarrow \R^{d}\), given by \(\xi(t) \defeq (1, t, t^2, \dots, t^{d-1})\).
For \(n \geqslant d\), let \(t_{1}, \dots, t_{n}\) be elements of \(\R\) such that \(t_{1} < \dots < t_{n}\).
We denote the origin by~\(\mathbf{0}\) and by~\(\overline{\mathbf{0}\xi(t_{i})}\) the line segment from \(\mathbf{0}\) to \(\xi({t_{i}})\).
The \(d\)-dimensional \defn{cyclic zonotope} in \(\R^{d}\) is the Minkowski sum of line segments
\[
Z(n, d) \defeq \overline{\mathbf{0}\xi(t_{1})} + \dots + \overline{\mathbf{0}\xi(t_{n})}.
\]
The combinatorial type of the zonotope does not depend on the specific choice of the parameters \(\{t_1, \dots, t_n\} \subset \R\).
In view of this fact, let us set \(t_i\defeq i\) for all \(1 \leq i \leq n\) for simplicity.

\SS\lemma
For any cyclic zonotope, the canonical frame is \(Z(n,d)\)-admissible.

\begin{proof}
	Since any face of a cyclic zonotope is itself a cyclic zonotope, it is enough to prove admissibility for the top face.
	This follows from the identity \(\pi_{k}(Z(n,d)) = Z(n,k)\) for \(k \leq d\).
\end{proof}

\SS We review some combinatorial language used to represent the faces of \(Z(n,d)\).
Let \([n]\) denote the set \(\set{1,\ldots,n}\).
Every face \(F\) of \(Z(n,d)\) is given by a Minkowski sum
\[
F = \sum_{l \in L} \overline{\mathbf{0}\xi(l)} \ +\ \sum_{a \in A} \xi(a)
\]
for some subsets \(L \subseteq [n]\) and \(A \subseteq [n] \setminus L\).
We call \(L\) the set of \defn{generating vectors} of \(F\), and refer to the set \(A\) as its \defn{initial vertex}.
We identify \(F\) with the pair \((L, A)\).

\SS \label{sec:cycliczonotopefacets} We will now describe combinatorially the source and target of the canonically framed cyclic zonotope \(Z(n,d)\).
Given a subset~\(L \subseteq [n]\) and~\(l \in [n] \setminus L\), we say that~\(l\) is \defn{even} or \defn{odd} in \(L\) depending on whether the cardinality of the set \(\set{\ell \in L \mid \ell > l}\) is either even or odd.

\medskip \lemma \label{l:facets-zonotope}
Let \(L \subseteq [n]\) and \(A \subseteq [n] \setminus L\).
Then \((L,A)\) indexes a facet~\(F\) of the canonically framed cyclic zonotope~\(Z(n,d)\) if and only if it satisfies one of the following two conditions
\begin{enumerate}
	\item for all \(a \in A\), \(a\) is odd in \(L\), and for all \(b \in [n] \setminus (L \cup A)\), \(b\) is even in \(L\),
	\item for all \(a \in A\), \(a\) is even in \(L\), and for all \(b \in [n] \setminus (L \cup A)\), \(b\) is odd in \(L\).
\end{enumerate}
Moreover, \(F\) is in the source \(\so(Z(n,d))\) (resp.\ \(\ta(Z(n,d))\)) if and only if it satisfies~(1) (resp.~(2)).

\begin{proof}
	This is Gale's evenness criterion for the alternating oriented matroid, see \cite[Cor.~2.5]{Athanasiadis01} or \cite[Lem.~2.1]{Thomas02} where the terminology of lower and upper facets is used instead of source and target.
	See also \cite[Prop.~8.1]{dkk19}.
\end{proof}

\SS For any \(d \in \N\), the polytope \(Z(d,d) \subset \R^d\) is combinatorially a \(d\)-cube.
We will refer to it as the \defn{cyclic cube} and simplify its notation to \(Z(d)\).

The facets of \((L, A)\) come in parallel pairs.
Explicitly, for any \(l \in L\), both \((L \setminus \set{l}, A)\) and \((L \setminus \set{l}, A \cup \set{l})\) are facets of \((L, A)\).

\SS\label{prop:up_low_facets}\label{def:upper-lower}

We now extend \cref{sec:cycliczonotopefacets} to describe combinatorially the sources and targets of faces in the canonically framed cyclic cube.

\medskip\lemma\label{l:cube_boundary}
Let \(L \subseteq [n]\) and \(A \subseteq [n] \setminus L\) such that \(F = (L,A)\) is a face of the canonically framed cyclic \(d\)-cube, and let \(E\) be one of the two facets of \(F\) determined by some \(l \in L\).
\begin{enumerate}
	\item If \(E = (L \setminus l, A)\) then \(E \in \ta(F)\) if and only if \(l\) is odd in \(L \setminus l\).
	\item If \(E = (L \setminus l, A \cup \{l\})\) then \(E \in \ta(F)\) if and only if \(l\) is even in \(L \setminus l\).
\end{enumerate}
A similar characterization of \(\so(F)\) holds using the complementary parity conditions.

\begin{proof}
	For the top face of the cyclic cube, this is a direct application of \cref{l:facets-zonotope} in the case \(n=d\).
	The condition for arbitrary faces follows since every face of a cyclic cube is itself a lower-dimensional cyclic cube.
\end{proof}

\subsection{Cubical orientals}

\SS\label{ss:cubical_oriental}\label{thm:KV-cubes}

In \cite[Ex.~3.10]{steiner2004omega}, Steiner demonstrated that the usual tensor product on based augmented chain complexes induces a product structure on \emph{strong Steiner diagrams}.
This product, known as the \defn{Gray tensor product} \(\otimes\), plays a pivotal role in higher category theory, particularly in constructing the so-called \emph{cubical orientals}, which were introduced by I. Aitchison \cite{aitchison2010cubes}, R. Street \cite{street1991paritycomplexes}, and M. Johnson \cite{johnson1987thesis} using different formalizations.

Specifically, let \(\chains(\interval)\) denote the based augmented chain complex of the canonically framed interval, which forms a strong Steiner diagram satisfying \(\bd[0,1] = [1]-[0]\).
The \(d\)-fold Gray tensor product \(\chains(\interval)^{\otimes d}\) is then referred to as the \defn{cubical \(d\)-oriental}.

Kapranov and Voevodsky conjectured that the cubical orientals could be recovered by applying their procedure to the canonically framed cyclic cubes (\cite[Thm~2.7]{kapranov1991polycategory}).
In what follows, we leverage our convex-geometric criterion (\cref{t:loop-free}) to confirm this conjecture, subject to a slight adjustment of the framing.

\medskip\theorem
\label{thm:cyclic-cube}
The cyclic \(d\)-cube framed by \((e_1, -e_2,\dots, (-1)^{d-1} e_d)\) has no inhomogeneous cellular loops.
Furthermore, the resulting strong Steiner diagram is isomorphic to the cubical \(d\)-oriental.

\begin{proof}
	It suffices to establish an isomorphism of based augmented chain complexes between \(\chains(\interval)^{\otimes d}\), the cubical \(d\)-oriental, and \(\chains(Z(d))\), the based augmented chain complex of the cyclic \(d\)-cube framed by \((e_1, -e_2,\dots, (-1)^{d-1} e_d)\).

	The basis of \(\chains(\interval)^{\otimes d}\) in degree \(k\) consists of tensor products \(x = x_1 \ot\dotsb\ot x_d\) where each \(x_j\) is either \([0]\), \([1]\) or one of exactly \(k\) instances of \([0,1]\).
	Denoting by \(1 \leq p_1 < \dots < p_k < d\) the indices of \(x\) with \(x_{p_i} = [0,1]\), the decomposition of the boundary \(\bd x = \bd^+x - \bd^-x\) into positive elements---which is determined by the Koszul sign convention---can be described explicitly as follows:
	\[
	\bd^+x = \sum_{i = 1}^{k} x_1 \ot\dotsb\ot \bd^{\varepsilon_i} x_{p_i} \ot\dotsb\ot x_d \,,
	\]
	where \(\varepsilon_i = +\) (resp.\ \(-\)) if \(i\) is odd (resp.\ even), and \(\bd^-x\) is obtained by the same expression under the complementary conditions.

	Let us identify the above basis element \(x = x_1 \ot\dotsb\ot x_d\) in \(\chains(\interval)^{\otimes d}_k\) with the face \(F_x\) in \(Z(d)\) determined by the pair \((L,A)\) with \(L = \set{p_1,\dots,p_k}\) and \(A = \set{j \in [d] \mid x_j = [1]}\) as in \cref{def:upper-lower}.
	This defines an isomorphism of augmented graded abelian groups, which we must verify that it extends to one of chain complexes.
	Explicitly, that, under this isomorphism, \(\bd^+ x = \ta(F_x)\) and \(\bd^-x = \so(F_x)\).
	We will only write the argument for the first identity, as the second is proven analogously.

	The basis elements \(x_1 \ot\dotsb\ot \bd^{-} x_{p_i} \ot\dotsb\ot x_d\) and \(x_1 \ot\dotsb\ot \bd^{+} x_{p_i} \ot\dotsb\ot x_d\) correspond respectively to the faces \(\bd_i^{-} F_x \defeq (L \setminus \set{p_i}, A)\) and \(\bd_i^{+} F_x \defeq (L \setminus \set{p_i}, A \cup \set{p_i})\).
	Adapting \cref{l:cube_boundary}, which was proven for the canonical frame, we have that if \(k\) is odd:
	\begin{enumerate}
		\item \(\bd_i^{+} F_x \in \ta(F_x)\) if and only if \(p_i\) is even in \(L \setminus \set{p_i} \) if and only if \(i\) is odd.
		\item \(\bd_i^{-} F_x \in \ta(F_x)\) if and only if \(p_i\) is odd in \(L \setminus \set{p_i}\) if and only if \(i\) is even.
	\end{enumerate}
	The opposite equivalences hold for \(k\) even.
	We therefore see that \(\bd^+ x = \ta(F_x)\) and \(\bd^-x = \so(F_x)\), as desired.
	This finishes the proof.
\end{proof}

\SS \remark
The previous proof naturally raises the question of whether one could define geometrically a product of framed polytopes that induces the Gray tensor product at the chain level.
This seems to be an interesting problem.
Note that the Gray product has been defined combinatorially at the level of oriented graded posets \cite[Def.~7.2.1]{Hadzihasanovic24}, see \cref{part:molecules} below.

\medskip\corollary
The canonical basis is \(Z(n,d)\)-admissible for any cyclic zonotope and its chain complex is a strong Steiner diagram.

\begin{proof}
	Since \(Z(n)\) is both framed and strongly loop-free, by applying \cref{cor:projected-loop-free} we conclude that the projection \(\pi_d(Z(n)) = Z(n,d)\) is also framed and strongly loop-free for any \(d \leq n\).
\end{proof}

\subsection{Higher Bruhat orders and cubillages} \label{sec:HBO}\label{def-notation-A_K}\label{cor:cubical-orientals}

As part of their study of the Yang--Baxter equation, in \cite{ManinSchechtman89} Manin and Schechtman introduced a family of posets \(B(n,d)\) generalizing the weak Bruhat order on the symmetric groups \(\sym_n \cong B(n, 1)\), a poset structure based on the number of inversions of a permutation.
As anticipated in \cite{kapranov1991polycategory} and further elaborated in \cite{Thomas02} (see also \cite{dkk19}), the poset \(B(n,d)\) is isomorphic to the poset of \(\pi_d\)-induced subdivisions, where \(\pi_d \colon Z(n) \to Z(n, d)\) is the \(d\)-projection associated to the canonical frame.

In page 21 of the same reference, they anticipate that the minimal and maximal elements of \(B(n, d)\) correspond respectively to the \(d\)-source and \(d\)-target of the top-dimensional basis element of the cubical \(n\)-oriental.
A claim that we intend to make precise in this subsection.

\SS Let \(\binom{[n]}{d}\) be the set of cardinality \(d\) subsets of \([n] = \set{1,\ldots,n}\).
Given a set \(M \in \binom{[n]}{d+2}\), the \defn{packet} of \(M\) is the set \(P(M) \defeq \{ M \setminus k \mid k \in M \} \subseteq \binom{[n]}{d+1}\).
A natural total order on the elements of \(P(M)\) is the lexicographic order.
A subset \(\cU \subseteq \binom{[n]}{d+1}\) is \defn{consistent} if for all \(M \in \binom{[n]}{d+2}\), the intersection \(P(M) \cap \cU\) is either a beginning or ending segment of \(P(M)\) in the lexicographic order.

\SS The \defn{\((n,d)\)-higher Bruhat order} is the poset \(B(n,d)\) whose elements are consistent subsets \(\cU \subseteq \binom{[n]}{d+1}\) and with cover relations \(\cU \lessdot \cU'\) if and only if \(\cU' = \cU \cup K\) for some \(K \in \binom{[n]}{d+1} \setminus \cU\).

\SS Recall that given \(L \in \binom{[n]}{d+1}\) and \(a \in [n] \setminus L\), we say that \(a\) is an \defn{even gap} if
\(\left|\{\,l \in L \mid l > a\,\}\right|\) is even, and that \(a\) is an \defn{odd gap}\ otherwise.
Given a consistent set \(\cU\) and a set of generating vectors \(L \in \binom{[n]}{d+1}\), we define the associated \defn{vertex set} \(A_L^\cU\) as follows.
Given \(a \in [n] \setminus L\), we say that \(a \in A_L^\cU\) if and only if
\begin{itemize}
	\item \(L \cup \{a\} \in \cU\) and \(a\) is an even gap in \(L\), or
	\item \(L \cup \{a\} \notin \cU\) and \(a\) is an odd gap in \(L\).
\end{itemize}

\SS The following theorem is stated in the unpublished manuscript \cite[Thm.~2.1]{Thomas02}.
Unfortunately, the proof presented there has a gap.
Specifically, at one point, the argument needs that certain sequence of faces is acyclic.
The justification given boils down to that if \(F_1,F_2\) form a cellular string, then the topmost point of \(F_1\) is lower than that of \(F_2\).
However, this does not hold in general.
As we have shown, loops are not only possible but likely.
Therefore, an alternative proof of acyclicity in this particular case is needed, and we provide a justification below.
\medskip

\theorem \label{thm:bijection-HBO-cubillages}
The assignment
\begin{equation}
	\begin{matrix}
		\phi & : & B(n,d) & \longrightarrow & Z(n,d) \\
		& & \cU & \mapsto & \bigcup_{L \in \binom{[n]}{d+1}}(L,A_L^\cU)
	\end{matrix}
\end{equation}
is a bijection between elements of the higher Bruhat orders and cubillages of the cyclic zonotope \(Z(n,d)\).

\begin{proof}
	We sketch the main ideas of the proof of \cite[Thm.~2.1]{Thomas02}, and give details on how to prove acyclicity.
	To prove that \(\phi\) is well-defined and injective, one proceeds by recurrence on the dimension, and uses the fundamental theorem of the higher Bruhat orders \cite[Thm.~3]{ManinSchechtman89} to express an element \(\cU\) of \(B(n,d)\) as a maximal chain of elements in \(B(n,d-1)\).
	In higher categorical terms, this corresponds to the choice of a \emph{layering} of \(\cU\) (see \cref{def:layering} below).
	To prove that \(\phi\) is surjective, one proceeds similarly by induction on the dimension \(d\).
	For the inductive step, one starts with a cubillage \(\cU\) of \(Z(n,d)\), and then needs to show that it can be reconstructed from a sequence of cubillages of \(Z(n,d-1)\) as above.
	One starts with the source \(\so(Z(n,d))\), and then needs to find a \(d\)-cube in \(\cU\) whose source is entirely contained in \(\so(Z(n,d))\).
	This is indeed always possible, and it is ensured by the fact that the cyclic cube has no cellular loops (\cref{thm:cyclic-cube}).
	In higher categorical terms, this amounts to the fact that the cubical \(n\)-oriental is \emph{freely generated} in the sense of polygraphs (cf. \cite{Street1976, Burroni1991}).
	Note that the argument given in the proof of \cite[Thm.~2.1]{Thomas02} for acyclicity is not valid: one can easily find a cubillage of \(Z(4,2)\) with two cubes \(C_1, C_2\) such that \(C_1 <_1 C_2\) but such that the topmost point of \(C_1\) is above that of \(C_2\), see \cite[Fig.~19]{RambauReiner12}.
	Alternative proofs of acyclicity for cubillages of cyclic zonotopes can be found in \cite[Lem.~3.2]{Athanasiadis01} and \cite[Prop.~9.1]{dkk19}.
\end{proof}

\SS \corollary
The bijection between the higher Bruhat order \(B(n,d)\) and cubillages of the cyclic zonotope \(Z(n,d)\) identifies the minimal and maximal elements of \(B(n,d)\) with the \(d\)-source and \(d\)-target of the top-dimensional basis element of the cubical \(n\)-oriental.

\begin{proof}
	Using the combinatorial description from \cref{l:facets-zonotope}, one sees directly that the bijection \(\phi\) from \cref{thm:bijection-HBO-cubillages} sends extremal elements of \(B(n,d)\), which correspond to the consistent sets \(\cU = \emptyset\) and \(\cU = \binom{[n]}{d+1}\), to the source and target of \(Z(n,d)\), or equivalently to the \(d\)-source and \(d\)-target of the cyclic \(n\)-cube \(Z(n)\).
	\cref{thm:cyclic-cube} then gives the correspondence with the sources and targets of the cubical \(n\)-oriental.
\end{proof}

%% file: sec/simplices.tex

\subsection{Cyclic simplices}\label{ss:simplex}

\SS\label{def:cyclicpolytope}\label{l:cyclic-simplex}
Recall from \cref{ss:cyclic_loops} the definition of \(\CP(n,d) \subset \R^d\), the \(d\)-dimensional cyclic polytope on~\(n\) vertices:
\(\CP(n,d)\) is the convex hull of \(n\) distinct points on the moment curve \(\nu(t) = (t,t^2,\dots,t^d)\).

\SS \lemma For any cyclic polytope \(\CP(n,d)\), the canonical frame is \(\CP(n,d)\)-admissible.

\begin{proof}
	Since every face of a cyclic polytope is itself a cyclic polytope, it is enough to prove admissibility for the top face.
	This follows from the identity \(\pi_k(\CP(n,d)) = \CP(n,k)\) for \(k \leq d\).
\end{proof}

\SS We review some combinatorial language used to represent the faces of \(\CP(n,d)\).
Let \([n]\) denote the set \(\{1,\ldots,n\}\).
Every face \(F\) of \(\CP(n,d)\) is given by the convex hull of a subset of vertices \(\set{\nu(l)}_{l \in L}\).
We will identify a face \(F\) with its set of vertices \(L\).

\SS We will now describe combinatorially the source and target of the canonically framed cyclic polytope.
As before, we say that \(b \in [n] \setminus L\) is \defn{even} in \(L\) if the cardinality of the set \(\set{\ell \in L \mid \ell > l}\) is an even number, and \defn{odd} in \(L\) otherwise.

\medskip \lemma \label{l:facets-cyclic-polytope}
Let \(F\) be a face of the canonically framed cyclic polytope \(\CP(n,d)\).
Then \(F\) is a facet of \(\CP(n,d)\) if and only if it satisfies one of the following two conditions
\begin{enumerate}
	\item for all \(b \in [n]\setminus F\), \(b\) is even in \(F\),
	\item for all \(b \in [n]\setminus F\), \(b\) is odd in \(F\).
\end{enumerate}
Moreover, \(F\) is in the source \(\so(\CP(n,d))\) (resp.\ \(\ta(\CP(n,d))\)) if and only if it satisfies~(1) (resp.~(2)).

\begin{proof}
	This is Gale's evenness criterion for cyclic configurations, see for instance \cite[Cor.~6.1.9]{deloeraTriangulations2010} where the terminology of lower and upper facets is used instead of source and target.
\end{proof}

\SS For any \(d \in \N\), the polytope~\(\CP(d+1, d) \subset \R^d\) is combinatorially a \(d\)-simplex.
We will refer to it as the \defn{cyclic simplex} and simplify its notation to~\(\CP(d)\).

\SS We will now describe combinatorially the sources and targets of faces in the canonically framed cyclic simplex.

\SS \lemma
\label{prop:up_low_facets-simplex}
Let \(L\) be the generating set of a face \(F\) of \(\CP(d)\), and let \(E\) be a facet of \(F\) with generating set \(L \setminus l\) for some \(l \in L\).
We have
\begin{enumerate}
	\item \(E \in \ta(F)\) if and only if \(l\) is odd in \(L\setminus l\),
	\item \(E \in \so(F)\) if and only if \(l\) is even in \(L\setminus l\).
\end{enumerate}

\begin{proof}
	This is a direct application of \cref{l:facets-cyclic-polytope} in the case \(n=d+1\), using the fact that any face of a cyclic simplex is itself a cyclic simplex.
\end{proof}

\SS The cyclic \(d\)-simplex is closely related to the cyclic cube~\(Z(d)\) as shown by the following.

\medskip\lemma
The cyclic simplex \(\CP(d)\) is affinely isomorphic to the vertex figure of~\(Z(d)\), either at the minimal or maximal vertex.

\begin{figure}[h!]
	\input{fig/cyclic-polytopes}
	\caption{Vertex figures of the cyclic cube at the minimal and maximal vertices, and intersection with the hyperplane \(x_1 = 1\).}
	\label{fig:cyclic-polytopes}
\end{figure}
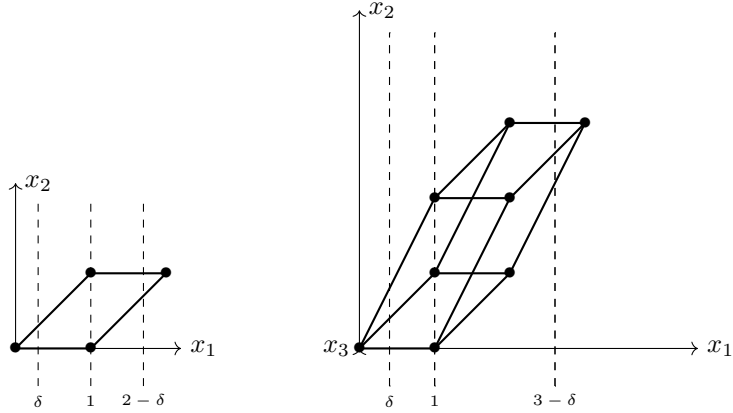

\begin{proof}
	The vertex figure of \(Z(d)\) at the minimal (resp.\ maximal) vertex is the intersection of \(Z(d)\) with a hyperplane \(x_1 = \delta\) (resp.\ \(x_1 = d-\delta\)) for a small \(\delta>0\), see \cref{fig:cyclic-polytopes}. Indeed, the intersection with \(x_1 = \delta\) is the convex hull of the points \(\delta \xi({t})\), which all share the same first coordinate.
	Removing this coordinate, we obtain a scaled copy of the cyclic simplex.
	Since \(Z(d)\) is centrally symmetric, the vertex figure of the maximal vertex is also affinely isomorphic to a cyclic simplex.
\end{proof}

\subsection{Street orientals}\label{s:street_orientals}

\SS \label{ss:streetorientals} In \cite[Ex.~3.8]{steiner2004omega}, Steiner shows that the following based augmented chain complex \(\chains(\triangle_d)\) is a strong Steiner diagram for each \(d \in \N\).
For \(0 \leq k \leq d\), let \(B_k\) be the set of all cardinality \(k+1\) subsets of \(\set{0,\dots,d}\), each denoted by the ordered tuple \([p_0,\dots,p_k]\) of its elements.
The \(\N\)-graded set \(B = \set{B_k}_{k \in \N}\) defines the basis of \(\chains(\triangle_d)\).
Its boundary is given by
\[
\bd\,[p_0,\dots,p_k] \defeq \sum_{i = 0}^k (-1)^i [p_0,\dots,p_{i-1},p_{i+1},\dots,p_k],
\]
and its augmentation by
\[
\varepsilon[p_0,\dots,p_k] \defeq
\begin{cases}
	1 & k = 0, \\
	0 & k \neq 0.
\end{cases}
\]
Furthermore, Steiner identified the free \(\omega\)-category generated by this diagram with the one constructed by Street in \cite{street1987orientals}.
%

\SS\label{l:simplicial-admissibility}\label{t:KV-simplex}

Let us now consider cyclic polytopes endowed with the frame
\[
\big(e_1, -e_2, \dots, (-1)^{d-1} e_d\big).
\]
This is just a reorientation of the canonical frame; therefore it is \(\CP(d)\)-admissible~(\cref{lem:reorientationequivalence}).

\medskip\theorem \label{thm:cyclic-simplex-oriental}
The cyclic \(d\)-simplex framed by \((e_1, -e_2,\dots, (-1)^{d-1} e_d)\) has no inhomogeneous cellular loops.
Furthermore, the resulting strong Steiner diagram is isomorphic to Street's \(d\)\textsuperscript{th} oriental.

\begin{proof}
	It suffices to establish an isomorphism of based augmented chain complexes between \(\chains(\triangle_d)\), Street's \(d\)-oriental, and \(\CP(d)\), the based augmented chain complex of the cyclic \(d\)-simplex framed by \((e_1, -e_2,\dots, (-1)^{d-1} e_d)\).

	The basis of \(\chains(\triangle_d)\) in degree \(k\) is given by ordered tuples of the form \(x=[p_0,\ldots,p_k]\).
	The decomposition of the boundary \(\bd x = \bd^{+} x - \bd^{-} x\) into positive elements can be described explicitly as follows:
	\[
	\bd^{+} x = \sum_{i \ \text{even}} [p_0,\dots,p_{i-1},p_{i+1},\dots,p_k].
	\]
	The other boundary component \(\bd^{-} x\) is given by the sum over \(i\) odd.

	Let us identify the above basis element \(x\) in \(\chains(\triangle_d)\) with the face \(F_x\) in \(\CP(d)\) determined by the set of vertices \(L=\set{p_0,\ldots,p_k}\) as in \cref{prop:up_low_facets-simplex}.
	This defines an isomorphism of augmented graded abelian groups, which we must verify that it extends to one of chain complexes.
	Explicitly, that, under this isomorphism, \(\bd^{+} x = \ta(F_x)\) and \(\bd^{-}x = \so(F_x)\).
	We will only write the argument for the first identity, as the second is proven analogously.

	The basis element \([p_0,\dots,p_{i-1},p_{i+1},\dots,p_k]\) corresponds to the face \(\bd_i F_x \defeq L \setminus \set{p_i}\).
	Adapting \cref{prop:up_low_facets-simplex}, which was proven for the canonical frame, we have that if \(k\) is odd: \(\bd_i F_x \in \ta(F_x)\) if and only if \(p_i\) is odd in \(L \setminus \set{p_i} \) if and only if \(i\) is even.
	The opposite equivalence holds for \(k\) even.
	We therefore see that \(\bd^+ x = \ta(F_x)\) and \(\bd^-x = \so(F_x)\), as desired.
	This finishes the proof.
\end{proof}

\SS\corollary
The canonical basis is \(\CP(n,d)\)-admissible for any cyclic polytope and its associated chain complex is a strong Steiner diagram.

\begin{proof}
	Since \(\CP(n)\) is both framed and strongly loop-free, by applying \cref{cor:projected-loop-free} we conclude that the projection \(\pi_d(\CP(n,n-1)) = \CP(n,d)\) is also framed and strongly loop-free for any \(d \leq n\).
\end{proof}

\subsection{Higher Stasheff--Tamari orders and triangulations}
\label{sec:HST}

The \emph{higher Stasheff--Tamari orders} are a family of posets \(S(n,d)\) introduced by Kapranov--Voevodsky in \cite{kapranov1991polycategory} generalising the Tamari order \(S(n,2)\) on triangulations of a convex \(n\)-gon.
They were the subject of much later study, see for instance the survey \cite{RambauReiner12} or the reference \cite[Sec.~6.1.4]{deloeraTriangulations2010}.
There are in fact two definitions of these orders, given by Edelman and Reiner in \cite{EdelmanReiner96}, which were only recently shown to coincide \cite{WilliamsTwoHigherStasheff2024}.

\SS The \((n,d)\)-\defn{higher Stasheff--Tamari order} is the poset \(S(n,d)\) whose elements are triangulations \(\cT\) of \(\CP(n,d)\) with cover relations \(\cT \lessdot \cT'\) if and only if \(\cT'\) can be obtained from \(\cT\) by an increasing bistellar flip.

Equivalently (cf. \cite{WilliamsTwoHigherStasheff2024}), we have \(\cT \lessdot \cT'\) if and only if the associated sections of the projection \(\CP(n,n-1) \to \CP(n,d)\) are such that \(\cT\) is weakly below \(\cT'\) with respect to the last coordinate.

\SS There are two preferred triangulations of \(\CP(n,d)\), given by the source and target of \(\CP(n,d+1)\).
These are called lower and upper triangulations in the literature, and define unique minimal and maximal elements of the poset \(S(n,d)\).

Since higher Stasheff--Tamari orders are defined directly in terms of triangulations, there is no analogue to bijection of the type \cref{thm:bijection-HBO-cubillages} to prove, and the following corollary is immediate.
However, there is an analogue of the fundamental theorem of the higher Bruhat orders \cite[Thm.~1.1]{Rambau1997} and an alternative proof to \cref{thm:cyclic-simplex-oriental} for acyclicity of triangulations is given in \cite[Cor.~5.9]{Rambau1997}.

\medskip \corollary The isomorphism from \cref{thm:cyclic-simplex-oriental} between the canonically framed cyclic \(n\)-simplex and Street's \(n\)-oriental identifies the minimal and maximal elements of the higher Stasheff--Tamari order \(S(n,d)\) with the \(d\)-source and \(d\)-target of the top-dimensional basis element of the \(n\)-oriental.

\medskip \remark Questions regarding maps between the higher Bruhat and higher Stasheff--Tamari orders are quite subtle \cite{Thomas02}, and some remain open \cite[Open~problem~8.10]{RambauReiner12}.

%% file: fig/cyclic-polytopes.tex
\begin{tikzpicture}
	\draw[->] (0, 0) -- (2.2, 0) node[anchor = west]{$x_1$};
	\draw[->] (0, 0) -- (0, 2.2) node[anchor = west]{$x_2$};
	\node (a) at (0, 0) {$\bullet$};
	\node (b) at (1, 0) {$\bullet$};
	\node (d) at (1, 1) {$\bullet$};
	\node (c) at (2, 1) {$\bullet$};
	\node at (0.3, -0.7) {\tiny $\delta$};
	\node at (1, -0.7) {\tiny $1$};
	\node at (1.7, -0.7) {\tiny $2-\delta$};
	\draw[-,thick] (0,0)--(1,0)--(2,1)--(1,1)--(0,0);
	\draw[-,dashed] (1.7, -0.5) -- (1.7, 2);
	\draw[-,dashed] (0.3, -0.5) -- (0.3, 2);
	\draw[-,dashed] (1, -0.5) -- (1, 2);
\end{tikzpicture}
\quad \quad
\tdplotsetmaincoords{0}{0}
\begin{tikzpicture}
	[tdplot_main_coords,
	cube/.style = {very thick,black},
	grid/.style = {very thin,gray},
	axis/.style = {->}]

	\node (a) at (0,0,0) {$\bullet$};
	\node (b) at (1,0,0) {$\bullet$};
	\node (c) at (1,1,1) {$\bullet$};
	\node (d) at (1,2,4) {$\bullet$};
	\node (bc) at (2,1,1) {$\bullet$};
	\node (bd) at (2,2,4) {$\bullet$};
	\node (cd) at (2,3,5) {$\bullet$};
	\node (bcd) at (3,3,5) {$\bullet$};

	\node at (0.4, -0.7,0) {\tiny $\delta$};
	\node at (1, -0.7,0) {\tiny $1$};
	\node at (2.6, -0.7,0) {\tiny $3-\delta$};

	\draw[axis] (0,0,0) -- (4.5,0,0) node[anchor = west]{$x_1$};
	\draw[axis] (0,0,0) -- (0,4.5,0) node[anchor = west]{$x_2$};
	\draw[axis] (0,0,0) -- (0,0,7.5) node[anchor = east]{$x_3$};

	\draw[-,thick] (0,0,0)--(1,0,0);
	\draw[-,thick] (0,0,0)--(1,1,1);
	\draw[-,thick] (0,0,0)--(1,2,4);

	\draw[-,thick] (2,1,1)--(1,0,0);
	\draw[-,thick] (2,1,1)--(1,1,1);

	\draw[-,thick] (2,2,4)--(1,0,0);
	\draw[-,thick] (2,2,4)--(1,2,4);

	\draw[-,thick] (2,3,5)--(1,1,1);
	\draw[-,thick] (2,3,5)--(1,2,4);

	\draw[-,thick] (2,1,1)--(3,3,5);
	\draw[-,thick] (2,3,5)--(3,3,5);
	\draw[-,thick] (2,2,4)--(3,3,5);

	\draw[-,dashed] (2.6,-0.5,-0.5) -- (2.6,4.2,-0.5);
	\draw[-,dashed] (2.6,-0.5,-0.5) -- (2.6,-0.5,6.2);
	\draw[-,dashed] (2.6,4.2,-0.5) -- (2.6,4.2,6.2);
	\draw[-,dashed] (2.6,-0.5,6.2) -- (2.6,4.2,6.2);

	\draw[-,dashed] (0.4,-0.5,-0.5) -- (0.4,4.2,-0.5);
	\draw[-,dashed] (0.4,-0.5,-0.5) -- (0.4,-0.5,6.2);
	\draw[-,dashed] (0.4,4.2,-0.5) -- (0.4,4.2,6.2);
	\draw[-,dashed] (0.4,-0.5,6.2) -- (0.4,4.2,6.2);

	\draw[-,dashed] (1,-0.5,-0.5) -- (1,4.2,-0.5);
	\draw[-,dashed] (1,-0.5,-0.5) -- (1,-0.5,6.2);
	\draw[-,dashed] (1,4.2,-0.5) -- (1,4.2,6.2);
	\draw[-,dashed] (1,-0.5,6.2) -- (1,4.2,6.2);
\end{tikzpicture}

%% file: sec/matroids.tex

\section{Framed simplices and oriented matroids}\label{part:matroids}

In this part, we demonstrate a close relationship between framed simplices and oriented matroids.
Specifically, we establish that the $f$-orientations of framed standard simplices~$(\ssimplex[d], B)$ correspond bijectively with the chirotopes of the point configurations~$\pi_k(\ssimplex[d])$ (\cref{lem:diagramsareflagoms}).
In more advanced terms, this implies that the $f$-orientations of framed simplices are in bijection with uniform acyclic realizable full flag chirotopes (\cref{thm:globular-flag-chirotopes}).
We will leverage this connection in the study of moduli spaces of frames in \cref{part:universality}.

\subsection{Chirotopes and oriented matroids}\label{s:matroids}

Oriented matroids are combinatorial objects abstracting geometric information of vector and point configurations and hyperplane arrangements.
While the ordinary matroid associated to a vector configuration is determined by the subsets forming a basis, the associated oriented matroid also stores the orientations of the bases.
Here we will only need to work with certain oriented matroids that are associated with vector and affine point configurations, which are in fact rare among all oriented matroids.
Therefore, to simplify the exposition, instead of giving the proper combinatorial definition of oriented matroids, we will only present realizable oriented matroids as certain equivalence classes of vector configurations.
We refer to \cite{BLSWZ} for a comprehensive reference on the topic.

\SS\label{def:chirotope}
The \defn{chirotope} associated to a vector configuration $V = (v_1,\dots,v_n)\in \R^{d \times n}$ is the map
\begin{align*}
	\chiro[V] \colon \set{1,\dots,n}^d & \to \{+,-,0\} \\
	(i_1,\dots,i_d) & \mapsto \sign(\det(v_{i_1},\dots,v_{i_d})).
\end{align*}
While the chirotope is defined on $\set{1,\dots,n}^d$, it is completely determined from~$\binom{n}{d}$ values, as the rest can be deduced from its alternating properties.

\SS A general \defn{chirotope} is defined as a non-zero alternating map $\chi \colon \set{1,\dots,n}^d \to \{+,-,0\}$ satisfying the \defn{chirotope axioms} based on the Grassmann-Pl\"ucker relations \cite[Def.~3.5.3]{BLSWZ}.
We will however only treat those obtained from vector configurations as above, which are called \defn{realizable}.

\SS\label{def:chirotope-point}
If $A = (a_1,\dots,a_n) \in \R^{d \times n}$ is a point configuration, then its \defn{chirotope} is the chirotope of its homogenization $\hat A = (\hat a_1,\dots,\hat a_n)\in \R^ {(d+1) \times n}$, where $\hat a = \binom{1}{a}$ is obtained by appending a~$1$ as a new coordinate. Note that computing the determinant of a homogenization is equivalent to computing the determinant of the differences of non-homogenized points:
\[\det(\hat a_{i_0},\dots,\hat a_{i_d}) = \det(a_{i_1}- a_{i_0},\dots,a_{i_d}- a_{i_0}).\]
Realizable chirotopes obtained this way are called \defn{acyclic}, because for oriented matroids arising from graphs, they correspond to acyclic orientations.
To avoid having too much notation, we will use $\chiro[A]$ to denote $\chiro[\hat A]$ when $A$ is a point configuration, which will always be clear from the context.

\SS A chirotope~$\chi$ is called \defn{uniform} if $\chi(i_1,\dots,i_d) \neq 0$ whenever $i_1,\dots,i_d$ are pairwise distinct.
A point configuration whose chirotope is uniform is in affine \defn{general position}, i.e., no $d+1$ points lie in a common affine hyperplane.

\SS\label{ss:chirovsom}
A realizable chirotope depends on a frame for the ground vector space, as an orientation-reversing change of basis results in a global sign change for the chirotope.
An \defn{oriented matroid} is an equivalence class $\pm\chi = \{\chi,-\chi\}$ of chirotopes up to global reorientation \cite[Prop.~3.5.2 and Thm.~3.5.5]{BLSWZ},
where $-\chi$ denotes the chirotope obtained from $\chi$ by reversing all the signs: $+ \mapsto -$, $- \mapsto +$, and $0 \mapsto 0$.
Despite this subtle difference, the two terms \emph{chirotope} and \emph{oriented matroid} are often used interchangeably in the literature.


\subsection{Chirotopes and framed simplices}\label{s:flagoms}

When we restrict to framed standard simplices~$(\ssimplex[d],B)$,
the relation between $f$-orientations and chirotopes is quite satisfying as the next statements show.


\SS\label{lem:diagramsareflagoms}
\lemma Let $(\ssimplex[d],B)$ be a framed standard $d$-simplex.
The $f$-orientation of $(\ssimplex[n-1],B)$ determines and is determined by the chirotopes {in the basis $B = (v_1,\ldots,v_n)$} of the point configurations $\pi_k(\ssimplex[n-1]) = (\pi_k(e_1),\dots,\pi_k(e_{d+1}))$ for all~$0 \leq k \leq d$.

\begin{proof}
	Consider the $f$-orientation $\beta$ of $(\ssimplex[d],B)$.
	Let $F = \conv(e_{i_0},e_{i_1},\dots,e_{i_k})$ be a $k$-face of $\ssimplex[d]$, and consider the family of vectors $w_j \defeq e_{i_j}-e_{i_0}$ for $1\leq j\leq k$.
	Then $ \bigwedge^k \Lin_F$ is generated by $w_1\wedge \cdots \wedge w_k$. Therefore, we have
	\[
	\beta_F = \sigma^F_k v_1 \wedge\dots\wedge \sigma^F_k v_k = \lambda w_1\wedge \cdots \wedge w_k,
	\]
	for some~$\lambda\in\R$, and the equivalence class of the orientation of $\Lin_F$ is determined by the sign of~$\lambda$.
	Now, by applying the projection $\pi_k$, we have
	\begin{equation}\label{eq:detwedge}
	\lambda (\pi_k w_1\wedge \cdots \wedge \pi_k w_k) = v_1 \wedge\dots\wedge v_k,
	\end{equation}
	with the same coefficient $\lambda$.
	We claim that the sign sign of $\lambda$ determines the sign of the determinant \[\det(\widehat{\pi_{k} e_{i_0}},\dots,\widehat{\pi_{k} e_{i_k}}) = \det(\pi_{k} e_{i_1}-\pi_{k} e_{i_0},\ldots,\pi_{k} e_{i_k}-\pi_{k} e_{i_0})\] corresponding to the subset $\{i_0,\dots,i_k\}$ in the definition of the chirotope of $\pi_{k}(\ssimplex[d])$ in the basis $B$.
	Indeed, using \eqref{eq:detwedge} and the exterior product definition of the determinant, we have that $\det(\pi_{k} e_{i_1}-\pi_{k} e_{i_0},\ldots,\pi_{k} e_{i_k}-\pi_{k} e_{i_0}) = \det(\pi_k w_1, \ldots, \pi_k w_k)$ in the basis $B$ is~$\frac{1}{\lambda}$.
	Thus the orientations of the $k$-faces can be determined by the chirotope of $\pi_{k}(\ssimplex[d])$, and reciprocally.
\end{proof}

\SS Flag matroids, introduced in \cite{BorovikGelfandWhite2003}, are a generalization of matroids, and admit an oriented version.
A \defn{flag chirotope} is a sequence $(\chi_1,\dots,\chi_s)$ of chirotopes related by strong maps, also called quotients, see \cite[Ex.~above Thm.~D]{JarraLorscheid2022} and \cite[Def.~4.1]{BoretskyEurWilliams2022}, and also \cite[Def.~3.5.3, Thms~3.5.5 and~3.6.2, and Def.~7.7.2]{BLSWZ} for more details on the definition and the relation with ordinary oriented matroids.
As before, we will not detail the full definition, as we are only interested in the particular case of realizable oriented flag matroids, which are those induced by vector configurations.

\SS\remark
In the literature, flag chirotopes are usually called \emph{oriented flag matroids}.
However, we think that the name flag chirotopes is more precise, in view of the (subtle) difference between the classical definitions (cf.~\cref{ss:chirovsom}).
We prefer to reserve the name oriented flag matroids to equivalence classes of flag chirotopes under reorientation of its components.

\SS\label{def:realizable-flag}
Let $B$ be a frame of $\R^d$, with associated system of projections $\pi_k:\R^d\to V_k$.
A \defn{realizable complete flag chirotope} is a sequence of chirotopes $(\chi_0,\dots,\chi_d)$, where $\chi_k$ is the chirotope of the vector configuration $\{\pi_k(e_1),\dots, \pi_k(e_d)\}$.
See \cite[Sec 1.7.5]{BorovikGelfandWhite2003}
A flag chirotope $(\chi_0,\dots,\chi_d)$ is \defn{uniform} (resp. \defn{acyclic}) if $\chi_k$ is uniform (resp. acyclic) for all~$0\leq k\leq d$.

\SS\label{thm:globular-flag-chirotopes}
\theorem
\(f\)-orientations induced by framed simplices are in bijection with uniform acyclic realizable full flag chirotopes.

\begin{proof}
	This is just a reformulation of \cref{lem:diagramsareflagoms} in the language of \cref{def:realizable-flag}.
\end{proof}

%% file: sec/modulispaces.tex

\section{Universality of $P$-equivalence classes}\label{part:universality}

Consider a fixed $d$-polytope $P \subset \R^d$, and let $\frames$ be the set of frames of $\R^d$ that are $P$-admissible.
We will consider $\frames$ as a subset of $\GLn$ and parametrize it as a subset of $\R^{d \times d}$.
The goal of this part is to show (\cref{thm:simplex-universality}) that equivalence classes of $P$-admissible frames defining equivalent $f$-orientations are universal in the sense of Mn\"ev~\cite{Mnev1988}.
This roughly means that they can behave as ``badly'' as any semialgebraic set or, in the language of \cite{Vakil2006}, that they satisfy Murphy's law in Algebraic Geometry.


\subsection{Realization spaces of oriented matroids}

We start by recalling some classical results on the equivalence classes of point configurations under the oriented matroid equivalence.
In the oriented matroid literature, they are known under the name of \emph{realization spaces}.

\SS\label{def:rsom}

Let $\chi$ be a chirotope.
Its \defn{realization space} is the set $\rs[\chi]\subset \R^{d\times n}$ of point configurations whose chirotope is $\pm\chi$, up to affine transformation:
\[
\rs[\chi] \defeq \bigslant{\set{A\in \R^{d\times n} \mid \chiro[A] = \pm\chi}}{\Affn[d]}.
\]
If $A = (a_1,\dots,a_n)\in \R^{d\times n}$ is a point configuration, we denote by $\rs[A] \defeq \rs[\chiro]$ the realization space of its
chirotope.

\SS

To mod out by the set of affine transformations, we will fix the coordinates of an affine basis ($d+1$ affinely independent points) by choosing the origin and the~$d$ canonical basis vectors in all realizations (c.f.\
\cite[Def.~2, p.~3]{RichterGebert1997}).
This specific choice is necessary but not important, as changing the affine basis gives rationally equivalent realization spaces~\cite[Lem.~2.5.4]{RichterGebert1997}.


\subsection{Mn\"ev's Universality}

A groundbreaking result of Mn\"ev, the Universality Theorem~\cite{Mnev1988} (see also~\cite{RichterGebert95,RichterGebert1999}), states that realization spaces can be as complicated as arbitrary semi-algebraic sets, modulo stable equivalence.
We need to introduce some concepts to give the precise statement.

\SS

A \defn{primary basic semi-algebraic} set is a subset $S$ of $\R^d$ defined by integer polynomial equations and strict inequalities
\[
S = \set{\mathbf{x}\in\R^d \mid f_1(\mathbf{x}) = 0, \dots, f_k(\mathbf{x}) = 0, f_{k+1}(\mathbf{x})>0, \dots, f_r(\mathbf{x})>0},
\]
where each $f_i \in \Z[\mathbf{x}]$.
Realization spaces of oriented matroids are examples of primary basic semi-algebraic sets.
They are cut by determinantal equations and inequalities, and the quotient by affine transformations can be done by fixing the coordinates of an affine basis.

\SS

Two semi-algebraic sets $S,S'$ are \defn{stably equivalent}, denoted $S \sim S'$, if they lie in the same equivalence class generated by stable projections and rational equivalence.
Here, a projection $\pi \colon S\to S'$ is called \defn{stable} if its fibers are relative interiors of non-empty polyhedra of the same dimension defined by polynomial functions on $S'$, see~\cite[Sec.~2.5]{RichterGebert1999} for details, and \cite{Verkama2023} for the constant dimension constraint.

\SS\label{thm:universality}

\textit{Theorem} (Universality theorem \cite{Mnev1988}).
For every primary basic semi-algebraic set $S$ defined over $\Z$ there is a $2$-dimensional point configuration whose realization space is (as an oriented matroid) stably equivalent to~$S$.
If moreover $S$ is open, then the oriented matroid may be chosen to be uniform.


\subsection{Universality of $P$-equivalence classes}

We now aim to prove our main result in this section, \cref{thm:simplex-universality}, which establishes that $\ssimplex[d]$-equivalence classes of $\ssimplex[d]$-admissible frames, and thus realization spaces of flag chirotopes, are universal.
This will be done by reducing it to the universality theorem for oriented matroids.

\SS\label{def:rsframe}

The \defn{realization space of a frame} $B$ of~$\R^n$ with respect to a polytope $P\subseteq \R^n$ is the set~$\msp[B][P]$ of $P$-admissible frames that are $P$-equivalent to $B$.
More generally, we define the \defn{realization space of an $f$-orientation} $\gamma$ of~$P$ as the set~$\ms[\gamma]$ of $P$-admissible frames inducing~$\gamma$.
We will only consider the case where $P$ is the standard simplex~$\ssimplex[d]$, and we will abbreviate $\msp[B][{\ssimplex[d]}]$ to $\ms[B]$ when $B$ is a frame of $\R^d$.

\SS The proof of universality of $P$-equivalence classes is a little technical, as it involves the manipulation of the realization spaces.
The main ideas are the following.
The starting point is \cref{lem:diagramsareflagoms}, which implies that given a framed simplex $(\ssimplex[d],B)$, we can read the chirotope of $\pi_k(\ssimplex[d])$ from the $f$-orientation of~$(\ssimplex[d], B)$.
We are particularly interested in the case $k = 2$ because it is the first dimension presenting universality for oriented matroids.
We know thus that there is a natural projection from the realization space of the frame~$\ms[B]$ to the realization space of the oriented matroid~$\rs[{\pi_2(\ssimplex[d])}]$.
Our goal will be to find a convenient section for this map.

More precisely, we will start from the universality theorem, which says that for any semialgebraic set~$S$ there is a planar point configuration~$A$ whose realization space encodes it.
Our goal will then be to find a framed simplex $(\ssimplex[d],B)$ such that~$\pi_2(\ssimplex[d]) = A$ and the projection $\ms[B] \to \rs[{\pi_2(\ssimplex[n])}]$ is stable.
This will imply that~$\ms[B]$ is stably equivalent to~$S$, and universality of $\ssimplex[d]$-equivalence classes will follow.

In terms of flag chirotopes, we will start with a rank~$3$ chirotope whose realization space is stably equivalent to~$S$, and we will complete it to a complete flag chirotope in a ``stable way''.

\SS As our focus is on framed simplices $(\ssimplex[d],B)$, in the study of moduli spaces of frames, it is sometimes helpful to do a linear transformation that maps~$B$ to the canonical basis (that is, to do a change of basis to the basis~$B$).
Note that, by~\cref{ss:invariancelinearautomorphism}, this does not affect the $f$-orientation.
This way, instead of studying the behavior of the $f$-orientations of a standard simplex under a varying basis, we study the $f$-orientations of a varying simplex under the canonical basis.

\SS For the technical part of the proof, we will need to work with a slight variant of the definition of realization space given above in \cref{def:rsom}: we will need point configurations that have points at infinity, in the projective sense.

Let $A = (a_1,\dots,a_n)\in \R^{d\times n}$ be a point configuration, and $V = (v_1,\dots,v_m)\in\R^{d\times m}$ be a vector configuration.
We define $\chiro[A,V]$ to be the chirotope of the vector configuration
\[\left(\tbinom{1}{a_1},\dots, \tbinom{1}{a_n},\tbinom{0}{v_1},\dots,\tbinom{0}{v_m}\right).\]
That is, points in $A$ are homogenized as usual, while the vectors in $V$ have a $0$ appended instead of a~$1$.
If $V = \{v\}$ is a singleton, we will abbreviate $\chiro[A,v]$.

Let $\chi \colon [n+1]^{d+1} \to \{+,-,0\}$ be a chirotope.
We define the \defn{extended realization space} of $\chi$ as the set $\rsinf[\chi]$ of all configurations $(A,a_{n+1}) = (a_1,\dots,a_{n},a_{n+1})\in \R^{d\times (n+1)}$ for which $\chi = \chiro[A,a_{n+1}]$, up to affine transformation (quotient done by fixing an affine basis).
Geometrically, this has to be thought of as the set of realizations of~$\chi$ in which the last point has been sent to the line at infinity.

\SS\label{lem:rsandrsinftystablyequivalent}
For the following lemma, let us recall that whether a point is a vertex of the convex hull of a point configuration can be read off the chirotope.
In fact, this is the case for the whole face lattice of the convex hull, see for example~\cite[Ch.~4]{BLSWZ}.

\medskip\lemma
Let $\chi \colon [n+1]^{d+1}\to \{+,-\}$ be the uniform chirotope of a point configuration such that the last point is a vertex of their convex hull.
Then, the realization spaces~$\rs[\chi]$ and~$\rsinf[\chi]$ are stably equivalent. That is, we have $$\rs[\chi]\sim \rsinf[\chi].$$

\begin{proof}
	As standard, to mod out by the set of affine transformations, we fix the coordinates of an affine basis (c.f.
	\cite[Def.~2, p.~3]{RichterGebert1997}): we impose that $a_1,\dots,a_{d+1}$, which form an affine basis since the chirotope is uniform, are the origin and the~$d$ canonical basis vectors in all realizations.

	The strategy is to define an auxiliary set $\cC\subset \rs[\chi]\times \R^d$, and show that $\rs[\chi] \sim \cC \sim \rsinf[\chi]$.
	This auxiliary set $\cC$ consists of pairs $(A,c)$ where $A$ is a point configuration which realizes $\chi$ and $c$ is a normal covector for a hyperplane supporting $a_{n+1}$ in the convex hull of~$A$, meaning that $\sprod{c}{a_i}<\sprod{c}{a_{n+1}}$ for any~$i\leq n$.
	We have that the projection $\cC\to \rs$ with $(A,c)\mapsto A$ is surjective.
	As the open normal cone of a vertex is an open convex cone defined by polynomial equations on~$A$, this projection is stable and we have~$\cC\sim \rs[\chi]$.
	It remains to show that $\cC\sim \rsinf[\chi]$, from which the result will follow.

	We construct a new point configuration~$\tilde A = (\tilde a_1,\ldots,\tilde a_{n+1})$
	via a projective transformation $\phi: a_i\mapsto \tilde a_i$ sending the hyperplane supporting~$a_{n+1}$ in direction~$c$ to infinity.
	See~\cite[App.~2.6]{Ziegler95} for more on projective transformations in this context.
	More precisely, we start from the matrix
	\[M \defeq \begin{pmatrix}
		\sprod{c}{a_{n+1}}-\sprod{c}{a_{2}}&0&\cdots&0\\
		0& \sprod{c}{a_{n+1}}-\sprod{c}{a_{3}}&\cdots&0\\
		\vdots & \vdots & \ddots & \vdots \\
		0& 0&\cdots &\sprod{c}{a_{n+1}}-\sprod{c}{a_{d+1}}
	\end{pmatrix}
	\]
	whose entries $\sprod{c}{a_{n+1}}-\sprod{c}{a_{i}}$, $2\leq i\leq d+1$ are all positive because $c$ is a normal covector. In particular, we have $\det(M)>0$.
	The points of the new configuration~$\tilde A$ are then defined by
	$$\tilde a_i \defeq
	\begin{cases}
		\frac{M\cdot a_i}{\sprod{c}{a_{n+1}}-\sprod{c}{a_i}} & \text{ if } i\leq n, \\
		M\cdot a_{i} & \text{ if }i = n+1.
	\end{cases}$$
	We claim that $\tilde A$ is a realization in~$\rsinf[\chi]$.
	Indeed, at the level of the determinants defining the chirotope, the change from $\left( \binom{1}{a_1},\dots, \binom{1}{a_n}, \binom{1}{a_{n+1}}\right)$ to $\left( \binom{1}{\tilde a_1},\dots, \binom{1}{\tilde a_n}, \binom{0}{\tilde a_{n+1}}\right)$ does not alter the signs.
	To see this, observe that we multiply first the points of $A$ by the matrix
	\[\tilde M \defeq
	\begin{pmatrix}
		\sprod{c}{a_{n+1}}&		-c^\top \\
		0_d &M
	\end{pmatrix},
	\]
	whose determinant is $\det(\tilde M) = \sprod{c}{a_{n+1}}\det(M)$.
	Since we chose $a_1 = 0$ when fixing the affine basis and $c$ is a normal covector, we have $\sprod{c}{a_{n+1}}> \sprod{c}{a_{1}} = 0$, which implies that $\det(\tilde M)>0$.
	Then, we rescale some columns of~$\tilde M$ by factors of the form $\frac{1}{\sprod{c}{a_{n+1}}-\sprod{c}{a_i}}$ with $i\leq n$, which are strictly positive, again because $c$ is a normal covector.
	Moreover, we have that $\tilde a_1 = 0$, $\tilde a_2 = e_1$, $\tilde a_3 = e_2$, etc, fixing the affine basis to mod out affine transformations as desired.
	Therefore, we have that $\tilde A\in \rsinf[\chi]$.

	It remains to show that the map $(A,c)\mapsto \tilde A$ is a stable projection from~$\cC$ to~$\rsinf[\chi]$.
	Indeed, to recover a preimage $A$ from $\tilde A$ one has to undo the projective transformation $\phi$ described above.
	To do so, we need to choose a hyperplane $H$ that is not intersecting $\phi(\conv(A)) = \conv(\tilde A\setminus \tilde a_{n+1})+\cone(\tilde a_{n+1})$, make a projective transformation sending~$H$ to infinity, and then perform an affine transformation sending the affine basis to the standard coordinates.
	Note that $c$ is determined as the normal covector to the image of the hyperplane at infinity.
	The set of preimages is therefore parametrized by the set of possible hyperplanes~$H$, which is a polyhedral cone --the polar cone of~$\phi(\conv(A))$-- defined by polynomial equations.
	Therefore, we have that $\cC\sim\rsinf[\chi]$.
	Combining this with the previous stable equivalence $\cC \sim \rs[\chi]$, we obtain~$\rs[\chi] \sim \cC \sim \rsinf[\chi]$, finishing the proof.
\end{proof}

\SS \label{lem:cyclicliftchiro}
We are almost ready to state and prove our main result.
The next step is to associate in \cref{prop:cyclicomegadiagramfromom} an $f$-orientation on a standard simplex $\ssimplex[d]$ to any planar point configuration $(a_1,\dots,a_{n+1})\in \R^{2\times (n+1)}$.

Let $A = (a_1,\dots,a_{n}) \in \R^{2 \times n}$ be a planar point configuration.
The \defn{cyclic lift} of~$A$ is the point configuration~$\cyclift$ in $\R^{n-1}$ given by the columns $\cyclift[a_1],\dots, \cyclift[a_n]$ of the following matrix
\[
\begin{pmatrix}
	(a_1)_1& (a_2)_1 & \cdots &(a_n)_1\\
	(a_1)_2& (a_2)_2 & \cdots &(a_n)_2\\
	K_3^1& K_3^2 & \cdots &K_3^n\\
	K_4^1& K_4^2 & \cdots &K_4^n\\
	\vdots&\vdots &\ddots &\vdots \\
	K_{n-1}^1& K_{n-1}^2 & \cdots &K_{n-1}^n
\end{pmatrix}
\]
for some $0\ll K_3\ll K_4\ll \cdots \ll K_{n-1}$ large enough so that the following lemma holds.
Here, we denote by~$(a_i)_j$ and $K_i^j$ the $j$th coordinates of the points~$a_i$ and~$K_i$ in the canonical frame, respectively.

\lemma
Let $A$ be a planar point configuration, and let $\pi_k$ denote the canonical system of projections in $\R^{n-1}$.
If $K_{k+1}$ is large enough with respect to~$K_3,K_4,\ldots,K_{k}$ and $A$, then the chirotope of $\pi_{k+1}(\cyclift[A])$ is determined by the chirotope of $\pi_{k}(\cyclift[A])$.

\begin{proof}
	Consider subset of $k+2$ indices $1\leq i_1<i_2<\dots<i_{k+2}\leq n$.
	Then, if $K_{k+1}$ is large enough we have that
	\[\sign
	\begin{vmatrix}
		1&1&\cdots&1&1\\
		(a_{i_1})_1& (a_{i_2})_1 & \cdots &(a_{i_{k+1}})_1&(a_{i_{k+2}})_1\\
		(a_{i_1})_2& (a_{i_2})_2 & \cdots &(a_{i_{k+1}})_2&(a_{i_{k+2}})_2\\
		K_3^{i_1}& K_3^{i_2} & \cdots &K_3^{i_{k+1}}&K_3^{i_{k+2}}\\
		K_4^{i_1}& K_4^{i_2} & \cdots &K_4^{i_{k+1}}&K_4^{i_{k+2}}\\
		\vdots&\vdots &\ddots &\vdots \\
		K_{k}^{i_1}& K_{k}^{i_2} & \cdots &K_{k}^{i_{k+1}}&K_{k}^{i_{k+2}}\\
		K_{k+1}^{i_1}& K_{k+1}^{i_2} & \cdots &K_{k+1}^{i_{k+1}} &K_{k+1}^{i_{k+2}}
	\end{vmatrix}
	 =
	\sign
	\begin{vmatrix}
		1&1&\cdots&1\\
		(a_{i_1})_1& (a_{i_2})_1 & \cdots &(a_{i_{k+1}})_1\\
		(a_{i_1})_2& (a_{i_2})_2 & \cdots &(a_{i_{k+1}})_2\\
		K_3^{i_1}& K_3^{i_2} & \cdots &K_3^{i_{k+1}}\\
		K_4^{i_1}& K_4^{i_2} & \cdots &K_4^{i_{k+1}}\\
		\vdots&\vdots &\ddots &\vdots \\
		K_{k}^{i_1}& K_{k}^{i_2} & \cdots &K_{k}^{i_{k+1}}\\
	\end{vmatrix}.
	\]
	Indeed, if we develop the determinant on the left hand side by the last row, we see that for $K_{k+1}$ large enough it is dominated by the factor with coefficient $K_{k+1}^{i_{k+2}}$, whose sign is that of the determinant on the right hand side.
	By definition (\cref{def:chirotope-point}), the signs of the determinants on left hand side for all possible subsets of indices~$\{i_1,\ldots,i_{k+2}\}$ define the chirotope of~$\pi_{k+1}(\cyclift[A])$, while the ones on the right hand side define the chirotope of~$\pi_{k}(\cyclift[A])$.
	Thus, the proof is complete.
\end{proof}

\SS\label{prop:cyclicomegadiagramfromom}
\theorem
Let $A = (a_1,\dots,a_n,a_{n+1}) \in \R^{2 \times (n+1)}$ be a planar point configuration in general position.
Then there is a $\ssimplex[n-1]$-admissible frame~$B$ of $\R^n$ whose realization space is stably equivalent to the extended realization space of $A$.
That is, we have $$\msp[B][{\ssimplex[n-1]}] \sim \rsinf[A].$$

\begin{proof}
	Relabeling the vertices if necessary, we can assume that $a_{n+1}$ is a vertex of the convex hull of~$A$.
	Consider the point configuration $\tilde A = (\tilde a_1,\dots,\tilde a_n)$ obtained by a projective transformation sending $a_{n+1}$ to infinity, as defined in the proof of~\cref{lem:rsandrsinftystablyequivalent}.
	Consider the cyclic lift~$\cyclift[\tilde A]$ of $\tilde A$ (see \cref{lem:cyclicliftchiro}), and denote by $P \defeq \conv(\cyclift[\tilde A])$ its convex hull.
	Since the points of $\tilde A$ are affinely independent, we have that $P$ is a simplex in~$\R^{n-1}$.
	We define a frame $B = (v_1,\dots,v_{n-1})$ by
	$$v_i \defeq
	\begin{cases}
		\tilde a_2-\tilde a_1 & \text{ if } i = 1, \\
		\tilde a_{n+1} & \text{ if } i = 2, \\
		e_i & \text{ if } i\geq 3.
	\end{cases}$$
	We have that $B$ is $P$-admissible: only the first vector $v_1$ is not generic with respect to~$P$, but the admissibility condition for the first vector of a frame is void.
	Let us denote by~$\gamma$ the $f$-orientation associated to the framed polytope~$(P,B)$.
	We can identify it with a $f$-orientation of the standard simplex, by applying an affine transformation sending $P$ to $\ssimplex[n-1]$ (\cref{ss:invariancelinearautomorphism}).

	Our first claim is that $\gamma$ is completely determined by the chirotope of~$A$.
	In view of \cref{lem:diagramsareflagoms}, it suffices to show that the chirotopes of $\pi_k(\cyclift[\tilde A])$ are completely determined by the chirotope of~$A$.
	Note first that, as argued in the proof of \cref{lem:rsandrsinftystablyequivalent}, the chirotope of $A$, which by definition is the chirotope of the vector configuration
	$\left( \binom{1}{a_1},\dots, \binom{1}{a_n}, \binom{1}{a_{n+1}}\right)$, coincides with the chirotope of the vector configuration $\left( \binom{1}{\tilde a_1},\dots, \binom{1}{\tilde a_n}, \binom{0}{\tilde a_{n+1}}\right)$.
	Now, by construction we have that~$\pi_2(\cyclift[\tilde A]) = \tilde A$, and thus the chirotope of $\cyclift[\tilde A]$ is just a subset of the chirotope of~$A$.
	Moreover, in virtue of \cref{lem:cyclicliftchiro}, the chirotopes of $\pi_k(\cyclift[\tilde A])$ with $k\geq 3$ are determined by the chirotope of~$\tilde A$.
	Finally, by reading off the signs of the determinants involving~$a_{n+1}$ in~$A$, we can read off the chirotope of $\pi_1(\cyclift[\tilde A])$, as done in \cref{lem:diagramsareflagoms}.
	There could be a global sign change, because the chirotope can change its sign if the basis is changed; however, by setting $v_1 = \tilde a_2-\tilde a_1$, we know that in the basis $B$ the segment $[\tilde a_1,\tilde a_2]$ has a positive orientation, which does fix the global sign.
	This finishes the proof that $\gamma$ is completely determined by the chirotope of~$A$.

	Now, using \cref{lem:diagramsareflagoms} again, we have that for any frame $B'\in \msp[B][{\ssimplex[n-1]}]$, the $2$-dimensional projection gives a point configuration $A' = (a_1',\dots,a_n')$ with the same chirotope as $A$.
	This, together with the vector~$v_2$, gives a realization in $\rsinf[A]$.
	Thus, we have a natural projection~$\msp[B][{\ssimplex[n-1]}] \to\rsinf[A]$.
	To prove our claim it then suffices to show that this projection is stable.
	By \cref{lem:cyclicliftchiro}, one can always find a cyclic lift of a given point configuration in $\rsinf[A]$.
	Moreover, by the preceding paragraph, the $f$-orientation associated to this lift will always be the same, and the associated frame $B'$ will be in $\msp[B][{\ssimplex[n-1]}]$.
	Therefore, the fibers of the projection~$\msp[B][{\ssimplex[n-1]}] \to\rsinf[A]$ are non-empty.

	To see that this projection is stable, we split it in a sequence of projections, each deleting one coordinate.
	The conditions are that at each step we have the chirotope of $\pi_k(P)$.
	The conditions on the $k$th coordinate, once the smaller coordinates are known, are given by determinants, which are linear in this coordinate, and thus the fibers are open polyhedra whose inequalities are given by polynomials on the fixed coordinates.
	All the inequalities are strict and, since the fibers are not empty, they are full-dimensional at each step.
	In particular, they are of constant dimension, which addresses the problem in~\cite{Verkama2023}.
	Therefore, the projection~$\msp[B][{\ssimplex[n-1]}] \to\rsinf[A]$ is stable, and we have $\msp[B][{\ssimplex[n-1]}] \sim \rsinf[A]$, as desired.
\end{proof}

\SS\label{thm:simplex-universality}
\theorem
For every open primary basic semi-algebraic set $S$ defined over $\Z$ there is a $f$-orientation on some standard simplex $\gsimplex_{d}$ whose realization space is stably equivalent to~$S$.

\begin{proof}
	By the universality theorem for oriented matroids \cref{thm:universality}, there is a planar point configuration~$A$ stably equivalent to $S$, that is, with $\rs[A]\sim S$.
	By \cref{prop:cyclicomegadiagramfromom}, there is a $f$-orientation~$\gamma$ on some standard simplex $\ssimplex[d]$ such that $\ms[\gamma]\sim \rsinf[A]$, and by \cref{lem:rsandrsinftystablyequivalent}, we also have $\rsinf[A]\sim \rs[A]$.
	Therefore, for any open primary basic semi-algebraic set $S$, there is a $f$-orientation $\gamma$ on some $\ssimplex[d]$ such that $\ms[\gamma]\sim \rsinf[A] \sim \rs[A] \sim S$, finishing the proof.
\end{proof}

\SS\theorem
For every open primary basic semi-algebraic set $S$ defined over $\Z$ there is
a full flag chirotope whose realization space is stably equivalent to~$S$.

\begin{proof}
	This is just a reformulation of \cref{thm:simplex-universality} in the language of \cref{def:realizable-flag}.
\end{proof}

%% file: sec/shelling.tex

\section{Framed polytopes and regular directed complexes}\label{part:molecules}

Oriented graded posets provide another model for pasting diagrams \cite{Hadzihasanovic24, Hadzihasanovic20, Steiner93}.
In this section we show that loop-free framed polytopes also define regular directed complexes.
The main difference compared to the model used in previous sections is that this requires the existence of layerings for sources and targets, which are explicit orders in which the cells are composed in the associated category.
The proof requires extending the notion of source and target to coherent subdivisions, which is an interesting topic by itself, and the layerings obtained are closely related to the line shellings of facets induced by the frame.
An interesting consequence is that the $k$-source and $k$-target of a framed polytope are always $(k-1)$-loop-free.

\subsection{Oriented graded posets}

An \defn{oriented graded poset} is a graded poset together with a labelling of the edges of its Hasse diagram by an element of~$\set{-,+}$.
These posets play a fundamental role in the theory of \emph{diagrammatic sets}, a model of higher categories developed in \cite{Hadzihasanovic24,Hadzihasanovic20,Steiner93} named after the work of Kapranov--Voevodsky \cite{KapranovVoevodskyGroupoids91}.

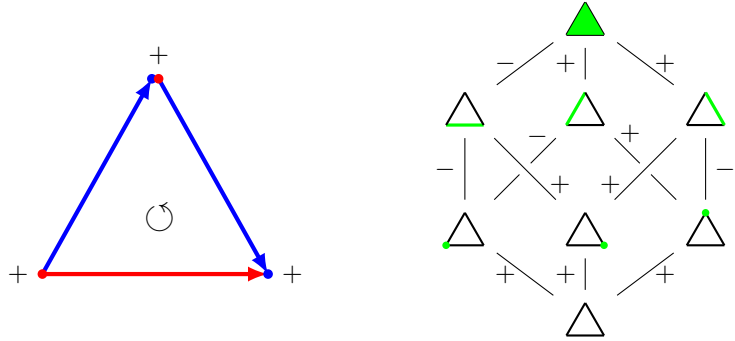
\begin{figure}
	\input{fig/orientations}
	\caption{The oriented graded poset associated to a $2$-simplex.}
	\label{fig:orientations-poly}
\end{figure}

\SS Let $P$ be a polytope.
A \defn{poset orientation} on $\faces$ is the assignment of an element in $\set{-,+}$ to each covering relation in~$\faces$ (\emph{i.e.} a pair $F \lessdot G$ in $\faces$ where $F$ is a facet of $G$).
If $\beta$ is an $f$-orientation, then the poset orientation on $\faces$ \defn{induced by $\beta$} is the one to the covering relation $F \lessdot G$ in $\faces$ the sign
\begin{equation}\label{eq:facial2facelattice}
	\beta_{F \lessdot G}
	\defeq
	\begin{cases}
		- &\text{ if }\beta_F \wedge v \sim -\beta_G,\\
		+ &\text{ if }\beta_F \wedge v \sim \beta_G \text{ or }F=\emptyset;
	\end{cases}
\end{equation}
where $v\in \Lin_G$ is an outer-pointing vector (\emph{i.e.} $\sprod{\normal_F^G}{v}>0$). The definition does not depend on the choice of $v$ by the same argument as in \cref{d:sources_and_targets}.
See \cref{fig:orientations-poly} for an example.

\SS\label{d:sources_and_targets}

The induced $f$-orientation of a framed polytope is completely determined by its induced poset orientation (and vice versa).
Explicitly, the sources and targets (which determine the $f$-orientation by \cref{ss:sotagiveforientation}) are
\[
\so(F) = \set{E \mid \beta_{E \lessdot F} = -}, \qquad\text{ and }\qquad
\ta(F) = \set{E \mid \beta_{E \lessdot F} = +}.
\]

\SS\label{l:face_vs_poset_orientation}\label{p:oriented-thin}
Not all poset orientations on $\faces$ are induced from an $f$-orientation on $P$.
We will now present necessary and sufficient conditions for a poset orientation on $\faces$ to be induced this way.

\begin{figure}
	\input{fig/thin}
	\caption{A thin poset orientation (left) and a poset orientation that is not thin (right).}
	\label{f:thin}
\end{figure}
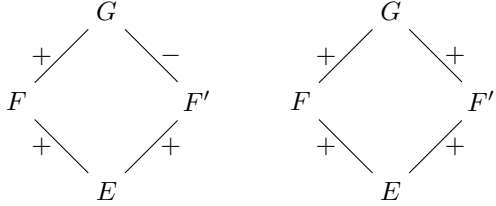

We start by remarking that the poset $\faces$ satisfies the \defn{diamond property} \cite[Thm.~2.7]{Ziegler95}, that is to say, for every pair of faces $E < G$ of codimension 2, there are exactly two faces $F$ and $F'$ between them.
A poset orientation $\beta$ on $\faces$ is said to be \defn{thin} if for any such diamond we have
\[
\beta_{E \lessdot F} \cdot \beta_{F \lessdot G} =
-\beta_{E \lessdot F'} \cdot \beta_{F' \lessdot G}.
\]
Please consult \cref{f:thin} for examples.
We say that a poset orientation is \defn{grounded} if $\beta_{\emptyset \lessdot v}=+$ for any vertex $v$ of $P$.

\medskip\theorem
The induced poset orientation construction defines a bijection between $f$-orientations on $P$ and grounded thin poset orientations on $\faces$.

\begin{proof}
	An induced poset orientation is grounded by definition.
	Let us show that it is thin.
	It suffices to take $P$ to be the top vertex of our diamond with the codimension 2 face $F$ as its bottom, and the facets $G$ and $G'$ in between.
	If $F = \emptyset$, then $G$ and $G'$ are the two vertices of a segment $P$, and thinness follows from the fact that if $v$ is outer-pointing for~$G$ then it is inner-pointing for $G'$.
	If $F\neq\emptyset$, let $v\in \Lin_G$ and $v'\in \Lin_{G'}$ be outer-pointing vectors for $F$ with respect to $G$ and $G'$, respectively.
	Therefore, $v$ and $v'$ are respectively outer-pointing for $G'$ and $G$ with respect to $P$.
	We have thus that
	\begin{align*}
		\beta_G &=	\beta_{F \lessdot G} \ \beta_F \wedge v,\\
		\beta_P &=	\beta_{G \lessdot P} \ \beta_G \wedge v'= \beta_{G \lessdot P}\ \beta_{F \lessdot G} \ \beta_F \wedge v \wedge v',
	\end{align*}
	and that
	\begin{align*}
	 	\beta_{G'} &= \beta_{F \lessdot G'} \ \beta_F \wedge v',\\
	 	\beta_P		&=	\beta_{G' \lessdot P}\ \beta_{G'} \wedge v= \beta_{G' \lessdot P}\ \beta_{F \lessdot G'} \ \beta_F \wedge v' \wedge v.
	\end{align*}
	By the antisymmetry of the wedge product, we conclude that \[\beta_{F \lessdot G}\cdot\beta_{G \lessdot P} =- \beta_{F \lessdot G'}\cdot\beta_{G' \lessdot P}\] as desired.

	For the converse, we construct now an $f$-orientation on $P$ inducing a given grounded thin poset orientation on~$\faces$.
	We start by orienting all the vertices with $1\in \R$ and recursively define the orientation of higher-dimensional faces.
	For $e\geq 1$, assume that all faces of dimension $e-1$ are oriented, and let $F$ be an $e$-face.
	Then we can define the orientation of~$F$ by applying \eqref{eq:facial2facelattice} with respect to any facet~$G$ of~$F$.
	It remains to see that this is well defined; that is, that we get the same orientation for~$F$ regardless of the facet~$G$ we choose.
	If $G$ and $G'$ are neighboring facets of~$F$, meaning that their intersection is an $(e-2)$-face of $P$, then thinness implies that both induce the same orientation on~$F$, just by reversing the argument above.
	We conclude by noting that the dual graph of~$P$ (which has a vertex for each facet and an edge when they share a codimension~$2$ face) is connected (it is the graph of the polar polytope~\cite[Cor.~2.14]{Ziegler95}, which is connected in dimensions~$\geq 1$~\cite[Thm.~3.14]{Ziegler95}).
	Therefore for any pair of facets $G$ and $G'$ of $F$, we can find a path $G_1=G,G_2,\dots,G_k=G'$ of facets of $F$ in which $G_i$ and $G_{i+1}$ are neighbors.
	By thinness all the $G_i$'s induce the same orientation on~$F$.
\end{proof}

In light of this theorem, on can regard $f$-oriented polytopes as either polytopes with oriented non-empty faces and positively oriented vertices, or as polytopes whose poset of non-empty faces is equipped with a grounded thin poset orientation.

\subsection{Regular directed complexes}

Our goal is now to define \emph{regular directed complexes}, which are oriented graded posets with specific properties \cite{Hadzihasanovic24,Hadzihasanovic20,Steiner93}.
For concreteness, we will give the definitions in the case where $P$ is a polytope with a poset orientation.

\SS The \defn{closure} of a subset $U \subset \faces$ is the set
\[
\overline{U} \defeq \set{F \in \faces \mid F \leq G \text{ and } G \in U}.
\]
We say that $U$ is \defn{closed} if $\overline{U} = U$.

\SS \label{def:k-source-target-bis}

For~$k\geq 0$, the \defn{$k$-source}~$\so_k(U)$ (resp.\ \defn{$k$-target}~$\ta_k(U)$) of a closed subset~$U \subset \faces$ is the set of $k$-faces~$F$ in~$U$ such that every covering relation~$F \lessdot G$ with~$G$ in~$U$ has label~$-$ (resp.\ $+$).

\medskip\remark
Note that in contrast to the preceding definition \cref{d:sources_and_targets}, the $k$-source~$\so_k(\overline{F})$ or target~$\ta_k(\overline{F})$ of a $k$-face $F$ is defined, and equal to~$F$ itself.

\SS The set of \defn{maximal elements} of a closed subset $U \subset \faces$ is
\[
\max(U) \defeq \set{F \in U \mid F \leq G \in U \text{ implies } F = G}.
\]
We denote by $\max_{< k}(U)$ the subset of $\max(U)$ containing faces of dimension~$<k$.

\SS\label{def:boundaries}
The \defn{$k$-input} and \defn{$k$-output boundaries} of a closed subset $U \subseteq \faces$ are defined by
\[
\bso_k(U) \defeq \overline{\so_k(U)} \cup \overline{\max\nolimits_{< k}(U)}
\quad \text{and} \quad
\bta_k(U) \defeq \overline{\ta_k(U)} \cup \overline{\max\nolimits_{< k}(U)}.
\]
We denote by $\bso_k(F)\defeq \bso_k(\overline{F})$ (resp.\ $\bta_k(F) \defeq \bta_k(\overline{F})$) the $k$-input (resp.\ output) boundary of a face~$F \in \faces$.
Note that in this case, we have $\bso_k(F) = \overline{\so_k(F)}$ (resp.\ $\bta_k(F) = \overline{\ta_k(F)}$).

\SS The \defn{dimension} of a subset $U \subset \faces$ is the maximum dimension of an element in $U$.
If $U$ is of dimension $n$, we abbreviate $\bso_{n-1}(U)$ and $\bta_{n-1}(U)$ by $\bso(U)$ and $\bta(U)$.

\SS \label{def:atoms-molecules}

An \defn{atom} consists of the closure $\overline{F}$ of a face $F \in \faces$.
For $U$ a closed subset of $\faces$, we say that $U$ is a \defn{molecule} if it is an atom, or if there exist molecules~$U_1,U_2 \subsetneq U$ and $k \in \N$ such that $U = U_1 \cup U_2$ and $U_1 \cap U_2 = \bta_k(U_1) = \bso_k(U_2)$.
In this situation, we write $U= U_1 \circ_k U_2$.

\SS\label{def:layering}
Let $U$ be a molecule, let $-1 \leq k < \dim U$, and define
$m \defeq \left| \cup_{i>k} (\max(U)_i)\right|.$
A \defn{$k$-layering} of $U$ is an ordered sequence of molecules $(U_1,\ldots,U_m)$ such that $$U=U_1 \circ_k \cdots \circ_k U_m$$
and $\dim U_i > k$ for all $1 \leq i \leq m$.

\SS\label{def:rdc}
A polytope $P$ with a poset orientation is a \defn{regular directed complex} if for all $d$-dimensional faces~$F \in \faces$, we have that
\begin{enumerate}
	\item \label{it:globularcond1} $\bso(F)$ and $\bta(F)$ are molecules,
	\item \label{it:globularcond2} $\bso(\bta(F)) = \bso(\bso(F))$ and $\bta(\bso(F)) = \bta(\bta(F))$ if $d>1$,
	\item \label{it:globularcond3} $\bso_k(F) \cap \bta_k(F) = \bso_{k-1}(F) \cup \bta_{k-1}(F)$ for all $k < d$.
\end{enumerate}
The identities~\eqref{it:globularcond2} are called the \defn{globular identities}.

\SS\remark
We use here the definition of regular directed complex from~\cite[Def.~1.29]{Hadzihasanovic20}.
Little work shows that it is equivalent to~\cite[Def.~5.3.1]{Hadzihasanovic24}.


\subsection{Generalized sources and targets}

Our goal is to prove in \cref{thm:gobular-structure} that framed polytopes are regular directed complexes.
Conditions \eqref{it:globularcond2} and \eqref{it:globularcond3} above are easy to verify.
To prove that sources and targets are molecules, we will slightly enlarge our setup: we generalize the definition of input and output boundaries, and prove the stronger statement that these satisfy \eqref{it:globularcond1}.

\SS\label{def:new-sources-and-targets}

Let~$P$ be a polytope in $\R^d$ and let~$B$ be a frame of~$\R^d$ which is \emph{not necessarily $P$-admissible}.
Let $F\in\faces$ be a face of dimension $\dim F>k$ such that $\dim(\pi_{k+1}(F)) = k+1$.
The \defn{$k$-source}~$\so_k(F)$ (resp. \defn{$k$-target}~$\ta_k(F)$) of~$F$ is the set of faces~$\pi_{k+1}^{-1}(G)\cap F$ where~$G$ is a facet of~$\pi_{k+1}(F)$ such that
\[
\sprod{\normal_{G}^{\pi_{k+1} F}}{v_{k+1}} < 0 \quad (\text{resp.} > 0).
\]
For $P$-admissible frames, this coincides with the definitions given in~\cref{lem:sources and targets in a framed polytope}.
Note that if the frame is not admissible, then the $k$-source and $k$-target can contain faces of dimension larger than~$k$.

\SS For a closed subset~$U \subset \faces$, the \defn{$k$-source} $\so_k(U)$ (resp.\ \defn{$k$-target} $\ta_k(U)$)
is the subset of faces $F\in U$ such that $\dim(\pi_{k+1}(F)) = k$ but that do not belong to any $\ta_k(F)$ (resp. $\so_k(F)$) for an $F\in U$ such that $\dim(\pi_{k+1}(F)) = k+1$.
We define $\max_{< k}(U)$ as the set of faces of $\max(U)$ such that the projection $\pi_{k+1}$ is of dimension smaller than~$k$.
We then define \defn{$k$-input} and \defn{$k$-output boundaries}, \defn{atoms} and \defn{molecules} exactly as in~\cref{def:boundaries} and~\cref{def:atoms-molecules}.
Again, note that for $P$-admissible frames, this coincides with the previous definitions.

\SS\remark\label{rem:no-longer-acyclic}
The above definition of $k$-source and target is just a reformulation of the concept of \emph{$\pi_k$-induced coherent subdivision} of $\pi_{k}(F)$ as introduced in~\cite{BilleraSturmfels1992}, see also \cite[Def.~9.2]{Ziegler95}.
Comparing with \cref{s:subdivision}, the $k$-source and target are still coherent subdivisions of $\pi_k(F)$, but they are no longer necessarily tight.
Using \cref{thm:globular-flag-chirotopes}, under this new definition we will have that $f$-orientations of framed simplices are in bijection with acyclic realizable full flag chirotopes, which are no longer necessarily uniform.


\subsection{Layerings via line shellings}

We will use a construction inspired by Bruggesser and Mani's line shellings \cite{BruggesserMani1971} to construct layerings of sources and targets.

\SS\label{def:Pk-admissible}
Let $P \subset \R^d$ be a polytope.
We say that a frame $(v_1,\dots,v_d)$ is \defn{$P_k$-admissible} if $(v_1,\dots,v_k)$ is {$\pi_k(P)$-admissible}.
That is, all possible degeneracies arise from the first projections $\pi_d,\dots\pi_k$, but the vectors $v_1,\dots,v_k$ can be perturbed without altering the orientation structure of~$P$.

\SS\label{def:stshelling}
Let $P$ be a $d$-dimensional polytope in $\R^d$ and $B = (v_1,\ldots,v_d)$ be a $P_k$-admissible frame for some $1\leq k<d$.
If needed, replace $v_k$ with a generic pertubation that will be detailed later, and define the vector $${w_{k+1}(\lambda)} \defeq \lambda v_{k+1}-v_{k}$$ depending on a parameter~$\lambda \in \R$.
For each facet~$F_i$ of $\pi_{k+1}(P)$, let $\lambda_i$ be the value for which $\sprod{\normal_{F_i}}{w_{k+1}(\lambda_i)} = 0$.
That is, we define
\[
\lambda_i \defeq \frac{\sprod{\normal_{F_i}}{v_{k}}}{\sprod{\normal_{F_i}}{v_{k+1}}}.
\]
The genericity condition that we impose is that all the $\lambda_i$'s are distinct. To this end, if needed we consider a perturbation of $v_{k}$ that is close enough so that it does not change the order of any two different $\lambda_i$'s, and generic enough so that they become all distinct.
Note that since the frame is $P_k$-admissible, this does not alter the orientation structure of~$P$.

\SS We use these numbers to define total orders $<_{\so}$ and $<_{\ta}$ on the faces in $\so_k(P)$ and $\ta_k(P)$.
Each element $G_i\in\so_k(P)$ or $G_i\in\ta_k(P)$ is of the form~$\pi_{k+1}^{-1}(F_i)\cap P$ for a facet~$F_i$ of~$\pi_{k+1}(P)$.
We set $G_i<_{\so} G_j$ if $\lambda_i>\lambda_j$ for $G_i,G_j\in \so_k(P)$, and $G_i<_{\ta} G_j$ if $\lambda_i<\lambda_j$ for $G_i,G_j\in \ta_k(P)$.

\SS\label{prop:orderst} \proposition
For any two faces $F_1,F_2$ in the $k$-source (resp.\ $k$-target) of a polytope $P$ with respect to a $P_k$-admissible frame,
if~$\ta_{k-1}(F_1)\cap \so_{k-1}(F_2)\neq \emptyset$ then we have $F_1 <_{\so} F_2$ (resp. $F_1 <_{\ta} F_2$).

\begin{proof}
We do the proof for the target, and the proof for the source is analogous.
Note that if a face $F_i$ belongs to the $k$-target of~$P$ then $G_i \defeq \pi_{k+1}(F_i)$ is a facet of $\pi_{k+1}(P)$ whose normal covector evaluates positively on $v_{k+1}$ (\cref{def:new-sources-and-targets}). Rescaling, we can thus choose a normal covector $\normal_i$ such that $\sprod{\normal_i}{v_{k+1}} = 1$.
Now let $G \defeq G_1\cap G_2$ be the $(k-1)$-face given by the intersection of the two facets. We have that $G=\pi_{k+1}(F)$, where $F=F_1\cap F_2=\ta_{k-1}(F_1)\cap \so_{k-1}(F_2)$.
We have that $\normal_2 - \normal_1$ is a normal covector of $\pi_{k}(G)$ with respect to $\pi_{k}(F_1)$ in~$V_{k}$. Indeed, by definition $\sprod{\normal_i}{x}\geq \sprod{\normal_i}{y}$ for any $x\in G_i$ and $y\in \pi_{k+1}(P)$; thus $\sprod{\normal_2-\normal_1}{x}\geq \sprod{\normal_2-\normal_1}{y}$ for any $x\in G=G_1\cap G_2$ and $y\in G_1$. Moreover, we have that $\sprod{\normal_2-\normal_1}{v_{k+1}} = 0$ and hence we can interpret it as a covector in $V_k^\ast$ that is a normal covector to~$\pi_{k}(G)$.
Now, since $F\in\ta_{k-1}(F_1)\cap \so_{k-1}(F_2)$, this means that $\sprod{\normal_2-\normal_1}{v_{k}}>0$.
Therefore, we have $\sprod{\normal_1}{v_{k}}<\sprod{\normal_2}{v_{k}}$.
Since $\sprod{\normal_i}{v_{k+1}} = 1$, this means that $\lambda_1<\lambda_2$ and that we have~$F_1 <_{\ta} F_2$.
\end{proof}

\SS\corollary\label{rem:loop-freeness-source-target}
This gives a proof of $(k-1)$-loop-freeness for the $k$-source and $k$-target of a framed polytope.
Indeed, since $<_{\so}$ and $<_{\ta}$ are total orders, there are no cellular $(k-1)$-strings of infinite length.
This generalizes \cref{prop:dimP-2}.

\SS\remark
The orderings $<_{\so}$ and $<_{\ta}$ do not necessarily give shellings of $\so(P)$ and $\ta(P)$, but give an ordering that is strongly related, tailored to the construction of layerings that we use in the sequel (see~\cref{def:layering}).
Their existence shows that sources and targets of framed polytopes are ``stackable'' in the sense of \cite[Sec.~3]{Athanasiadis01}.

More precisely, the polytopal complex $\bd F_j\cap (\cup_{i=1}^{j-1}\bd F_i)$ might be not pure, or not connected, and thus not shellable (see for example \cite[Def.~8.1]{Ziegler95} for the definition of shelling).
However, we could solve this issue if we considered some stratified version of shellability in which not only the maximal cells were considered.
The key property \cref{prop:orderst} shows that for any $F_j\in\so(P)$, we have that $\bd F_j\cap (\cup_{i=1}^{j-1}\bd F_i)$ is the subset of $\so(F_j)$ containing the faces of the form $\ta(F_i)\cap \so(F_j)$ for some $F_i\in \so(P)$.
However, this does not cover all of $\so(F_j)$, as some faces in $\so(F_j)$ belong to $\so(\so(P))$ and are not in any target.
If in our order we had not only the facets of $P$, but we also allowed faces of lower dimensions, then we could add all the cells of $\so(\so(P))$ before starting adding the facets in $\so(P)$, then we would have $\bd F_j\cap (\cup_{i=1}^{j-1}\bd F_i)=\so(F_j)$.
This is the key property that allows us to construct a layering in \cref{cor:st-molecules}.
For reasons of scope, we omit more details and a formal definition of such a stratified shellability.


\subsection{Framed polytopes are regular directed complexes}

We now combine the results of the preceding two sections in order to prove our main result \cref{thm:gobular-structure}.
The most important part of the proof is to show that sources and targets of atoms form molecules.
This is accomplished by the construction of layerings of $\bso_k(F)$ and $\bta_k(F)$ induced by the frame.

\SS\theorem\label{cor:st-molecules}
Let $P\subset \R^d$ be a polytope, let $F$ be a face of~$P$, and let $k \leq \dim F$.
Let $B = (v_1,\ldots,v_d)$ be a $P_k$-admissible frame.
Then, $\bso_k(F)$ and $\bta_k(F)$ are molecules.

\begin{proof}
	We lose no generality by assuming $F = P$.
	Moreover, it suffices to show the statement for $\bso_k(P)$, the case of $\bta_k(P)$ is analogous.
	We will assume that the $\lambda_i$'s from~\cref{def:stshelling} are all distinct, as otherwise we can replace $v_k$ by a generic perturbation~$v_k'$, as in~\cref{def:stshelling}, without altering the orientation structure.
	We proceed by induction on~$k \geq 0$.
	For~$k = 0$, we have that $\bso_0(P)$ is always a face of $P$, thus an atom and a molecule.
	Let $P$ be a polytope of dimension $\dim P \geq k+1$.
	We suppose that $\bso_{k-1}(P)$ is a molecule (for any frame) and show that~$\bso_{k}(P)$ (for this frame) is a molecule.

	We consider the total order $G_1,\ldots,G_\ell$ on the maximal faces of~$\bso_{k}(P)$ from \cref{def:stshelling}.
	We now show by induction on $j \geq 0$ that $\cG_j \defeq \bso_{k-1}(P) \cup \bigcup_{i = 1}^{j} \overline{G_i}$ is a molecule.
	For $j = \ell$, this will prove that $\bso_{k}(P)$ is a molecule, completing the induction step.
	For $j = 0$, there is nothing to show, since $\cG_0 =\bso_{k-1}(P)$, which is a molecule by the induction hypothesis on $k$.
	Suppose now that $\cG_m$ is a molecule in order to show that~$\cG_{m+1}$ is a molecule.
	For this, it suffices to show that $\cH_{m+1} \defeq \overline{G_{m+1}} \cup \bta_{k-1}(\cG_m)$ is a molecule.
	Indeed, by \cref{prop:orderst}, we have the inclusion~$\bso_{k-1}(G_{m+1}) \subset \bta_{k-1}(\cG_m)$, and thus~$\bta_{k-1}(\cG_m) = \bso_{k-1}(G_{m+1} \cup \bta_{k-1}(\cG_m))$.
	Assuming that $\overline{G_{m+1}} \cup \bta_{k-1}(\cG_m)$ is a molecule allows one to write $\cG_{m+1} = \cG_m \cup \overline{G_{m+1}}$ which is a gluing of two molecules ($\cG_m$ and $\cH_{m+1}$) satisfying the recursive definition \cref{def:atoms-molecules}.
	Therefore, we are left with showing that~$\cH_{m+1}$ is a molecule.

	To prove that~$\cH_{m+1}=\overline{G_{m+1}} \cup \bta_{k-1}(\cG_m)$ is a molecule, we will show that it is precisely the $(k-1)$-output boundary of $P$ for the frame
	\[\tilde B(\lambda_{m+1}) = (v_1,\dots,v_{k-1},\tilde v_k \defeq v_{k+1},\tilde v_{k+1} \defeq w_{k+1}(\lambda_{m+1}),v_{k+2},\ldots),\]
	and conclude by induction.
	To avoid confusion, we denote by $\tilde \pi_1,\ldots,\tilde \pi_n $ the system of projections associated to~$\tilde B$.
	Our induction will prove that $\overline{G_{m+1}} \cup \bta_{k-1}(\cG_m)$ is a molecule for~$\tilde B$.
	However, note that $\pi_j = \tilde \pi_j$ for all~$j$ except for~$j = k$ and~$j = k+1$.
	In particular, we have $\pi_{k-1} = \tilde \pi_{k-1}$.
	As the first~$k-1$ vectors of both frames coincide, all $j$-sources and targets with $j < k-1$ coincide for both frames.
	Therefore, proving that $\overline{G_{m+1}} \cup \bta_{k-1}(\cG_m)$ is a molecule for~$\tilde B$ implies that is is a molecule for~$B$, which is the result that we need.

	We start by describing recursively $\bta_{k-1}(\cG_m)$.
	Its maximal faces are, by definition, those that are in the $(k-1)$-target of some face in $\cG_m$, but not in the $(k-1)$-source of any other face in~$\cG_m$.
	For~$m = 0$, we have $\cG_0 = \bso_{k-1}(P)$ and $\bta_{k-1}(\cG_0) = \cG_0 = \bso_{k-1}(P)$.
	For $m\geq 0$, we have that $\bso_{k-1}(G_{m+1})$ is included in $\bta_{k-1}(\cG_{m})$.
	Indeed, by~\cref{prop:orderst} any maximal face in $\bso_{k-1}(G_{m+1})$ belongs to either~$\bso_{k-1}(P)$ or to~$\bta_{k-1}(G_{i})$, for some face~$G_i$ with~$i<m$, which belongs to $\cG_m$ by definition.
	Hence, these maximal faces of~$\bso_{k-1}(G_{m+1})$ do not belong to~$\bta_{k-1}(\cG_{m+1})$.
	In contrast, the maximal faces of~$\bta_{k-1}(G_{m+1})$ are in~$\bta_{k-1}(\cG_{m+1})$: by \cref{prop:orderst} again, they do not belong to any~$G_i$ with~$i<m$.
	Thus, at the level of maximal faces, the output boundary~$\bta_{k-1}(\cG_{m+1})$ is the symmetric difference of $\bta_{k-1}(\cG_{m})$ with the $(k-1)$-boundary of~$G_{m+1}$, replacing its source by its target.

	Let us see that the same recursive description holds for the $(k-1)$-output boundary of~$P$ with respect to the frame $\tilde B(\lambda)$ as $\lambda$ increases.
	When $\lambda$ has a very large negative value, the direction of $w_{k+1}(\lambda)$ is very close to~$-v_{k+1}$, and thus $\tilde \pi_k(P)$ is very close to $\pi_k(P)$.
	With a triangular transformation as in \cref{ss:lowertriangular}, we can add $\frac{1}{|{\lambda}|}\tilde v_{k+1} = - v_{k+1}-\frac{1}{|{\lambda}|}v_{k}$ to $\tilde v_k = v_{k+1}$ to obtain $- \frac{1}{|{\lambda}|}v_{k}$.
	So, we see that the $(k-1)$-target (or equivalently its closure, the $(k-1)$-output boundary) of~$P$ with respect to~$\tilde B$ coincides with the $(k-1)$-source (or equivalently its closure, the $(k-1)$-input boundary) of $P$ with respect to~$B$.
	Similarly, we can see that if $\lambda$ has very large positive value, then the $(k-1)$-target for $\tilde B$ coincides with the $(k-1)$-target for~$B$.
	The only possible changes in the $(k-1)$-target of $\tilde B(\lambda)$ appear for~$\lambda = \lambda_i$.
	For~$\lambda = \lambda_i-\varepsilon$, with $\varepsilon$ a small positive number, the source of~$G_i$ is in the $(k-1)$-target of $\tilde B(\lambda)$.
	For~$\lambda = \lambda_i$ the whole face~$G_i$ enters the $(k-1)$-target of~$\tilde B(\lambda)$, while for $\lambda = \lambda_i+\varepsilon$, only its target remains in the $(k-1)$-target of $\tilde B(\lambda)$.

	Therefore, we have that~$\overline{G_{m+1}} \cup \bta_{k-1}(\cG_m)$ is the $(k-1)$-output boundary of $P$ for the frame $\tilde B(\lambda_{m+1})$.
	By the induction hypothesis on $k$, we have that~$\overline{G_{m+1}} \cup \bta_{k-1}(\cG_m)$ is a molecule, which finishes the proof.
\end{proof}

\SS\remark
It is clear that several different layerings of sources and targets exist.
In the case of cyclic cubes and simplices, they correspond to linear extensions of the higher Bruhat and Stasheff--Tamari orders, respectively, see \cref{sec:HBO,sec:HST}.

\SS \label{thm:gobular-structure}\theorem
Framed polytopes are regular directed complexes.

\begin{proof}
	Let $(P,B)$ be a framed polytope with frame~$B = (v_1,v_2,\ldots)$, and let $F$ be a face of dimension $\dim F = d$.
	We need to check the $3$ conditions from \cref{def:rdc}.
	Condition \eqref{it:globularcond1} is a special case of~\cref{cor:st-molecules}.
	It suffices to prove Condition~\eqref{it:globularcond2} for the $k$-input boundary, the output case is similar.
	From \cref{s:sk_tk} and \cref{s:subdivision} we get directly $$\bso(\bso(F)) = \bso(\pi_{d-1}(\bso(F))) = \bso(\pi_{d-1}(F)) = \bso(\pi_{d-1}(\bta(F))) = \bso(\bta(F)),$$ as desired.
	Alternatively, one can use \cref{p:oriented-thin} for a combinatorial, instead of geometric proof.
	It remains to prove Condition~\eqref{it:globularcond3}.
	From \cref{s:sk_tk}, we know that $\bso_k(F)\cup \bta_k(F) = \overline{\bd_k(F)}$.
	We thus need to show that $\bso_{k+1}(F)\cap \bta_{k+1}(F) = \overline{\bd_k(F)}$.
	Since $\bso_k(F) = \bso(\pi_{k+1}(F))$ for any $k$, it suffices to show that $\bso(F) \cap \bta(F) = \overline{\bd_{d-2}(F)}$, and the result follows by induction.
	Returning to the definition \cref{def:extended-st}, the maximal faces of $\overline{\bd_{d-2}(F)}$ are the facets of $\pi_{d-1}(F)$.
	Since the frame~$B$ is~$P$-admissible, these are precisely the $(d-2)$-faces of $F$ which are the intersection $G \cap G'$ of two facets~$G,G'$ of~$F$ with $\sprod{\normal_G}{v_{d}}\sprod{\normal_{G'}}{v_{d}}<0$.
	According to~\cref{lem:sources and targets in a framed polytope}, these are exactly the maximal faces of the intersection~$\bso(F) \cap \bta(F)$.
	Thus, Condition~\eqref{it:globularcond3} holds and the proof is complete.
\end{proof}

\SS\remark\label{ss:iterative-source-target}
The same properties which allowed us to show the globular identities in the proof of~\cref{thm:gobular-structure} show that input and output boundaries can be computed iteratively.
Using \cref{s:sk_tk} and \cref{s:subdivision} we compute that
$$\bso_k(F) = \bso(\pi_{k+1}(F)) = \bso(\pi_{k+1}(\bso(F))) = \bso(\bso(\pi_{k+2}(F))) = \cdots = \bso(\bso(\cdots \bso(F))).$$
Therefore, we recover the fact that the $f$-orientation of a framed polytope is completely determined by its source and target sets.

\SS\remark\label{rem:open-questions}
Any oriented graded poset $P$ has an associated $\omega$-category $\Mol_{/P}$ of molecules under $P$ \cite[Thm.~5.2.12]{Hadzihasanovic24}.
If $P$ is a regular directed complex, then $\Mol_{/P}$ is generated by the atoms of~$P$ \cite[Cor.~5.3.10]{Hadzihasanovic24}.
There are moreover four notions of acyclicity for $\Mol_{/P}$, related by strict implications as follows
\[ \text{acyclic} \implies \text{strongly dim-wise acyclic} \implies \text{dim-wise acyclic} \implies \text{frame-acyclic},\]
see \cite{HadzihasanovicKessler2024} or \cite[Ch.~8]{Hadzihasanovic24}.
If $P$ is frame-acyclic, then the $\omega$-category $\Mol_{/P}$ is free on the atoms of $P$, that is, $\Mol_{/P}$ is a polygraph \cite[Thm.~8.2.14]{Hadzihasanovic24}.
It would be interesting to know if the strict chain of implications above holds for framed polytopes.
In particular, it seems desirable to determine if all framed polytopes are frame-acyclic, and thus always define polygraphs.

\SS\remark\label{thm:comparison-molecules-chains}
The augmented chain complex $\chains(P)$ from~\cref{d:associated_chain_complex} can be defined for any oriented graded poset $P$ which is oriented thin~\cite[Par.~11.1.6]{Hadzihasanovic24}.
Moreover, the former satisfies~\cref{ss:atoms} and~\cref{t:unital} whenever the latter is a regular directed complex \cite[Lem.~11.2.9 \& Prop.~11.2.11]{Hadzihasanovic24}.
One can then apply Steiner's~$\nu$ functor and consider the $\omega$-category $\nu(\chains(P))$.
If $P$ is dimension-wise acyclic, then the $\omega$-categories $\Mol_{/P}$ and $\nu(\chains(P))$ are isomorphic \cite[Thm.~11.2.18]{Hadzihasanovic24}.
Moreover, we have that $\nu(\chains(P))$ is loop-free (resp.\ strongly loop-free) if and only if $\Mol_{/P}$ is dimension-wise acyclic (resp.\ acyclic).
It would be interesting to compare $\Mol_{/\faces}$ and $\nu(\chains(P))$ for non loop-free framed polytopes \cite[Prop.~11.2.17~\&~11.2.35]{Hadzihasanovic24}.

%% file: fig/orientations.tex
\raisebox{.75cm}{\begin{tikzpicture}%
	[
 	scale = 1.5,
	sourceedge/.style = {color = red, ultra thick},
	targetedge/.style = {color = blue, ultra thick},
	targetvertex/.style = {inner sep = 1pt,circle,draw = blue,fill = blue,thick},
	sourcevertex/.style = {inner sep = 1pt,circle,draw = red,fill = red,thick }]

	\coordinate (A) at (-1,0);
	\coordinate (B) at (1,0);
	\coordinate (C) at (0,1.732);
	\coordinate (C1) at (-0.03,1.732);
	\coordinate (C2) at (0.03,1.732);
	
	\draw[sourceedge,-latex] (A) -- (B);
	\draw[targetedge,-latex] (A) -- (C1);
	\draw[targetedge,-latex] (C2) -- (B);
	
	\node[sourcevertex, label = left:{+}] at (A) {};
	\node[targetvertex,label = right:{+}] at (B) {};
	\node[targetvertex] at (C1) {};
	\node[sourcevertex, label = above:{+}] at (C2) {};
	
	\node at (0,0.5) {{ \scalebox{1.5}{$\circlearrowleft$}}};
	
	\begin{scope}[shift = {(-2.5,0)}, scale = .5]
	
	\end{scope}

\end{tikzpicture}}
\qquad\qquad
\begin{tikzpicture}[scale = .8]
	 \node (max) at (0,3.5)
	 {
	 \begin{tikzpicture}%
		[
	 	scale = 0.25,
		edge/.style = {thick},
		vertex/.style = {inner sep = 1pt,circle,thick }]
	\coordinate (A) at (-1,0);
	\coordinate (B) at (1,0);
	\coordinate (C) at (0,1.732);
	
	\draw[edge] (A) -- (B);
	\draw[edge] (A) -- (C);
	\draw[edge] (C) -- (B);
	
	\node[vertex] at (A) {};
	\node[vertex] at (B) {};
	\node[vertex] at (C) {};

	\fill[fill = green] (A)--(B)--(C)--(A);
	\end{tikzpicture}
	};
	 \node (a) at (-2,2) {
	 \begin{tikzpicture}%
		[
	 	scale = 0.25,
		edge/.style = {thick},
		vertex/.style = {inner sep = 1pt,circle,thick }]
	\coordinate (A) at (-1,0);
	\coordinate (B) at (1,0);
	\coordinate (C) at (0,1.732);
	
	\draw[edge, green,very thick] (A) -- (B);
	\draw[edge] (A) -- (C);
	\draw[edge] (C) -- (B);
	
	\node[vertex] at (A) {};
	\node[vertex] at (B) {};
	\node[vertex] at (C) {};

	\end{tikzpicture}
	};
	 \node (b) at (0,2) {
	 \begin{tikzpicture}%
		[
	 	scale = 0.25,
		edge/.style = {thick},
		vertex/.style = {inner sep = 1pt,circle,thick }]
	\coordinate (A) at (-1,0);
	\coordinate (B) at (1,0);
	\coordinate (C) at (0,1.732);
	
	\draw[edge] (A) -- (B);
	\draw[edge, green, very thick] (A) -- (C);
	\draw[edge] (C) -- (B);
	
	\node[vertex] at (A) {};
	\node[vertex] at (B) {};
	\node[vertex] at (C) {};
	
	\end{tikzpicture}
	};
	 \node (c) at (2,2) {
	 \begin{tikzpicture}%
		[
	 	scale = 0.25,
		edge/.style = {thick},
		vertex/.style = {inner sep = 1pt,circle,thick }]
	\coordinate (A) at (-1,0);
	\coordinate (B) at (1,0);
	\coordinate (C) at (0,1.732);
	
	\draw[edge] (A) -- (B);
	\draw[edge] (A) -- (C);
	\draw[edge, green,very thick] (C) -- (B);
	
	\node[vertex] at (A) {};
	\node[vertex] at (B) {};
	\node[vertex] at (C) {};

	\end{tikzpicture}
	};
	 \node (d) at (-2,0) {
	 \begin{tikzpicture}%
		[
	 	scale = 0.25,
		edge/.style = {thick},
		vertex/.style = {inner sep = 1pt,circle,thick }]
	\coordinate (A) at (-1,0);
	\coordinate (B) at (1,0);
	\coordinate (C) at (0,1.732);
	
	\draw[edge] (A) -- (B);
	\draw[edge] (A) -- (C);
	\draw[edge] (C) -- (B);
	
	\node[vertex, circle, fill = green,inner sep = 1pt ] at (A) {};
	\node[vertex] at (B) {};
	\node[vertex] at (C) {};
	
	\end{tikzpicture}
	 };
	 \node (e) at (0,0) {\begin{tikzpicture}%
		[
	 	scale = 0.25,
		edge/.style = {thick},
		vertex/.style = {inner sep = 1pt,circle,thick }]
	\coordinate (A) at (-1,0);
	\coordinate (B) at (1,0);
	\coordinate (C) at (0,1.732);
	
	\draw[edge] (A) -- (B);
	\draw[edge] (A) -- (C);
	\draw[edge] (C) -- (B);
	
	\node[vertex] at (A) {};
	\node[vertex, circle, fill = green,inner sep = 1pt ] at (B) {};
	\node[vertex] at (C) {};
	
	\end{tikzpicture}};
	 \node (f) at (2,0) {\begin{tikzpicture}%
		[
	 	scale = 0.25,
		edge/.style = {thick},
		vertex/.style = {inner sep = 1pt,circle,thick }]
	\coordinate (A) at (-1,0);
	\coordinate (B) at (1,0);
	\coordinate (C) at (0,1.732);
	
	\draw[edge] (A) -- (B);
	\draw[edge] (A) -- (C);
	\draw[edge] (C) -- (B);
	
	\node[vertex] at (A) {};
	\node[vertex ] at (B) {};
	\node[vertex, circle, fill = green,inner sep = 1pt] at (C) {};
	
	\end{tikzpicture}};

 \node (0) at (0,-1.5) {\begin{tikzpicture}%
		[
		scale = 0.25,
		edge/.style = {thick},
		vertex/.style = {inner sep = 1pt,circle,thick }]
		\coordinate (A) at (-1,0);
		\coordinate (B) at (1,0);
		\coordinate (C) at (0,1.732);
		
		\draw[edge] (A) -- (B);
		\draw[edge] (A) -- (C);
		\draw[edge] (C) -- (B);
		
		\node[vertex] at (A) {};
		\node[vertex ] at (B) {};
		\node[vertex] at (C) {};
		
\end{tikzpicture}};
	
	 \draw (d) --node[left]{$-$} (a) --node[left]{$-$} (max) --node[left]{$+$} (b) --node[pos = 0.25,above = .8pt]{$+$} (f)
	 (f) --node[right]{$-$} (c) --node[right]{$+$} (max)
	 (d) --node[pos = 0.7,above = .8pt]{$-$} (b);
	 \draw[preaction = {draw = white, -,line width = 6pt}] (a) --node[pos = .75,right]{$+$} (e) --node[pos = .25,left]{$+$} (c);
	 
	 \draw (d) --node[left]{$+$} (0);
 \draw (e) --node[left]{$+$} (0);
 \draw (f) --node[right]{$+$} (0);
\end{tikzpicture}

%% file: fig/thin.tex
\begin{tikzpicture}[scale = .6]
	\node (A) at (0,2) {$F$};
	\node (B) at (2,4) {$G$};
	\node (C) at (4,2) {$F'$};
	\node (D) at (2,0) {$E$};

	\draw[-] (A) -- (B) node[midway, left] {$+$};
	\draw[-] (A) -- (D) node[midway, left] {$+$};
	\draw[-] (B) -- (C) node[midway, right] {$-$};
	\draw[-] (D) -- (C) node[midway, right] {$+$};
\end{tikzpicture}
\qquad
\begin{tikzpicture}[scale = .6]
	\node (A) at (0,2) {$F$};
	\node (B) at (2,4) {$G$};
	\node (C) at (4,2) {$F'$};
	\node (D) at (2,0) {$E$};

	\draw[-] (A) -- (B) node[midway, left] {$+$};
	\draw[-] (A) -- (D) node[midway, left] {$+$};
	\draw[-] (B) -- (C) node[midway, right] {$+$};
	\draw[-] (D) -- (C) node[midway, right] {$+$};
\end{tikzpicture}